\documentclass[3pt,bookmark=true]{article}

 \RequirePackage{geometry}
\usepackage{graphicx}				

\usepackage{amssymb}
\usepackage{amsmath}
\usepackage{amsthm}

 \usepackage{lineno}

\newtheorem{theorem}{\indent Theorem}[section]

\newtheorem{definition}[theorem]{\indent Definition}
\newtheorem{lemma}[theorem]{\indent Lemma}
\newtheorem{remark}[theorem]{\indent Remark}

\begin{document}
\title{Large deviation principle for the three dimensional planetary geostrophic equations of large-scale ocean circulation with small multiplicative noise}
\author{Bo You\footnote{Email address: youb2013@xjtu.edu.cn} \\
{\small School of Mathematics and Statistics, Xi'an Jiaotong University, Xi'an, 710049, P. R. China}
}


\date{July 14, 2020}							
\maketitle

\begin{center}
\begin{abstract}
We demonstrate the large deviation principle in the small noise limit for the three dimensional stochastic planetary geostrophic equations of large-scale ocean circulation. In this paper, we first prove the well-posedness of weak solutions to this system by the method of monotonicity. As we know, a recently developed method, weak convergent method, has been employed in studying the large deviations and this method is essentially based on the main result of \cite{ba2} which discloses the variational representation of exponential integrals with respect to the Brownian noise. The It\^{o} inequality and Burkholder-Davis-Gundy inequality are the main tools in our proofs, and the weak convergence method introduced by Budhiraja, Dupuis and Ganguly in \cite{ba3} is also used to establish the large deviation principle.


\end{abstract}
\end{center}

\section{Introduction}
\def\theequation{1.\arabic{equation}}\makeatother
\setcounter{equation}{0}
This paper is concerned with the large deviations for the following three dimensional stochastic planetary geostrophic equations of the large-scale ocean circulation driven by the small multiplicative noise:
\begin{equation}\label{1.1}
\begin{cases}
&\nabla p^\epsilon+f{v^\epsilon}^\bot+L_1v^\epsilon=0,\\
&\frac{\partial p^\epsilon}{\partial z}+\theta^\epsilon=0,\\
&\nabla\cdot v^\epsilon+\frac{\partial w^\epsilon}{\partial z}=0,\\
&d\theta^\epsilon+v^\epsilon\cdot\nabla \theta^\epsilon\,dt+w^\epsilon\frac{\partial \theta^\epsilon}{\partial z}\,dt+L_2\theta^\epsilon\,dt=g(x,y,z,t)\,dt+\sqrt{\epsilon}\sigma(t,\theta^\epsilon(t))dW(t)
\end{cases}
\end{equation}
in the domain
\begin{align*}
\mathcal{O}=\mathcal{M}\times (-h,0)\subset\mathbb{R}^3.
\end{align*}
Here, $\mathcal{M}\subset\mathbb{R}^2$ is a bounded domain with smooth boundary $\partial\mathcal{M},$ $h>0$ is the depth of the ocean, $\epsilon\in(0,1)$ is a small parameter. $(v^\epsilon,w^\epsilon)=(v^\epsilon_1,v^\epsilon_2,w^\epsilon)=(v^\epsilon_1(x,y,z,t),v^\epsilon_2(x,y,z,t),w^\epsilon(x,y,z,t))$ is the velocity field, $p^\epsilon(x,y,z,t)$ is the pressure function, $\theta^\epsilon(x,y,z,t)$ is the temperature function, $g(x,y,z,t)$ is a deterministic heat source, ${v^\epsilon}^\bot=(-v_2^\epsilon,v_1^\epsilon),$ $f$ is the Coriolis parameter defined by $f=f_0+\beta y$ with the constants $ f_0$ and $\beta$, and the noise coefficient $\sigma$ satisfies some assumptions specified in the sequel. $W$ is a $U$-valued cylindrical Wiener process with respect to a complete filtered probability space $(\Omega,\mathcal{F},\{\mathcal{F}_t\}_{t\geq 0}, \mathcal{P})$ with the covariance operator $Q,$ where $U$ is a Hilbert space, $Q$ is a positive, symmetric, trace class operator on $U.$

Throughout this paper, we denote the two-dimensional horizontal gradient and Laplacian by $\nabla$ and $\Delta,$ respectively, and define $L_1$ and $L_2$ by
\begin{align*}
&L_1=-A_h\Delta-A_{\nu}\frac{\partial^2}{\partial z^2},\\
&L_2=-K_h\Delta-K_{\nu}\frac{\partial^2}{\partial z^2},
\end{align*}
where $A_h,$ $A_{\nu}$ are positive eddy viscosity constants and $K_h,$ $K_{\nu}$ are positive conductivity constants.
For the sake of simplicity, we define
\begin{align*}
&\Gamma_u:=\{(x,y,z)\in\bar{\mathcal{O}}:z=0\},\\
&\Gamma_l:=\{(x,y,z)\in\bar{\mathcal{O}}:(x,y)\in \partial\mathcal{M},-h\leq z\leq 0\},\\
&\Gamma_b=\{(x,y,z)\in\bar{\mathcal{O}}:z=-h\}.
\end{align*}
Equations \eqref{1.1} is subject to the following boundary conditions with the wind-driven on the top surface, nonslip and non-flux on the side walls and bottom (cf \cite{pj1, sm2})
\begin{equation}\label{1.2}
\begin{cases}
&\left.A_{\nu}\frac{\partial v^\epsilon}{\partial z}\right|_{\Gamma_u}=\mu, \left.w^\epsilon\right|_{\Gamma_u}=0, \left.\left(K_{\nu}\frac{\partial \theta^\epsilon}{\partial z}+\beta(\theta^\epsilon-\theta^*)\right)\right|_{\Gamma_u}=0,\\
&\left.\frac{\partial v^\epsilon}{\partial z}\right|_{\Gamma_b}=0, \left.w^\epsilon\right|_{\Gamma_b}=0, \left.\frac{\partial \theta^\epsilon}{\partial z}\right|_{\Gamma_b}=0,\\
&\left.v^\epsilon\cdot\vec{n}\right|_{\Gamma_l}=0,\left.\frac{\partial v^\epsilon}{\partial\vec{n}}\times\vec{n}\right|_{\Gamma_l}=0,\left.\frac{\partial \theta^\epsilon}{\partial\vec{n}}\right|_{\Gamma_l}=0
\end{cases}
\end{equation}
and initial data
\begin{align}
\label{1.3}\theta^\epsilon(x,y,z,0)=\theta_0(x,y,z).
\end{align}
Here $\mu(x,y)$ is the wind stress, $\vec{n}$ is the external unit normal vector on $\Gamma_l,$ $\beta>0$ is a positive constant, $\theta^*(x,y)$ is the typical temperature of the top surface satisfying the compatibility boundary condition
\begin{align*}
\left.\frac{\partial \theta^*}{\partial\vec{n}}\right|_{\partial\mathcal{M}}=0.
\end{align*}

The inviscid planetary geostrophic equations are derived from the Boussinesq equations for the planetary scale ocean by standard scale analysis like \cite{pj, pj1, pna, ra, wp}.  The well-posedness and the long-time behavior of the solutions for the three dimensional planetary geostrophic viscous equations of the large-scale ocean circulation were widely considered during the past several decades in e.g., \cite{ccs, sm, sm1, yb3, yb1, yb2, yb}. In particular, it was shown in \cite{yb2} that the existence of a random attractor in $H$ with small multiplicative noise was proved by verifying the pullback flattening property.  The similar results were later obtained in \cite{yb3} for small additive noise case by Sobolev compactness embedding theorem and the pullback flattening property. Recently, the authors in \cite{dz} established the connection between the invariant measure of the corresponding Markovian semigroup and the random attractor. To the best of our knowledge, however, the large deviations for it with any kind of processes were few results.

The large deviation principle arises in the theory of statistical inference quite naturally and offers a precise estimation associated with the law of large number. Moreover, the moderate deviation provides us with the rate of convergence, which is a useful method for constructing asymptotic confidence intervals by giving further estimations related to center limit theorem and the law of iterated logarithm. We refer the reader to \cite{em, gfq, it, kw}. There mainly exist three typical approaches to analyze large deviation principle for stochastic partial differential equations under small perturbations in the literature since the first research was formulated by Varadhan \cite{vsrs}.  Freidlin and Wentzell \cite{ff1} developed one way to deal with some semi-linear stochastic partial differential equations in finite dimensional case, based on discretization approximations and the contraction principle. Substantial progress has been made for this approach since then; see some papers like \cite{cwc, cs, cpl, ps, pgd, rm1, xt, zts} in both finite and infinite dimensional cases. But the situation became much complicated in the infinite dimensional case since there is no uniform way to deal with nonlinear stochastic partial differential equations. The second one was developed in \cite{fj} by using nonlinear semigroup theory and infinite-dimensional Hamilton-Jacobi equations. However, this approach relies on the uniqueness theory for the infinite dimensional Hamilton-Jacobi equation and some exponential tightness estimates. The last one is the so-called weak convergence method which is originally established by Budhiraja, Dupuis and Ganguly \cite{ba3}. The main idea is based on a variational representation for certain functionals of infinite dimensional Brownian motion, whose advantage is to avoid some exponential probability estimates that might be very difficult to be derived for infinite dimensional models. Hence, it was extensively used to investigate the large deviation principle. For further researches on this approach we may refer to \cite{ba3, ci, djq, dp, rj, rj1, rm2, sss, zjl} for detail discussions. It should be noticed that some technical difficulties need to be overcome in the variational framework in the implement of weak convergence approach.

 Inspired by the above works, we intend to take advantage of in this paper the weak convergence approach to study the large deviation principle of the problem \eqref{1.1}-\eqref{1.3}. Precisely, let $\theta^\epsilon$ and $\theta$ be the solutions for the problem \eqref{1.1}-\eqref{1.3} with $\epsilon\in (0,1)$ and  $\epsilon=0$, respectively. The large deviation principle deals with deviations of the asymptotic behavior of the trajectory $\theta^\epsilon(t)$ as $\epsilon \rightarrow 0.$ To establish this result, we need to verify the condition (A) of Lemma \ref{2.3.1}. There exists mainly two difficulties: On one hand, due to $\sigma(s,\theta(s))$ is a Hilbert-Schmidt operator from $U_0$ to $H_2$ for any fixed $s\in[0,T],$ we can only obtain $\sigma(s,\theta(s))\varphi^n\rightarrow\sigma(s,\theta(s))\varphi$ in $H_2$ for any fixed $s\in[0,T],$ if $\varphi^n\rightharpoonup \varphi$ in $U_0.$ Therefore, it is very difficult to  prove that the integral $\sup_{t\in[0,T]}\left|\int_0^t\int_{\mathcal{O}}\sigma(s,\theta(s))(\chi_n(s)-\chi(s))\eta_n(s)\,dxdydzds\right|$ tends to zero as $n\rightarrow+\infty.$ Inspired by the idea in \cite{ci}, we will in this paper carry out some a priori estimates and establish a technical result (i.e., Theorem \ref{5.2}) to conclude the desired conclusion. On the other hand, as for the second condition of condition (A), there is no any information but the assumption that $\chi^\epsilon$ converges to $\chi$ in distribution as $\epsilon\rightarrow 0.$ To verify this condition, we will combining Theorem \ref{5.2} with the Skorokhod representation Theorem as well as the similar proof of the first condition of condition (A), we can complete the verification of the second condition of condition (A).

This paper is organized as follows. In Section 2, we give the preliminaries including the variational formulation of the problem \eqref{1.1}-\eqref{1.3}, some useful lemma, as well as some stand definitions and results from large deviation principle. In Section 3, we establish the well-posedness of weak solutions, and Section 4 is devoted to the large deviation principle.

\textsl{Notation.} Denote by $X$ a Banach space with the norm $\|\cdot\|_X$ and let $C$ be a positive generic constant which can change from one line to the next. If it is essential, we will write the dependence of the constant on parameters explicitly.
 \section{Preliminaries}
\def\theequation{2.\arabic{equation}}\makeatother
\setcounter{equation}{0}
\subsection{New formulation}
We can reformulate problem \eqref{1.1}-\eqref{1.3} by integrating the second and the third equation of \eqref{1.1} with respect to $z$ and combining the boundary conditions \eqref{1.2} as follows just like in \cite{ccs}:
\begin{equation}\label{2.1.1}
\begin{cases}
&\nabla p^\epsilon_s(x,y,t)-\int_{-h}^z\nabla \theta^\epsilon(x,y,\zeta,t)\,d\zeta+f{v^\epsilon}^\bot+L_1v^\epsilon=0,\\
&\int_{-h}^0\nabla\cdot v^\epsilon(x,y,\zeta,t)\,d\zeta=0,\\
&d\theta^\epsilon+v^\epsilon\cdot\nabla\theta^\epsilon\,dt-\left(\int_{-h}^z\nabla\cdot v^\epsilon(x,y,\zeta,t)\,d\zeta\right)\frac{\partial \theta^\epsilon}{\partial z}\,dt+L_2\theta^\epsilon\,dt\\
&=g(x,y,z,t)\,dt+\sqrt{\epsilon}\sigma(t,\theta^\epsilon(t))dW(t),\\
&\left.A_{\nu}\frac{\partial v^\epsilon}{\partial z}\right|_{\Gamma_u}=\mu, \left.\frac{\partial v^\epsilon}{\partial z}\right|_{\Gamma_b}=0, \left.v^\epsilon\cdot\vec{n}\right|_{\Gamma_l}=0,\left.\frac{\partial v^\epsilon}{\partial\vec{n}}\times\vec{n}\right|_{\Gamma_l}=0,\\
&\left.\left(K_{\nu}\frac{\partial \theta^\epsilon}{\partial z}+\beta(\theta^\epsilon-\theta^*)\right)\right|_{\Gamma_u}=0, \left.\frac{\partial \theta^\epsilon}{\partial z}\right|_{\Gamma_b}=0,\left.\frac{\partial \theta^\epsilon}{\partial\vec{n}}\right|_{\Gamma_l}=0,\\
&\theta^\epsilon(x,y,z,0)=\theta_0(x,y,z).
\end{cases}
\end{equation}
In order to further recast problem \eqref{2.1.1} into an abstract form, we need to introduce some notations of function space and operators. \\
Define
\begin{align*}
\mathcal {V}_1=&\left\{v\in (C^{\infty}(\bar{\mathcal{O}}))^{2}:\frac{\partial
v}{\partial z}|_{\Gamma_u}=0,\frac{\partial v}{\partial
z}|_{\Gamma_b}=0,v\cdot \vec{n}|_{\Gamma_l}=0,\frac{\partial
v}{\partial \vec{n}}\times\vec{
n}|_{\Gamma_l}=0,\int_{-h}^0\nabla\cdot v(x,y,\zeta)\,d\zeta=0\right\},\\
 \mathcal{V}_2=&\left\{ \theta\in C^\infty(\bar{\mathcal{O}}):\left.(K_{\nu}\frac{\partial \theta}{\partial z}+\beta
\theta)\right|_{\Gamma_u}=0,\left.\frac{\partial \theta}{\partial
z}\right|_{\Gamma_b}=0,\left.\frac{\partial \theta}{\partial
\vec{n}}\right|_{\Gamma_l}=0\right\}.
\end{align*}
For any $v\in\mathcal{V}_1,$ $\theta\in\mathcal{V}_2,$ denote by $H_1,$ $H_2$ the closure of $\mathcal{V}_1,$ $\mathcal{V}_2,$ respectively, with respect to the following norms
\begin{align*}
\|v\|_{L^2(\mathcal{O})}=&\left(\int_{\mathcal{O}}|v(x,y,z)|^2\,dxdydz\right)^{\frac{1}{2}},\\
\|\theta\|_{L^2(\mathcal{O})}=&\left(\int_{\mathcal{O}}|\theta(x,y,z)|^2\,dxdydz\right)^{\frac{1}{2}},
\end{align*}
and $V_1,$ $V_2$ the closure of $\mathcal{V}_1,$ $\mathcal{V}_2,$ respectively, with respect to the following norms
\begin{align*}
\|v\|_{H^1(\mathcal{O})}=&\left(A_h\int_{\mathcal{O}}|\nabla
v(x,y,z)|^2\,dxdydz+A_{\nu}\int_\mathcal{O}|v_z(x,y,z)|^2\,dxdydz\right)^{\frac{1}{2}},\\
\|\theta\|=&\left(K_h\int_{\mathcal{O}}|\nabla
\theta(x,y,z)|^2\,dxdydz+K_{\nu}\int_\mathcal{O} |\theta_z(x,y,z)|^2\,dxdydz+\beta\int_{\mathcal{M}} |\theta(x,y,0)|^2\,dxdy\right)^{\frac{1}{2}}.
\end{align*}
By Riesz isomorphism, we have $V_i\subset H_i=H_i'\subset V_i',$ where $H_i'$ is identified with $H_i$ for $i=1,2$ and $V_i'$ is the dual space of  $V_i$ with the dual action $\langle\cdot,\cdot\rangle.$\\
Define the operator $A_1: V_1\rightarrow V_1'$ associated with the bilinear form given by
\begin{align*}
\langle A_1v_1, v_2\rangle=A_h\int_{\mathcal{O}}\nabla
v_1(x,y,z)\cdot\nabla v_2(x,y,z)\,dxdydz+A_{\nu}\int_\mathcal{O}\partial_z v_1(x,y,z)\cdot\partial_z v_2(x,y,z)\,dxdydz
\end{align*}
for any $v_1,$ $v_2\in V_1,$
and the operator $A_2: V_2\rightarrow V_2'$ associated with the bilinear form given by
\begin{align*}
\langle A_2\theta_1, \theta_2\rangle=&K_h\int_{\mathcal{O}}\nabla
\theta_1(x,y,z)\cdot\nabla \theta_2(x,y,z)\,dxdydz+K_{\nu}\int_\mathcal{O}\partial_z \theta_1(x,y,z)\partial_z \theta_2(x,y,z)\,dxdydz\\
&+\beta\int_{\mathcal{M}} \theta_1(x,y,0)\theta_2(x,y,0)\,dxdy
\end{align*}
for any $\theta_1,$ $\theta_2\in V_2.$\\
Introduce a trilinear form $b$ on $\left(H^2(\mathcal{O})\cap V_1\right)\times V_2\times V_2$ by
\begin{align*}
b(v,\theta,\eta)=\int_{\mathcal{O}}\left(v(x,y,z)\cdot\nabla\theta(x,y,z)-\left(\int_{-h}^z\nabla\cdot v(x,y,\zeta)\,d\zeta\right)\partial_z \theta(x,y,z)\right)\eta(x,y,z)\,dxdydz
\end{align*}
and a bilinear form $B: \left(H^2(\mathcal{O})\cap V_1\right)\times V_2\rightarrow V_2'$ by
\begin{align*}
\langle B(v,\theta),\eta\rangle=b(v,\theta,\eta).
\end{align*}
Integrating by parts, we obtain
\begin{align*}
b(v,\theta,\eta)=-b(v,\eta,\theta)
\end{align*}
for any $v\in H^2(\mathcal{O})\cap V_1,$ $\theta,$ $\eta\in V_2.$\\
Hence, we can formally recast problem \eqref{2.1.1} into the following abstract equation
\begin{equation}\label{2.1.2}
\begin{cases}
&\nabla p^\epsilon_s+f{v^\epsilon}^\bot+L_1v^\epsilon=\int_{-h}^z\nabla \theta^\epsilon(x,y,\zeta,t)\,d\zeta,\\
&\int_{-h}^{0}\nabla\cdot v^\epsilon(x,y,\zeta,t)\,d\zeta=0,\\
&d\theta^\epsilon+B(v^\epsilon, \theta^\epsilon)\,dt+A_2(\theta^\epsilon-\theta^*)\,dt-K_h\Delta \theta^*\,dt=g(x,y,z,t)\,dt+\sqrt{\epsilon}\sigma(t,\theta^\epsilon(t))dW(t),\\
&\left.A_{\nu}\frac{\partial v^\epsilon}{\partial z}\right|_{\Gamma_u}=\mu, \left.\frac{\partial v^\epsilon}{\partial z}\right|_{\Gamma_b}=0, \left.v^\epsilon\cdot\vec{n}\right|_{\Gamma_l}=0,\left.\frac{\partial v^\epsilon}{\partial\vec{n}}\times\vec{n}\right|_{\Gamma_l}=0,\\
&\theta^\epsilon(x,y,z,0)=\theta_0(x,y,z).
\end{cases}
\end{equation}

\begin{remark}
The well-posedness of solutions $(v,p_s,\theta)\in \mathcal{C}([0,T]; V_1\times L^2(\mathcal{M})\times H_2)\cap L^2(0,T; H^2(\mathcal{O})\times H^1(\mathcal{M})\times V_2)$ for problem \eqref{2.1.2} with $\epsilon =0$ has been established in \cite{ccs} and the well-posedness of solutions $(v^\epsilon,p_s^\epsilon,\theta^\epsilon)\in L^2(\Omega,\mathcal{C}([0,T]; V_1\times L^2(\mathcal{M})\times H_2))\cap L^2(\Omega,L^2(0,T; H^2(\mathcal{O})\times H^1(\mathcal{M})\times V_2))$ for problem \eqref{2.1.2} with $0<\epsilon<1$ will be established in this paper under the following assumptions $(A_1)$-$(A_3)$ stated in subsection 2.3.
\end{remark}

\subsection{Some useful Lemmas}

In this subsection, we recall and derive some lemmas which will be required in the rest of this paper.

\begin{lemma}\label{2.2.1}(\cite{ccs})
There exists a positive constant $K_1$ such that
\begin{align*}
\frac{1}{K_1}\|\theta\|^2\leq\|\theta\|_{H^1(\mathcal{O})}^2\leq K_1\|\theta\|^2
\end{align*}
for any $\theta\in V_2.$ Moreover, we have
\begin{align*}
K_2\|\theta\|^2_{L^2(\mathcal{O})}\leq\|\theta\|^2
\end{align*}
for any $\theta\in V_2,$ where
\begin{align*}
K_2=\min\{\frac{\beta}{2h},\frac{K_{\nu}}{2h^2}\}.
\end{align*}
\end{lemma}

\begin{lemma}\label{2.2.2}(\cite{ccs})
Assume that $\mu\in H_0^1(\mathcal{M})$ and $\theta\in H^\gamma(\mathcal{O})$ with $\gamma=0$ or $1.$ Then there exists a unique solution $(v,p_s)\in H^{\gamma+1}(\mathcal{O})\times H^\gamma(\mathcal{M})$ ($p_s$ is unique up to a constant) to the following problem:
\begin{equation*}
\begin{cases}
&\nabla p_s(x,y,t)-\int_{-h}^z\nabla \theta(x,y,\zeta,t)\,d\zeta+fv^\bot+L_1v=0,\\
&\int_{-h}^0\nabla\cdot v(x,y,\zeta,t)\,d\zeta=0,\\
&\left.A_{\nu}\frac{\partial v}{\partial z}\right|_{\Gamma_u}=\mu, \left.\frac{\partial v}{\partial z}\right|_{\Gamma_b}=0, \left.v\cdot\vec{n}\right|_{\Gamma_l}=0,\left.\frac{\partial v}{\partial\vec{n}}\times\vec{n}\right|_{\Gamma_l}=0.
\end{cases}
\end{equation*}
Moreover, there exists a positive constant $\mathcal{K}_1=\mathcal{K}_1(A_\nu, A_h)$ such that
\begin{align*}
\|v\|_{H^{\gamma+1}(\mathcal{O})}^2+\|p_s\|_{H^\gamma(\mathcal{M})}^2\leq \mathcal{K}_1\left(\|\theta\|_{H^\gamma(\mathcal{O})}^2+\|\mu\|_{H^1(\mathcal{M})}^2\right).
\end{align*}
\end{lemma}
In the following, we give the estimates of the trilinear form $b$ and a new version of Gronwall inequality which will be critical in the proof that follows. For the sake of brevity, we only list the results as follows whose proofs will be stated in Appendix for the readers' convenience.
\begin{lemma}\label{2.2.3}
Assume that $v\in H^2(\mathcal{O})\cap V_1$ and $\theta,$ $\eta\in V_2.$ Then there exists a positive constant $\mathcal{K}$ such that
\begin{align*}
|b(v,\theta,\eta)|\leq \mathcal{K}\|v\|_{H^1(\mathcal{O})}^{\frac{1}{2}}\|v\|_{H^2(\mathcal{O})}^{\frac{1}{2}}\|\theta\|\|\eta\|_{L^2(\mathcal{O})}^{\frac{1}{2}}\|\eta\|^{\frac{1}{2}}.
\end{align*}
\end{lemma}
\begin{lemma}\label{2.2.4}
Let $Y(t)\in\mathcal{C}([t_0,t_1])$ be a non-negative function satisfying the following inequality:
\begin{align*}
Y(t)+\int_0^t X(s)\,ds\leq \int_0^t a(s)Y(s)\,ds+Z(t),
\end{align*}
where $a(t),$ $X(t)\in\mathcal{C}([t_0,t_1])$ are non-negative functions and $Z(t)\in\mathcal{C}([t_0,t_1])$ is non-negative, non-decreasing function. Then
\begin{align*}
Y(t)+\int_0^t X(s)\,ds\leq Z(t)\,e^{\int_0^ta(r)\,dr}.
\end{align*}
\end{lemma}
\subsection{Large deviation principle}
In this subsection, let us recall some standard definitions and results from the large deviation theory given like in \cite{ba}.

Let $X$ be a Polish space with the Borel $\sigma$-field $\mathcal{B}(X)$ and $(\Omega, \mathcal{F}, \{\mathcal{F}_t\}_{t\geq 0}, \mathcal{P})$  be a probability space with an increasing family $\{\mathcal{F}_t\}_{0\leq t\leq T}$ of the sub-$\sigma$-fields of $\mathcal{F}$ satisfying the usual conditions, that is, $\{\mathcal{F}_t\}_{t\in\mathbb{R}}$ is an increasing right continuous family of sub-$\sigma$-algebras of $\mathcal{F}$ that contains all $\mathcal{P}$-null sets. Let $W$ be a $U$-valued cylindrical Wiener process with respect to a complete filtered probability space $(\Omega,\mathcal{F},\{\mathcal{F}_t\}_{t\geq 0}, \mathcal{P})$ with the covariance operator $Q,$ where $U$ is a Hilbert space and $Q$ is a positive, symmetric, trace class operator on $U.$ Denote by $U_0=Q^{\frac{1}{2}}U$. Then $U_0$ is a Hilbert space with the inner product
\begin{align*}
(u,v)_{U_0}=(Q^{-\frac{1}{2}}u, Q^{-\frac{1}{2}}v)_{U}, \,\,\,\,\forall\,\,\,u,\,\,\,v\in U_0,
\end{align*}
where $(\cdot,\cdot)_U$ represents the inner product  in $U.$ Denote by $\mathcal{L}_2(U_0,H_2)$ the space of linear operators $S$ satisfying that $SQ^{\frac{1}{2}}$ is a Hilbert-Schmidt operator from $U$ to $H_2$ with norm $\|S\|_{\mathcal{L}_2(U_0,H_2)}=\left(tr(SQS^*)\right)^{\frac{1}{2}}$ for any $S\in \mathcal{L}_2(U_0,H_2).$

Define the Cameron-Martin space associated with the Wiener process $\{W(t):t\in[0,T]\}$ by
\begin{align*}
\mathcal{H}_0=\left\{\chi:[0,T]\rightarrow U_0:\chi\,\,\textit{is\,\,absolutely\,\,continuous\,\,and\,\,}\,\,\int_0^T\|\dot{\chi}(s)\|_{U_0}^2\,ds<+\infty\right\}.
\end{align*}
Hence, the space $\mathcal{H}_0$ is a Hilbert space with inner product
\begin{align*}
(\chi_1,\chi_2)_{\mathcal{H}_0}=\int_0^T(\dot{\chi}_1(s),\dot{\chi}_2(s))_{U_0}\,ds.
\end{align*}
Let $\mathcal{A}$ be the class of $U_0$-valued $\{\mathcal{F}_t\}$-predictable processes $\chi$ belonging to $\mathcal{H}_0$ a.s. Define
\begin{align*}
S_N=\left\{\chi\in\mathcal{H}_0:\int_0^T\|\dot{\chi}(s)\|_{U_0}^2\,ds\leq N\right\},
\end{align*}
then the set $S_N$ endowed with the weak topology is a Polish space. Denote by 
\begin{align*}
\mathcal{A}_N=\left\{\chi\in\mathcal{A}:\chi(\omega)\in S_N,\mathcal{P}-a.s\right\}.
\end{align*}

Roughly speaking, the large deviation theory concerns itself with the exponential decay of the probability measures of certain kinds of extreme or tail events. The rate of such exponential decay is expressed by the "rate function".
\begin{definition}
(Rate function). A function $I:X\rightarrow[0,+\infty]$ is called a good rate function on $X,$ if for each $M<+\infty,$ the level set $\{x\in X:I(x)\leq M\}$ is a compact subset of $X.$
\end{definition}

\begin{definition}
(Large deviation principle). Let $I$ be a good rate function on $X.$ A family $\{X^\epsilon\}$ of $X$-valued random elements is said to satisfy the large deviation principle on $X$ with rate function $I,$ if the following two conditions hold:
\begin{itemize}
\item [(i)] (Large deviation upper bound) For each closed subset $F$ of $X,$
\begin{align*}
\limsup_{\epsilon\rightarrow 0}\epsilon\log\mathcal{P}\left(X^\epsilon\in F\right)\leq -\inf_{x\in F}I(x).
\end{align*}
\item [(ii)] (Large deviation lower bound) For each open subset $G$ of $X,$
\begin{align*}
\liminf_{\epsilon\rightarrow 0}\epsilon\log\mathcal{P}\left(X^\epsilon\in G\right)\geq -\inf_{x\in G}I(x).
\end{align*}
\end{itemize}
\end{definition}

%

Now we state the following sufficient condition for large deviation principle given by Budhiraja and Dupuis in \cite{ba}.
\begin{itemize}
\item [(A)]  There exists a measurable mapping $\Phi^0:\mathcal{C}([0,T]; U)\rightarrow X$ such that the following two conditions hold:
\begin{itemize}
\item [(i)] for every $N<+\infty,$ the set
\begin{align*}
\left\{\Phi^0\left(\int_0^\cdot\dot{\chi}(s)\,ds\right):\chi\in S_N\right\}
\end{align*}
is a compact subset of $X.$
\item [(ii)] let $\{\chi^\epsilon:\epsilon>0\}\subset\mathcal{A}_N$ for some $N<+\infty$. If  $\chi^\epsilon$ converge to $\chi$ in distribution as $S_N$-valued random elements, then
\begin{align*}
\Phi^\epsilon(W(\cdot)+\frac{1}{\sqrt{\epsilon}}\int_0^\cdot\dot{\chi}^\epsilon(s)\,ds)\rightarrow \Phi^0(\int_0^\cdot\dot{\chi}(s)\,ds)
\end{align*}
in distribution as $\epsilon\rightarrow 0.$
\end{itemize}
\end{itemize}
\begin{lemma}(\cite{ba})\label{2.3.1}
For $\epsilon>0,$ let $\Phi^\epsilon$ be a measurable mapping from $\mathcal{C}([0,T];U)$ into $X$ and $X^\epsilon=\Phi^\epsilon(W(\cdot)).$ If $\{\Phi^\epsilon:\epsilon>0\}$ satisfies the assumption (A),
then the family $\{X^\epsilon:\epsilon>0\}$ satisfies a large deviation principle in $X$ with the rate function $I$ given by
\begin{align}\label{2.3.2}
I(f)=\inf_{\left\{\chi\in\mathcal{H}_0:f=\Phi^0\left(\int_0^\cdot\dot{\chi}(s)\,ds\right)\right\}}\left(\frac{1}{2}\int_0^T\|\dot{\chi}(s)\|_{U_0}^2\,ds\right),\,\,\,\forall\,\,f\in X,
\end{align}
with the convention that the infimum of an empty set is infinity.
\end{lemma}

  We conclude this section by giving some basic assumptions used in this paper. Assume that $\mu\in H_0^1(\mathcal{M}),$ $\theta^*\in H^2(\mathcal{M})$, $g\in L_{loc}^2(\mathbb{R}; L^2(\mathcal{O}))$ and the following assumption holds:
\begin{itemize}
\item [(B)] For any $T>0,$ the diffusion coefficient $\sigma:[0,T]\times H_2\rightarrow \mathcal{L}_2(U_0,H_2)$ is progressively measurable and satisfies the following conditions:
\begin{itemize}
\item [($A_1$)]  $\sigma\in\mathcal{C}([0,T]\times H_2;\mathcal{L}_2(U_0,H_2)).$
\item [($A_2$)] There exists a positive constant $K$ such that for all  $t\in(0,T),$ $\theta\in H_2,$
\begin{align*}
\|\sigma(t,\theta)\|_{\mathcal{L}_2(U_0,H_2)}^2\leq K(1+\|\theta\|_{L^2(\mathcal{O})}^2).
\end{align*}
\item [($A_3$)] There exists a positive constant $L$ such that for all $t\in(0,T),$ $\theta_1,$ $\theta_2\in H_2,$
\begin{align*}
 \|\sigma(t,\theta_1)-\sigma(t,\theta_2)\|_{\mathcal{L}_2(U_0,H_2)}^2\leq L\|\theta_1-\theta_2\|_{L^2(\mathcal{O})}^2.
 \end{align*}
 \item [($A_4$)] (Time H\"{o}lder regularity of $\sigma$) There exist two positive constants $\gamma>0$ and $L_1\geq 0$ such that for all $t_1,$ $t_2\in [0,T]$ and $\theta\in H_2,$
\begin{align*}
 \|\sigma(t_1,\theta)-\sigma(t_2,\theta)\|_{\mathcal{L}_2(U_0,H_2)}\leq L_1(1+\|\theta\|_{L^2(\mathcal{O})})|t_1-t_2|^\gamma. 
 \end{align*} 
\end{itemize}
\end{itemize}
\begin{definition}\label{2.3.3}
An $V_1\times L^2(\mathcal{M})\times H_2$-valued c\`{a}dl\`{a}g $\mathcal{F}_t$-measurable process $(v^\epsilon(x,y,z,t), p^\epsilon_s(x,y,t), \theta^\epsilon(x,y,z,t))$ is said to be a weak solution of problem \eqref{2.1.2} on $[0,T]$ for any $T>0,$ if the following conditions are satisfied
\begin{itemize}
\item [(i)]
\begin{equation*}
\begin{cases}
&p^\epsilon_s(x,y,t)\in L^2(\Omega,C([0,T];L^2(\mathcal{M})))\cap L^2(\Omega,L^2(0,T;H^1(\mathcal{M}))),\\
&v^\epsilon(x,y,z,t)\in L^2(\Omega,C([0,T];V_1))\cap L^2(\Omega,L^2(0,T;H^2(\mathcal{O}))),\\
&\theta^\epsilon(x,y,z,t)\in L^2(\Omega,C([0,T]; L^2(\mathcal{O})))\cap L^2(\Omega,L^2(0,T;V_2));
\end{cases}
\end{equation*}
\item [(ii)] For any $\phi\in (H^1(\mathcal{O}))^2,$
\begin{align}\label{2.3.4}
\nonumber&\int_\mathcal{O} \nabla p^\epsilon_s(x,y,t)\cdot\phi\,dxdydz-\int_{\mathcal{O}}\left(\int_{-h}^z
\nabla \theta^\epsilon(x,y,\zeta,t)d\zeta\right)\cdot\phi\,dxdydz+\int_{\mathcal{O}}f{v^\epsilon}^\bot\cdot\phi\,dxdydz+\langle A_1v^\epsilon, \phi\rangle\\
&=\int_{\Gamma_u} \kappa\mu\cdot\phi\, dxdy,\,\,\,\,d\mathcal{P}\otimes dt-a.e.\,\,\,\textit{on}\,\,\,\Omega\times(0,T),
\end{align}
and for any $t\in [0,T]$ and $\mathcal{F}_0$-measurable $H_2$-valued initial data $\theta_0,$ the following equality holds $\mathcal{P}$-a.s.
\begin{align}\label{2.3.5}
\nonumber\theta^\epsilon(t)=&\theta_0-\int_0^t A_2(\theta^\epsilon(s)-\theta^*)\,ds-\int_0^t B(v^\epsilon(s),\theta^\epsilon(s))\,ds
+K_ht\Delta \theta^* +\int_0^t g(s)
\,ds\\
&+\sqrt{\epsilon}\int_0^t\sigma(s,\theta^\epsilon(s))\,dW(s).
\end{align}
\end{itemize}
\end{definition}
\section{The well-posedness of solutions}
\def\theequation{3.\arabic{equation}}\makeatother
\setcounter{equation}{0}
In this section, we will prove the existence and uniqueness of weak solutions for the three dimensional stochastic planetary geostrophic equations of large-scale ocean circulation \eqref{2.1.2} by the method of monotonicity as in \cite{cpl, mm}.
\begin{theorem}\label{3.1}
Suppose that assumptions $(A_1)$-$(A_3)$ hold, $\theta_0\in L^2(\Omega, H_2)$ is $\mathcal{F}_0$-measurable and $g\in L^2(0,T; L^2(\mathcal{O}))$ for every $T>0.$ Then problem \eqref{2.1.2} has a unique weak solution $(v^\epsilon, p^\epsilon_s,\theta^\epsilon)$ in the sense of Definition \ref{2.3.3}. Moreover, it satisfies the following estimates:
\begin{align*}
&\|(v^\epsilon, p^\epsilon_s, \theta^\epsilon)\|_{L^2(\Omega,\mathcal{C}(0,T; V_1\times L^2(\mathcal{M})\times H_2))}+\|(v^\epsilon, p^\epsilon_s, \theta^\epsilon)\|_{L^2(\Omega,L^2(0,T; H^2(\mathcal{O})\times H^1(\mathcal{M})\times V_2))}\\
\leq &K_1(T)\left(1+\|\theta_0\|_{L^2(\Omega, H_2)}^2+\|g\|_{L^2(0,S: H_2)}^2\right),
\end{align*}
where $K_1(T)$ is a positive constant depending only on $T.$ 
\end{theorem}
\textbf{Proof.} In what follows, we will first prove the existence of weak solutions for problem \eqref{2.1.2} by using Galerkin approximation methods (see \cite{tr}) and the method of monotonicity as in \cite{cpl, mm, pe, sss}. We will do this in two steps.

{\bf Step 1.} Assume that $\theta_0\in L^4(\Omega, H_2).$

It is well-known (see \cite{ccs}) that for the eigenvalue problem $A_2\omega=\lambda \omega,$ there exists a sequences of non-decreasing numbers $\{\lambda_n\}_{n=1}^{\infty}$ and a sequences of functions $\{\omega_n\}_{n=1}^{\infty}$ such that for every $k\geq 1,$ we have
\begin{align*}
A_2\omega_k=\lambda_k \omega_k
\end{align*}
 and
\begin{align*}
\lim_{k\rightarrow+\infty}\lambda_k=+\infty.
\end{align*}
Moreover, the eigenfunctions form an orthonormal basis of $H_2,$ which are also orthogonal basis of $V_2.$ For any $n\geq1,$ we introduce a finite-dimensional space $H_n=span\{\omega_1, ...,\omega_n\}.$ Let $P_n$ be
the orthogonal projector from $L^2(\Omega)$ to $H_n.$

We are looking for an approximate solution $\theta^\epsilon_n(t)$ having the form
\begin{align*}
\theta^\epsilon_n(t)-\theta^*=\sum_{i=1}^n \alpha_i(t)\omega_i.
\end{align*}
Such an approximate solution satisfies the problem
\begin{equation}\label{3.2}
\begin{cases}
&\int_\mathcal{O} \nabla p^\epsilon_{sn}(x,y,t)\cdot\phi\,dxdydz-\int_{\mathcal{O}}\left(\int_{-h}^z
\nabla \theta^\epsilon_n(x,y,\zeta,t)d\zeta\right)\cdot\phi\,dxdydz+\langle A_1v^\epsilon_n, \phi\rangle\\
&=\int_{\Gamma_u} \kappa\mu\cdot\phi\, dxdy-\int_{\mathcal{O}}f{v^\epsilon_n}^\bot\cdot\phi\,dxdydz,\\
&\int_{-h}^0\nabla\cdot v^\epsilon_n(x,y,\zeta,t)\,d\zeta=0,\\
&\int_{\mathcal{O}}\theta^\epsilon_n(t)\psi
\,dxdydz+\int_0^t\langle A_2(\theta^\epsilon_n-\theta^*), \psi\rangle\,ds-K_ht\int_\mathcal{O}\Delta \theta^* \psi\,dxdydz+\int_0^tb(v^\epsilon_n,\theta^\epsilon_n,\psi)\, ds\\
&=\int_{\mathcal{O}}\theta_{0n}\psi\, dxdydz+\int_0^t\int_{\mathcal{O}}g(s)\psi
\,dxdydzds+\sqrt{\epsilon}\int_0^t\int_{\mathcal{O}}\psi\sigma_n(s,\theta_n^\epsilon(s))\,dW(s)dxdydz,\\
&\theta_n(0)=P_n\theta_0=\sum_{i=1}^n\left(\int_\mathcal{O} \theta_0(x,y,z)\omega_i(x,y,z)\,dxdydz\right)\omega_i(x,y,z)
\end{cases}
\end{equation}
for any $\phi\in H^1(\mathcal{O})$ and $\psi\in H_n,$ where $\sigma_n=P_n\sigma.$

From Lemma \ref{2.2.2}, we deduce that for any fixed and given $\theta^\epsilon_n,$ there is a unique $(v^\epsilon_n, p^\epsilon_{sn})=(v^\epsilon(\theta_n^\epsilon), p^\epsilon_s(\theta^\epsilon_n))$ ($p^\epsilon_{sn}$ is unique up to a constant) such that for $\gamma=0$ or $1,$
\begin{align*}
v^\epsilon_n\in H^{\gamma+1}(\mathcal{O}),\,\,\,p^\epsilon_{sn}\in H^\gamma(\mathcal{M}).
\end{align*}
Furthermore, we have
\begin{align}\label{3.3}
\|v^\epsilon_n(t)\|_{H^{\gamma+1}(\mathcal{O})}^2+\|p^\epsilon_{sn}\|_{H^\gamma(\mathcal{M})}^2\leq \mathcal{K}_1\left(\|\theta^\epsilon_n\|_{H^\gamma(\mathcal{O})}^2+\|\mu\|_{H^1(\mathcal{M})}^2\right).
\end{align}
Let $v_n=v(\theta^\epsilon_n)$ in the third equation of \eqref{3.2}, we get an ordinary differential equations of the unknown $\theta^\epsilon_n,$ that is, the third equation of \eqref{3.2} is an ordinary differential equations with the unknown $\alpha_i(t), \,\,i=1,2,\cdots, n.$ It is easy to check that each term of the third equation of \eqref{3.2} is locally Lipschitz in $\theta^\epsilon_n.$ Therefore, from the theory of stochastic differential equations (see, for instance, the existence results given in \cite{sav} ), there exists a local solution $\theta^\epsilon_n$  to the equation the third equation of \eqref{3.2} defined on an interval $[0, T_n],$ which implies that there exists a unique solution $v^\epsilon_n(x,y,z,t)$ and $p^\epsilon_{sn}(x,y,t)$ of the first equation of \eqref{3.2} on $[0, T_n].$ From the estimates below, we will conclude that $T_n = T.$

It follows from It\^{o}'s Lemma that
\begin{align*}
&d\|\theta^\epsilon_n(t)\|_{L^2(\mathcal{O})}^2+2\|\theta^\epsilon_n(t)\|^2\,dt-2\beta\int_{\Gamma_u}\theta^*\theta^\epsilon_n(t)\, dxdydt\\
=&2\int_{\mathcal{O}}g(t)\theta^\epsilon_n(t)
\,dxdydzdt+2\sqrt{\epsilon}\int_{\mathcal{O}}\theta^\epsilon_n(t)\sigma_n(t,\theta^\epsilon_n(t))\,dW(t)dxdydz+\epsilon\|\sigma_n(t,\theta^\epsilon_n(t))\|_{\mathcal{L}_2(U_0,H_2)}^2\,dt.
\end{align*}
Define $\tau^n_N=\inf\{t>0:\|\theta^\epsilon_n(t)\|_{L^2(\mathcal{O})}^2+\int_0^t\|\theta^\epsilon_n(s)\|^2\,ds>N\},$ it follows from Young's inequality, H\"{o}lder's inequality, inequality \eqref{3.3} and Lemma \ref{2.2.1} that 
\begin{align*}
&\|\theta^\epsilon_n(t\land\tau^n_N)\|_{L^2(\mathcal{O})}^2+2\int_0^{t\land\tau^n_N}\|\theta^\epsilon_n(s)\|^2\,ds\\
\leq&\|P_n\theta_0\|_{L^2(\mathcal{O})}^2+2\int_0^{t\land\tau^n_N}\|g(s)\|_{L^2(\mathcal{O})}\|\theta^\epsilon_n(s)
\|_{L^2(\mathcal{O})}\,ds+2\sqrt{\epsilon}\int_0^{t\land\tau^n_N}\int_{\mathcal{O}}\theta^\epsilon_n(s)\sigma_n(s,\theta^\epsilon_n(s))\,dW(s)dxdydz\\
&+\epsilon\int_0^{t\land\tau^n_N}\|\sigma_n(s,\theta_n^\epsilon(s))\|_{\mathcal{L}_2(U_0,H_2)}^2\,ds+2\beta\int_0^{t\land\tau^n_N}\|\theta^*\|_{L^2(\Gamma_u)}\|\theta^\epsilon_n(s)\|_{L^2(\Gamma_u)}\,ds\\
\leq&\|\theta_0\|_{L^2(\mathcal{O})}^2+C\int_0^{t\land\tau^n_N}\|g(s)\|_{L^2(\mathcal{O})}^2\,ds+2\sqrt{\epsilon}\int_0^{t\land\tau^n_N}\int_{\mathcal{O}}\theta^\epsilon_n(s)\sigma_n(s,\theta^\epsilon_n(s))\,dW(s)dxdydz\\
&+\epsilon\int_0^{t\land\tau^n_N}\|\sigma(s,\theta^\epsilon_n(s))\|_{\mathcal{L}_2(U_0,H_2)}^2\,ds+Ct\|\theta^*\|_{L^2(\mathcal{M})} ^2+\int_0^{t\land\tau^n_N}\|\theta^\epsilon_n(s)
\|^2\,ds,
\end{align*}
which implies that 
\begin{align}\label{3.4}
\nonumber&\|\theta^\epsilon_n({t\land\tau^n_N})\|_{L^2(\mathcal{O})}^2+\int_0^{t\land\tau^n_N}\|\theta^\epsilon_n(s)\|^2\,ds\\
\nonumber\leq&\|\theta_0\|_{L^2(\mathcal{O})}^2+C\int_0^{t\land\tau^n_N}\|g(s)\|_{L^2(\mathcal{O})}^2\,ds+2\sqrt{\epsilon}\int_0^{t\land\tau^n_N}\int_{\mathcal{O}}\theta^\epsilon_n(s)\sigma_n(s,\theta^\epsilon_n(s))\,dW(s)dxdydz\\
&+\epsilon\int_0^{t\land\tau^n_N}\|\sigma(s,\theta^\epsilon_n(s))\|_{\mathcal{L}_2(U_0,H_2)}^2\,ds+Ct\|\theta^*\|_{L^2(\mathcal{M})} ^2,
\end{align}
Taking the supremum up to time $T$ in inequality \eqref{3.4} and taking the expectation on both hand sides of the resulting inequality, we obtain
\begin{align}\label{3.5}
\nonumber &E\sup_{0\leq t\leq T\land\tau^n_N}\left(\|\theta^\epsilon_n(t)\|_{L^2(\mathcal{O})}^2+\int_0^{t}\|\theta^\epsilon_n(s)\|^2\,ds\right)\\
\nonumber\leq&E\|\theta_0\|_{L^2(\mathcal{O})}^2+C\int_0^T\|g(s)\|_{L^2(\mathcal{O})}^2\,ds+2\sqrt{\epsilon}E\left(\sup_{0\leq r\leq T\land\tau^n_N}\left|\int_0^r\int_{\mathcal{O}}\theta^\epsilon_n(s)\sigma_n(s,\theta^\epsilon_n(s))\,dW(s)dxdydz\right|\right)\\
&+\epsilon E\left(\int_0^T\|\sigma(s,\theta^\epsilon_n(s))\|_{\mathcal{L}_2(U_0,H_2)}^2\,ds\right)+CT\|\theta^*\|_{L^2(\mathcal{M})}^2.
\end{align}
By the Burkholder-Davis-Gundy inequality, we obtain 
\begin{align}\label{3.6}
\nonumber&2\sqrt{\epsilon}E\left(\sup_{0\leq r\leq T\land\tau^n_N}\left|\int_0^r\int_{\mathcal{O}}\theta^\epsilon_n(s)\sigma_n(s,\theta^\epsilon_n(s))\,dW(s)dxdydz\right|\right)\\
\nonumber\leq &C\sqrt{\epsilon}E\left(\left(\int_0^{T\land\tau^n_N}\|\theta^\epsilon_n(s)\|_{L^2(\mathcal{O})}^2\|\sigma_n(s,\theta^\epsilon_n(s))\|_{\mathcal{L}_2(U_0,H_2)}^2\,ds\right)^{\frac{1}{2}}\right)\\
\nonumber\leq &C\sqrt{\epsilon}E\left(\sup_{0\leq r\leq T\land\tau^n_N}\|\theta^\epsilon_n(r)\|_{L^2(\mathcal{O})}\left(\int_0^{T\land\tau^n_N}\|\sigma_n(s,\theta^\epsilon_n(s))\|_{\mathcal{L}_2(U_0,H_2)}^2\,ds\right)^{\frac{1}{2}}\right)
\end{align}
\begin{align}
\leq &\frac{1}{2}E\left(\sup_{0\leq r\leq T\land\tau^n_N}\|\theta^\epsilon_n(r)\|_{L^2(\mathcal{O})}^2\right)+C\epsilon E\left(\int_0^{T\land\tau^n_N}\|\sigma(s,\theta^\epsilon_n(s))\|_{\mathcal{L}_2(U_0,H_2)}^2\,ds\right).
\end{align}
It follows from inequalities \eqref{3.5}-\eqref{3.6} and assumption ($A_2$) that 
\begin{align*}
\nonumber &E\left(\sup_{0\leq r\leq T\land\tau^n_N}\|\theta^\epsilon_n(r)\|_{L^2(\mathcal{O})}^2\right)+2E\left(\int_0^{T\land\tau^n_N}\|\theta^\epsilon_n(s)\|^2\,ds\right)\\
\nonumber\leq&2E\|\theta_0\|_{L^2(\mathcal{O})}^2+C\int_0^T\|g(s)\|_{L^2(\mathcal{O})}^2\,ds+CT\|\theta^*\|_{L^2(\mathcal{M})}^2+C\epsilon E\left(\int_0^T\|\sigma(s,\theta^\epsilon_n(s))\|_{\mathcal{L}_2(U,H_2)}^2\,ds\right)\\
\nonumber\leq&2E\|\theta_0\|_{L^2(\mathcal{O})}^2+C\int_0^T\|g(s)\|_{L^2(\mathcal{O})}^2\,ds+CT\|\theta^*\|_{L^2(\mathcal{M})}^2+CK\epsilon E\left(\int_0^T(1+\|\theta^\epsilon_n(s)\|_{L^2(\mathcal{O})}^2)\,ds\right)\\
\leq&2E\|\theta_0\|_{L^2(\mathcal{O})}^2+C\int_0^T\|g(s)\|_{L^2(\mathcal{O})}^2\,ds+CT(\|\theta^*\|_{L^2(\mathcal{M})}^2+1)+CK E\left(\int_0^T\sup_{0\leq r\leq s\land\tau^n_N}\|\theta^\epsilon_n(r)\|_{L^2(\mathcal{O})}^2\,ds\right),
\end{align*}
we conclude from Lemma \ref{2.2.4} that
\begin{align}\label{3.7}
\nonumber &E\left(\sup_{0\leq r\leq T\land\tau^n_N}\|\theta^\epsilon_n(r)\|_{L^2(\mathcal{O})}^2\right)+2E\left(\int_0^{T\land\tau^n_N}\|\theta^\epsilon_n(s)\|^2\,ds\right)\\
\nonumber\leq&\left(2E\|\theta_0\|_{L^2(\mathcal{O})}^2+C\int_0^T\|g(s)\|_{L^2(\mathcal{O})}^2\,ds+CT(\|\theta^*\|_{L^2(\mathcal{M})}^2+1)\right)e^{CKT}\\
=:&D(T),
\end{align}
which implies that for each natural number $n,$ $T\land\tau^n_N$ increases to $T$ a.s. as $N\rightarrow+\infty.$

Taking the limit in inequality \eqref{3.7} as $N\rightarrow +\infty,$ we infer from the Fatou's Lemma and inequality \eqref{3.7} that
\begin{align}\label{3.8}
E\left(\sup_{0\leq r\leq T}\|\theta^\epsilon_n(r)\|_{L^2(\mathcal{O})}^2\right)+E\left(\int_0^T\|\theta^\epsilon_n(s)\|^2\,ds\right)\leq D(T).
\end{align}
Define $\tau^n_R=\inf\{t>0:\|\theta^\epsilon_n(t)\|_{L^2(\mathcal{O})}^4>R\},$ we apply the finite dimensional It\^{o}'s formula (see \cite{in}) to the function $\|\theta^\epsilon_n(t)\|_{L^2(\mathcal{O})}^4,$ yields
\begin{align*}
&\|\theta^\epsilon_n(t\land\tau_R^n)\|_{L^2(\mathcal{O})}^4+4\int_0^{t\land\tau_R^n}\|\theta^\epsilon_n(s)\|_{L^2(\mathcal{O})}^2\|\theta^\epsilon_n(s)\|^2\,ds\\
\leq&\|P_n\theta_0\|_{L^2(\mathcal{O})}^4+4\int_0^{t\land\tau_R^n}\int_{\mathcal{O}}\|\theta^\epsilon_n(s)\|_{L^2(\mathcal{O})}^2g(s)\theta^\epsilon_n(s)
\,dxdydzds\\
&+4\sqrt{\epsilon}\int_0^{t\land\tau_R^n}\int_{\mathcal{O}}\|\theta^\epsilon_n(s)\|_{L^2(\mathcal{O})}^2\theta^\epsilon_n(s)\sigma_n(s,\theta^\epsilon_n(s))\,dW(s)dxdydz\\
&+2\epsilon\int_0^{t\land\tau_R^n}\|\theta^\epsilon_n(s)\|_{L^2(\mathcal{O})}^2\|\sigma_n(s,\theta^\epsilon_n(s))\|_{\mathcal{L}_2(U_0,H_2)}^2\,ds+4\epsilon\int_0^{t\land\tau_R^n}\|\theta^\epsilon_n(s)\|_{L^2(\mathcal{O})}^2\|\sigma_n(s,\theta^\epsilon_n(s))\|_{\mathcal{L}_2(U_0,H_2)}^2\,ds\\
&+4\beta\int_0^{t\land\tau_R^n}\int_{\Gamma_u}\|\theta^\epsilon_n(s)\|_{L^2(\mathcal{O})}^2\theta^* \theta^\epsilon_n(s)\,dxdyds.
\end{align*}
We deduce from Young's inequality, H\"{o}lder's inequality, inequality \eqref{3.3} and Lemma \ref{2.2.1} that 
\begin{align*}
&\|\theta^\epsilon_n(t\land\tau_R^n)\|_{L^2(\mathcal{O})}^4+4\int_0^{t\land\tau_R^n}\|\theta^\epsilon_n(s)\|_{L^2(\mathcal{O})}^2\|\theta^\epsilon_n(s)\|^2\,ds\\
\leq&\|\theta_0\|_{L^2(\mathcal{O})}^4+4\int_0^{t\land\tau_R^n}\|\theta^\epsilon_n(s)\|_{L^2(\mathcal{O})}^3\|g(s)\|_{L^2(\mathcal{O})}\,ds\\
&+4\sqrt{\epsilon}\int_0^{t\land\tau_R^n}\int_{\mathcal{O}}\|\theta^\epsilon_n(s)\|_{L^2(\mathcal{O})}^2\theta^\epsilon_n(s)\sigma_n(s,\theta^\epsilon_n(s))\,dW(s)dxdydz\\
&+6\epsilon\int_0^{t\land\tau_R^n}\|\theta^\epsilon_n(s)\|_{L^2(\mathcal{O})}^2\|\sigma_n(s,\theta^\epsilon_n(s))\|_{\mathcal{L}_2(U_0,H_2)}^2\,ds+4\beta\int_0^{t\land\tau_R^n}\|\theta^\epsilon_n(s)\|_{L^2(\mathcal{O})}^2\|\theta^*\|_{L^2(\mathcal{M})}\| \theta^\epsilon_n(s)\|_{L^2(\Gamma_u)}\,ds\\
\leq&\|\theta_0\|_{L^2(\mathcal{O})}^4+C\int_0^{t\land\tau_R^n}\|\theta^\epsilon_n(s)\|_{L^2(\mathcal{O})}^2\|\theta^\epsilon_n(s)\|\|g(s)\|_{L^2(\mathcal{O})}\,ds\\
&+4\sqrt{\epsilon}\int_0^{t\land\tau_R^n}\int_{\mathcal{O}}\|\theta^\epsilon_n(s)\|_{L^2(\mathcal{O})}^2\theta^\epsilon_n(s)\sigma_n(s,\theta^\epsilon_n(s))\,dW(s)dxdydz\\
&+6\epsilon\int_0^{t\land\tau_R^n}\|\theta^\epsilon_n(s)\|_{L^2(\mathcal{O})}^2\|\sigma_n(s,\theta^\epsilon_n(s))\|_{\mathcal{L}_2(U_0,H_2)}^2\,ds+4\beta\int_0^{t\land\tau_R^n}\|\theta^\epsilon_n(s)\|_{L^2(\mathcal{O})}^2\|\theta^*\|_{L^2(\mathcal{M})}\| \theta^\epsilon_n(s)\|_{L^2(\Gamma_u)}\,ds\\
\leq&\|\theta_0\|_{L^2(\mathcal{O})}^4+C\int_0^{t\land\tau_R^n}\|\theta^\epsilon_n(s)\|_{L^2(\mathcal{O})}^2\|g(s)\|_{L^2(\mathcal{O})}^2\,ds\\
&+4\sqrt{\epsilon}\int_0^{t\land\tau_R^n}\int_{\mathcal{O}}\|\theta^\epsilon_n(s)\|_{L^2(\mathcal{O})}^2\theta^\epsilon_n(s)\sigma_n(s,\theta^\epsilon_n(s))\,dW(s)dxdydz\\
&+6\epsilon\int_0^{t\land\tau_R^n}\|\theta^\epsilon_n(s)\|_{L^2(\mathcal{O})}^2\|\sigma(s,\theta^\epsilon_n(s))\|_{\mathcal{L}_2(U_0,H_2)}^2\,ds+C\int_0^{t\land\tau_R^n}\|\theta^\epsilon_n(s)\|_{L^2(\mathcal{O})}^2\|\theta^*\|_{L^2(\mathcal{M})}^2\,ds\\
&+2\int_0^{t\land\tau_R^n}\|\theta^\epsilon_n(s)\|_{L^2(\mathcal{O})}^2\|\theta^\epsilon_n(s)\|^2\,ds,
\end{align*}
which entails that
\begin{align}\label{3.9}
\nonumber&\|\theta^\epsilon_n(t\land\tau_R^n)\|_{L^2(\mathcal{O})}^4+2\int_0^{t\land\tau_R^n}\|\theta^\epsilon_n(s)\|_{L^2(\mathcal{O})}^2\|\theta^\epsilon_n(s)\|^2\,ds\\
\nonumber\leq&\|\theta_0\|_{L^2(\mathcal{O})}^4+C\int_0^{t\land\tau_R^n}\|\theta^\epsilon_n(s)\|_{L^2(\mathcal{O})}^2\|g(s)\|_{L^2(\mathcal{O})}^2\,ds+6\epsilon\int_0^{t\land\tau_R^n}\|\theta^\epsilon_n(s)\|_{L^2(\mathcal{O})}^2\|\sigma(s,\theta_n^\epsilon(s))\|_{\mathcal{L}_2(U_0,H_2)}^2\,ds\\
&+4\sqrt{\epsilon}\int_0^{t\land\tau_R^n}\int_{\mathcal{O}}\|\theta^\epsilon_n(s)\|_{L^2(\mathcal{O})}^2\theta^\epsilon_n(s)\sigma_n(s,\theta_n^\epsilon(s))\,dW(s)dxdydz+C\int_0^{t\land\tau_R^n}\|\theta^\epsilon_n(s)\|_{L^2(\mathcal{O})}^2\,ds.
\end{align}
Taking the supremum up to time $T$ in inequality \eqref{3.9} and taking the expectation on both hand sides of the resulting inequality, we obtain
\begin{align}\label{3.10}
\nonumber &E\sup_{0\leq t\leq T\land\tau^n_R}\left(\|\theta^\epsilon_n(t)\|_{L^2(\mathcal{O})}^4+2\int_0^{t}\|\theta^\epsilon_n(s)\|_{L^2(\mathcal{O})}^2\|\theta^\epsilon_n(s)\|^2\,ds\right)\\
\nonumber\leq&E\|\theta_0\|_{L^2(\mathcal{O})}^4+CE\left(\int_0^{T\land\tau_R^n}\|\theta^\epsilon_n(s)\|_{L^2(\mathcal{O})}^2\|g(s)\|_{L^2(\mathcal{O})}^2\,ds\right)
\end{align}
\begin{align}
\nonumber&+6\epsilon E\left(\int_0^{T\land\tau_R^n}\|\theta^\epsilon_n(s)\|_{L^2(\mathcal{O})}^2\|\sigma(s,\theta^\epsilon_n(s))\|_{\mathcal{L}_2(U_0,H_2)}^2\,ds\right)\\
\nonumber&+4\sqrt{\epsilon}E\left(\sup_{0\leq r\leq T\land\tau^n_R}\left|\int_0^r\int_{\mathcal{O}}\|\theta^\epsilon_n(s)\|_{L^2(\mathcal{O})}^2\theta^\epsilon_n(s)\sigma_n(s,\theta^\epsilon_n(s))\,dW(s)dxdydz\right|\right)\\
&+CE\left(\int_0^{T\land\tau_R^n}\|\theta^\epsilon_n(s)\|_{L^2(\mathcal{O})}^2\,ds\right).
\end{align}
It follows from the Burkholder-Davis-Gundy inequality that
\begin{align}\label{3.11}
\nonumber&4\sqrt{\epsilon}E\left(\sup_{0\leq r\leq T\land\tau^n_R}\left|\int_0^r\int_{\mathcal{O}}\|\theta^\epsilon_n(s)\|_{L^2(\mathcal{O})}^2\theta^\epsilon_n(s)\sigma_n(s,\theta_n^\epsilon(s))\,dW(s)dxdydz\right|\right)\\
\nonumber\leq &C\sqrt{\epsilon}E\left(\left(\int_0^{T\land\tau^n_R}\|\theta^\epsilon_n(s)\|_{L^2(\mathcal{O})}^4\|\theta^\epsilon_n(s)\|_{L^2(\mathcal{O})}^2\|\sigma_n(s,\theta^\epsilon_n(s))\|_{\mathcal{L}_2(U_0,H_2)}^2\,ds\right)^{\frac{1}{2}}\right)\\
\nonumber\leq &C\sqrt{\epsilon}E\left(\sup_{0\leq r\leq T\land\tau^n_R}\|\theta^\epsilon_n(r)\|_{L^2(\mathcal{O})}^2\left(\int_0^{T\land\tau_R^n}\|\theta^\epsilon_n(s)\|_{L^2(\mathcal{O})}^2\|\sigma_n(s,\theta^\epsilon_n(s))\|_{\mathcal{L}_2(U_0,H_2)}^2\,ds\right)^{\frac{1}{2}}\right)\\
\leq &\frac{1}{2}E\left(\sup_{0\leq r\leq T\land\tau^n_R}\|\theta^\epsilon_n(r)\|_{L^2(\mathcal{O})}^4\right)+C\epsilon E\left(\int_0^{T\land\tau_R^n}\|\theta^\epsilon_n(s)\|_{L^2(\mathcal{O})}^2\|\sigma(s,\theta^\epsilon_n(s))\|_{\mathcal{L}_2(U_0,H_2)}^2\,ds\right).
\end{align}
Therefore, we conclude from inequalities \eqref{3.10}-\eqref{3.11} and assumption ($A_2$) that
\begin{align}\label{3.12}
\nonumber &E\left(\sup_{0\leq t\leq T\land\tau^n_R}\|\theta^\epsilon_n(t)\|_{L^2(\mathcal{O})}^4\right)+4E\left(\int_0^{T\land\tau_R^n}\|\theta^\epsilon_n(s)\|_{L^2(\mathcal{O})}^2\|\theta^\epsilon_n(s)\|^2\,ds\right)\\
\nonumber\leq&2E\|\theta_0\|_{L^2(\mathcal{O})}^4+CE\left(\int_0^{T\land\tau_R^n}\|\theta^\epsilon_n(s)\|_{L^2(\mathcal{O})}^2\|g(s)\|_{L^2(\mathcal{O})}^2\,ds\right)+CE\left(\int_0^{T\land\tau_R^n}\|\theta^\epsilon_n(s)\|_{L^2(\mathcal{O})}^2\,ds\right)\\
\nonumber&+C\epsilon E\left(\int_0^{T\land\tau_R^n}\|\theta^\epsilon_n(s)\|_{L^2(\mathcal{O})}^2\|\sigma(s,\theta^\epsilon_n(s))\|_{\mathcal{L}_2(U_0,H_2)}^2\,ds\right)\\
\nonumber\leq&2E\|\theta_0\|_{L^2(\mathcal{O})}^4+CE\left(\sup_{0\leq r\leq T\land\tau^n_R}\|\theta^\epsilon_n(r)\|_{L^2(\mathcal{O})}^2\int_0^{T\land\tau_R^n}\|g(s)\|_{L^2(\mathcal{O})}^2\,ds\right)\\
\nonumber&+CE\left(\int_0^{T\land\tau_R^n}\|\theta^\epsilon_n(s)\|_{L^2(\mathcal{O})}^2\,ds\right)
+CK\epsilon E\left(\int_0^{T\land\tau_R^n}\|\theta^\epsilon_n(s)\|_{L^2(\mathcal{O})}^2(1+\|\theta^\epsilon_n(s)\|_{L^2(\mathcal{O})}^2)\,ds\right)\\
\nonumber\leq&2E\|\theta_0\|_{L^2(\mathcal{O})}^4+C\left(\int_0^T\|g(s)\|_{L^2(\mathcal{O})}^2\,ds\right)^2+CE\left(\int_0^{T\land\tau_R^n}\|\theta^\epsilon_n(s)\|^2\,ds\right)\\
&+C E\left(\int_0^{T\land\tau_R^n}(1+\|\theta^\epsilon_n(s)\|_{L^2(\mathcal{O})}^4)\,ds\right)+\frac{1}{2}E\left(\sup_{0\leq t\leq T\land\tau^n_R}\|\theta^\epsilon_n(t)\|_{L^2(\mathcal{O})}^4\right),
\end{align}
which implies that
\begin{align*}
&E\left(\sup_{0\leq t\leq T\land\tau^n_R}\|\theta^\epsilon_n(t)\|_{L^2(\mathcal{O})}^4\right)+8E\left(\int_0^{T\land\tau_R^n}\|\theta^\epsilon_n(s)\|_{L^2(\mathcal{O})}^2\|\theta^\epsilon_n(s)\|^2\,ds\right)\\
\leq&4E\|\theta_0\|_{L^2(\mathcal{O})}^4+C\left(\int_0^T\|g(s)\|_{L^2(\mathcal{O})}^2\,ds\right)^2+CE\left(\int_0^T\|\theta^\epsilon_n(s)\|^2\,ds\right)+CT\\
&+C E\left(\int_0^T\sup_{0\leq r\leq s\land\tau^n_R}\|\theta^\epsilon_n(r)\|_{L^2(\mathcal{O})}^4\,ds\right),
\end{align*}
we infer from Lemma \ref{2.2.4} and inequality \eqref{3.8} that
\begin{align}\label{3.13}
\nonumber &E\left(\sup_{0\leq t\leq T\land\tau^n_R}\|\theta^\epsilon_n(t)\|_{L^2(\mathcal{O})}^4\right)+8E\left(\int_0^{T\land\tau_R^n}\|\theta^\epsilon_n(s)\|_{L^2(\mathcal{O})}^2\|\theta^\epsilon_n(s)\|^2\,ds\right)\\
\nonumber\leq&\left(4E\|\theta_0\|_{L^2(\mathcal{O})}^4+C\left(\int_0^T\|g(s)\|_{L^2(\mathcal{O})}^2\,ds\right)^2+CD(T)+CT\right)e^{CT}\\
=:&D_1(T).
\end{align}
Thus, for each natural number $n,$ $T\land\tau^n_R$ increases to $T$ a.s. as $R\rightarrow+\infty.$

Taking the limit in inequality \eqref{3.13} as $R\rightarrow +\infty,$ we infer from the Fatou's Lemma and inequality \eqref{3.8} that
\begin{align}\label{3.14}
E\left(\sup_{0\leq r\leq T}\|\theta^\epsilon_n(r)\|_{L^2(\mathcal{O})}^4\right)+E\left(\int_0^T\|\theta^\epsilon_n(s)\|_{L^2(\mathcal{O})}^2\|\theta^\epsilon_n(s)\|^2\,ds\right)\leq D_1(T),
\end{align}
where
\begin{align*}
D_1(T)=\left(4E\|\theta_0\|_{L^2(\mathcal{O})}^4+C\left(\int_0^T\|g(s)\|_{L^2(\mathcal{O})}^2\,ds\right)^2+CD(T)+CT\right)e^{CT}.
\end{align*}
For brevity, let us put in the following notation:
\begin{align*}
F(\theta^\epsilon)=-B(v^\epsilon,\theta^\epsilon)-A_2(\theta^\epsilon-\theta^*)+K_h\Delta \theta^*+g(x,y,z,t),
\end{align*}
where $(v^\epsilon,p^\epsilon_s)=(v^\epsilon(\theta^\epsilon),p^\epsilon_s(\theta^\epsilon))$ is established by the first equation and the second equation of equations \eqref{2.2.1} for any given $\theta^\epsilon,$ then $v$ satisfies the following estimates:
\begin{align}\label{3.15}
\|v^\epsilon(t)\|_{H^{\gamma+1}(\mathcal{O})}^2+\|p^\epsilon_{s}\|_{H^\gamma(\mathcal{M})}^2\leq \mathcal{K}_1\left(\|\theta^\epsilon\|_{H^\gamma(\mathcal{O})}^2+\|\mu\|_{H^1(\mathcal{M})}^2\right)
\end{align}
for $\gamma=0,1.$

For any $w\in L^2(\Omega\times(0,T); V_2),$ let $w_n=P_nw$ and $v^\epsilon_n=v^\epsilon(\theta^\epsilon_n),$ we conclude from H\"{o}lder inequality and interpolation inequality that
\begin{align}\label{3.16}
\nonumber&E\left(\int_0^T\left|\langle P_nF(\theta^\epsilon_n(s)),w(s)\rangle\right|\,ds\right)\\
\nonumber\leq&\mathcal{K}E\left(\int_0^T\|v^\epsilon_n(s)\|_{H^1(\mathcal{O})}^{\frac{1}{2}}\|v^\epsilon_n(s)\|_{H^2(\mathcal{O})}^{\frac{1}{2}}\|\theta^\epsilon_n(s)\|_{L^2(\mathcal{O})}^{\frac{1}{2}}\|\theta^\epsilon_n(s)\|^{\frac{1}{2}}\|w_n(s)\|\,ds\right)+E\left(\int_0^T\|\theta^\epsilon_n(s)\|\|w_n(s)\|\,ds\right)\\
\nonumber&+\beta E\left(\int_0^T\|\theta^*\|_{L^2(\mathcal{M})}\|w_n(s)\|_{L^2(\Gamma_u)}\,ds\right)+E\left(\int_0^T\|g(s)\|_{L^2(\mathcal{O})}\|w_n(s)\|_{L^2(\mathcal{O})}\,ds\right)\\
\nonumber\leq&\mathcal{K}\left(E\int_0^T\|v^\epsilon_n(s)\|_{H^1(\mathcal{O})}^2\|v^\epsilon_n(s)\|_{H^2(\mathcal{O})}^2\,ds\right)^{\frac{1}{4}}\left(E\int_0^T\|\theta^\epsilon_n(s)\|_{L^2(\mathcal{O})}^2\|\theta^\epsilon_n(s)\|^2\,ds\right)^{\frac{1}{4}}\left(E\int_0^T\|w(s)\|^2\,ds\right)^{\frac{1}{2}}\\
\nonumber&+E\left(\int_0^T\|\theta^\epsilon_n(s)\|^2\,ds\right)^{\frac{1}{2}}\left(E\int_0^T\|w(s)\|^2\,ds\right)^{\frac{1}{2}}+\beta\|\theta^*\|_{L^2(\mathcal{M})}T^{\frac{1}{2}} E\left(\int_0^T\|w(s)\|_{L^2(\Gamma_u)}^2\,ds\right)^{\frac{1}{2}}\\
&+\left(\int_0^T\|g(s)\|_{L^2(\mathcal{O})}^2\right)^{\frac{1}{2}}E\left(\int_0^T\|w(s)\|_{L^2(\mathcal{O})}^2\,ds\right)^{\frac{1}{2}}.
\end{align}
We infer from Lemma \ref{2.2.1}, inequalities \eqref{3.8}, \eqref{3.13}, \eqref{3.15}-\eqref{3.16} and assumption ($A_2$) that
\begin{align}\label{3.17}
\nonumber&\|P_nF(\theta^\epsilon_n)\|_{L^2(\Omega\times(0,T); V_2')}=\sup_{\|w\|_{L^2(\Omega\times(0,T); V_2)}\leq 1}\left|E\left(\int_0^T\langle P_nF(\theta^\epsilon_n(s)),w(s)\rangle\,ds\right)\right|\\
\nonumber\leq&C\left(E\int_0^T\|v^\epsilon_n(s)\|_{H^1(\mathcal{O})}^4\,ds\right)^{\frac{1}{4}}\left(E\int_0^T\|\theta^\epsilon_n(s)\|_{L^2(\mathcal{O})}^2\|\theta^\epsilon_n(s)\|^2\,ds\right)^{\frac{1}{4}}\\
\nonumber&+C\left(E\int_0^T\|v^\epsilon_n(s)\|_{H^1(\mathcal{O})}^2\|v^\epsilon_n(s)\|_{H^2(\mathcal{O})}^2\,ds\right)^{\frac{1}{4}}\left(E\int_0^T\|\theta^\epsilon_n(s)\|_{L^2(\mathcal{O})}^2\|\theta^\epsilon_n(s)\|^2\,ds\right)^{\frac{1}{4}}\\
\nonumber&+E\left(\int_0^T\|\theta^\epsilon_n(s)\|^2\,ds\right)^{\frac{1}{2}}+C\beta\|\theta^*\|_{L^2(\mathcal{M})}T^{\frac{1}{2}} +C\left(\int_0^T\|g(s)\|_{L^2(\mathcal{O})}^2\right)^{\frac{1}{2}}\\
<&+\infty
\end{align}
and
\begin{align}\label{3.18}
\nonumber &E\left(\int_0^T\|p^\epsilon_{sn}(s)\|_{H^1(M)}^2\,ds\right)+E\left(\int_0^T\|v^\epsilon_n(s)\|_{H^1(\mathcal{O})}^4\,ds\right)+E\left(\int_0^T\|\theta^\epsilon_n(s)\|_{L^3(\mathcal{O})}^4\,ds\right)\\
\nonumber&+E\left(\int_0^T\|\sigma_n(s,\theta^\epsilon_n(s))\|_{\mathcal{L}_2(U_0,H_2)}^2\,ds\right)\\
\nonumber\leq&CE\left(\int_0^T\|\theta^\epsilon_n(s)\|_{L^2(\mathcal{O})}^2\|\theta^\epsilon_n(s)\|^2\,ds\right)+TE\left(\sup_{0\leq r\leq T}\|\theta^\epsilon_n(r)\|_{L^2(\mathcal{O})}^4\right)+CT(1+\|\mu\|_{H^1(\mathcal{M})}^4)\\
\leq &CD_1(T)+CT(1+\|\mu\|_{H^1(\mathcal{M})}^4).
\end{align}
 Therefore, there exists a subsequence of $\{(v^\epsilon_n,p^\epsilon_{sn},\theta^\epsilon_n):n\in\mathbb{N}\}$ (still denoted by using the same notation) of processes and elements
 \begin{align*}
 &\bar{\theta}^\epsilon\in L^2(\Omega\times(0,T);V_2)\cap L^4(\Omega\times (0,T);L^3(\mathcal{O}))\cap L^4(\Omega;L^{\infty}(0,T;H_2)),\\
 &v^\epsilon\in L^4(\Omega\times (0,T);V_1),\\
 &u^\epsilon\in L^2(\Omega;H_2),\\
 &p^\epsilon_s\in L^2(\Omega\times(0,T);H^1(\mathcal{M})) ,\\
 &F\in L^2(\Omega\times (0,T);V_2'),\\
 &S\in L^2(\Omega\times (0,T);\mathcal{L}_2(U_0,H_2)),
 \end{align*}
such that
\begin{align*}
&\theta^\epsilon_n\rightharpoonup \bar{\theta}^\epsilon\;\textit{weakly in}\; L^2(\Omega\times(0,T);V_2),\\
&\theta^\epsilon_n\rightharpoonup \bar{\theta}^\epsilon\;\textit{weakly in}\; L^4(\Omega\times (0,T);L^3(\mathcal{O})),\\
&\theta^\epsilon_n\rightharpoonup \bar{\theta}\;\textit{weakly star in}\; L^4(\Omega;L^{\infty}(0,T;H_2)),\\
&\theta^\epsilon_n(T)\rightharpoonup u^\epsilon\;\textit{weakly in}\; L^2(\Omega;H_2),\\
&v^\epsilon_n\rightharpoonup v^\epsilon\;\textit{weakly in}\; L^4(\Omega\times (0,T);V_1),\\
&p^\epsilon_{sn}\rightharpoonup p^\epsilon_s\;\textit{weakly in}\; L^2(\Omega\times (0,T);H^1(\mathcal{M})),\\
&P_nF(\theta^\epsilon_n)\rightharpoonup F\;\textit{weakly in}\; L^2(\Omega\times (0,T);V_2'),\\
&\sigma_n(\cdot,\theta^\epsilon_n(\cdot))\rightharpoonup S\;\textit{weakly in}\; L^2(\Omega\times (0,T);\mathcal{L}_2(U_0,H_2)).
\end{align*}
For any $\chi\in L^2(\Omega\times(0,T);\mathbb{R})$ and $\psi\in H^1(\mathcal{O}),$ we obtain
\begin{align}\label{3.19}
\nonumber&E\int_0^T\int_\mathcal{O}\chi(s) \nabla p^\epsilon_{sn}(x,y,s)\cdot\psi\,dxdydzds-E\int_0^T\int_{\mathcal{O}}\chi(s)\left(\int_{-h}^z
\nabla \theta^\epsilon_n(x,y,\zeta,s)d\zeta\right)\cdot\psi\,dxdydzds\\
\nonumber&+E\int_0^T\int_{\mathcal{O}}\chi(s)f{v^\epsilon_n}^\bot\cdot\psi\,dxdydzds+E\int_0^T\chi(s)\langle A_1v^\epsilon_n(s), \phi\rangle\,ds\\
&=\kappa E\int_0^T\int_{\Gamma_u} \chi(s)\mu\cdot\phi\, dxdyds.
\end{align}
Taking the limit in equality \eqref{3.19} as $n\rightarrow+\infty,$ we conclude that
\begin{align*}
&E\int_0^T\int_\mathcal{O}\chi(s) \nabla p^\epsilon_s(x,y,s)\cdot\psi\,dxdydzds-E\int_0^T\int_{\mathcal{O}}\chi(s)\left(\int_{-h}^z
\nabla \bar{\theta}^\epsilon(x,y,\zeta,s)d\zeta\right)\cdot\psi\,dxdydzds\\
&+E\int_0^T\int_{\mathcal{O}}\chi(s)f{v^\epsilon}^\bot\cdot\psi\,dxdydzds+E\int_0^T\chi(s)\langle A_1v^\epsilon(s), \phi\rangle\,ds\\
&=\kappa E\int_0^T\int_{\Gamma_u} \chi(s)\mu\cdot\phi\, dxdyds
\end{align*}
for any $\chi\in L^2(\Omega\times(0,T);\mathbb{R})$ and $\psi\in H^1(\mathcal{O}),$ which implies that
\begin{align*}
&\int_\mathcal{O} \nabla p^\epsilon_s(x,y,t)\cdot\phi\,dxdydz-\int_{\mathcal{O}}\left(\int_{-h}^z
\nabla \bar{\theta}^\epsilon(x,y,\zeta,t)d\zeta\right)\cdot\phi\,dxdydz+\int_{\mathcal{O}}f{v^\epsilon}^\bot\cdot\phi\,dxdydz
+\langle A_1v^\epsilon, \phi\rangle\\
&=\int_{\Gamma_u} \kappa\mu\cdot\phi\, dxdy,\,\,\,\,d\mathcal{P}\otimes dt-a.e.\,\,\,\textit{on}\,\,\,\Omega\times(0,T),
\end{align*}
Similarly, we can obtain
\begin{align*}
\bar{\theta}^\epsilon(t)=\theta_0+\int_0^t F(s)\,ds+\sqrt{\epsilon}\int_0^tS(s)\,dW(s),\,\,\,\,d\mathcal{P}\otimes dt-a.e.\,\,\,\textit{on}\,\,\,\Omega\times(0,T)
\end{align*}
and
\begin{align*}
u^\epsilon=\theta_0+\int_0^T F(s)\,ds+\sqrt{\epsilon}\int_0^T S(s)\,dW(s),\,\,\,\,d\mathcal{P}-a.s.
\end{align*}
Define a $V_2'$-valued process $\theta^\epsilon$ by 
\begin{align*}
\theta^\epsilon(t)=\theta_0+\int_0^t F(s)\,ds+\sqrt{\epsilon}\int_0^t S(s)\,dW(s),\,\,\,\,t\in(0,T)\,\,\,\textit{in}\,\,\,V_2',
\end{align*}
then, $\theta^\epsilon$ is a $V_2'$-valued modification of the $V_2$-valued process $\bar{\theta}^\epsilon\in L^2(\Omega\times (0,T); V_2)$ and
\begin{align*}
\theta^\epsilon(T)=u^\epsilon,\,\,\,d\mathcal{P}-a.s.
\end{align*}
Therefore, $\theta^\epsilon$ is an $H_2$-valued c\`{a}dl\`{a}g $(\mathcal{F}_t)$-adapted process, and for every $t\in[0,S],$ the following formula holds $\mathcal{P}$-a.s.
\begin{align}\label{3.20}
\nonumber\|\theta^\epsilon(t)\|_{L^2(\mathcal{O})}^2=&\|\theta_0\|_{L^2(\mathcal{O})}^2+2\int_0^t \langle F(s),\theta^\epsilon(s)\rangle\,ds+\epsilon\int_0^t\|S(s)\|_{\mathcal{L}_2(U_0,H_2)}^2\,ds\\
&+2\sqrt{\epsilon}\int_0^t\int_{\mathcal{O}}\theta^\epsilon(s)S(s)\,dxdydzdW(s).
\end{align}
In what follows, we only need to prove that
\begin{align*}
F(s,\omega)=&F(\theta^\epsilon(s,\omega))\,\,\,\textit{for}\,\,dt\otimes d\mathcal{P} -a.e.\,\,\,(s,\omega)\in (0,T)\times\Omega,\\
S(s,\omega)=&\sigma(s,\theta^\epsilon(s,\omega))\,\,\,\textit{for}\,\,dt\otimes d\mathcal{P} -a.e.\,\,\,(s,\omega)\in (0,T)\times\Omega.
\end{align*}
To establish these relation, we use the same idea as in \cite{bz, ghj}. For any natural number $m\leq n,$ let $\eta$ be a progressively measurable process belonging to $L^2(\Omega\times(0,T);V_2\cap H_m)\cap L^4(\Omega\times (0,T);L^3(\mathcal{O})\cap H_m)\cap L^4(\Omega;L^{\infty}(0,T;H_2\cap H_m)).$ For any $t\in [0,T],$ define
\begin{align*}
r(t)=\int_0^t(L_1\|\eta(s)\|^2+L)\,ds,\,\,\,t\in[0,T],
\end{align*}
where $L_1\geq \frac{\mathcal{K}^2\mathcal{K}_1}{2},$ then we infer from assumption ($A_3$), Lemma \ref{2.2.3} and H\"{o}lder's inequality that the following conclusion holds:
\begin{align}\label{3.21}
\nonumber &E\left[-\int_0^Te^{-r(s)}r'(s)\|\theta^\epsilon_n(s)-\eta(s)\|_{L^2(\mathcal{O})}^2\,ds+2\int_0^Te^{-r(s)}\langle P_nF(\theta^\epsilon_n(s))-P_nF(\eta(s)),\theta^\epsilon_n(s)-\eta(s)\rangle\,ds\right.\\
\nonumber&\left.+\int_0^Te^{-r(s)}\|\sigma_n(s,\theta^\epsilon_n(s))-\sigma_n(s,\eta(s))|_{\mathcal{L}_2(U_0,H_2)}^2\,ds
\right]\\
\nonumber \leq&E\left[-L_1\int_0^Te^{-r(s)}\|\eta(s)\|^2\|\theta^\epsilon_n(s)-\eta(s)\|_{L^2(\mathcal{O})}^2\,ds-2\int_0^Te^{-r(s)}\|\theta^\epsilon_n(s)-\eta(s)\|^2\,ds\right]
\end{align}
\begin{align}
\nonumber&-2E\left[\int_0^Te^{-r(s)}b(v^\epsilon_n(\theta^\epsilon_n)-v^\epsilon_n(\eta),\eta(s),\theta^\epsilon_n(s)-\eta(s))\,ds\right]\\
\nonumber \leq&E\left[-L_1\int_0^Te^{-r(s)}\|\eta(s)\|^2\|\theta^\epsilon_n(s)-\eta(s)\|_{L^2(\mathcal{O})}^2\,ds-2\int_0^Te^{-r(s)}\|\theta^\epsilon_n(s)-\eta(s)\|^2\,ds\right]\\
&+2\mathcal{K}E\left[\int_0^Te^{-r(s)}\|v^\epsilon_n(\theta^\epsilon_n)-v^\epsilon_n(\eta)\|_{H^1(\mathcal{O})}^{\frac{1}{2}}\|v^\epsilon_n(\theta^\epsilon_n)-v^\epsilon_n(\eta)\|_{H^2(\mathcal{O})}^{\frac{1}{2}}\|\theta^\epsilon_n(s)-\eta(s)\|_{L^2(\mathcal{O})}^{\frac{1}{2}}\|\theta^\epsilon_n(s)-\eta(s)\|^{\frac{1}{2}}\|\eta(s)\|\,ds\right].
\end{align}
It follows from the definition of $v^\epsilon_n(\theta^\epsilon_n)$ and $v^\epsilon_n(\eta)$ as well as Lemma \ref{2.2.2} that
\begin{align}\label{3.22}
\|v^\epsilon_n(\theta^\epsilon_n)(t)-v^\epsilon_n(\eta)(t)\|^2_{H^{\gamma+1}(\mathcal{O})}\leq \mathcal{K}_1\|\theta^\epsilon_n(t)-\eta(t)\|^2_{H^\gamma(\mathcal{O})}
\end{align}
for $\gamma=0,1.$

We deduce from inequalities \eqref{3.21}-\eqref{3.22} that
\begin{align}\label{3.23}
\nonumber &E\left[-\int_0^Te^{-r(s)}r'(s)\|\theta^\epsilon_n(s)-\eta(s)\|_{L^2(\mathcal{O})}^2\,ds+2\int_0^Te^{-r(s)}\langle P_nF(\theta^\epsilon_n(s))-P_nF(\eta(s)),\theta^\epsilon_n(s)-\eta(s)\rangle\,ds\right]\\
\nonumber \leq&E\left[-L_1\int_0^Te^{-r(s)}\|\eta(s)\|^2\|\theta^\epsilon_n(s)-\eta(s)\|_{L^2(\mathcal{O})}^2\,ds-2\int_0^Te^{-r(s)}\|\theta^\epsilon_n(s)-\eta(s)\|^2\,ds\right]\\
\nonumber&+2\mathcal{K}\sqrt{\mathcal{K}_1}E\left[\int_0^Te^{-r(s)}\|\theta^\epsilon_n(s)-\eta(s)\|_{L^2(\mathcal{O})}\|\theta^\epsilon_n(s)-\eta(s)\|\|\eta(s)\|\,ds\right]\\
\leq&0.
\end{align}

Since $\theta^\epsilon_n(T)\rightharpoonup u^\epsilon$ weakly in $L^2(\Omega;H_2)$ and $u^\epsilon=\theta^\epsilon(T)$ as well as $E\|P_n\theta_0\|_{L^2(\mathcal{O})}^2\leq E\|\theta_0\|_{L^2(\mathcal{O})}^2,$ we obtain
\begin{align}\label{3.24}
\left(\|\theta^\epsilon(T)\|_{L^2(\mathcal{O})}^2e^{-r(T)}\right)-E\|\theta_0\|_{L^2(\mathcal{O})}^2\leq\liminf_{n\rightarrow+\infty}\left[E\left(\|\theta^\epsilon_n(T)\|_{L^2(\mathcal{O})}^2e^{-r(T)}\right)-E\|P_n\theta_0\|_{L^2(\mathcal{O})}^2\right]
\end{align}
We apply the It\^{o}'s formula to the process $\|\theta^\epsilon(s)\|_{L^2(\mathcal{O})}^2e^{-r(s)},$ yield
\begin{align}\label{3.25}
\nonumber\|\theta^\epsilon(T)\|_{L^2(\mathcal{O})}^2e^{-r(T)}=&\|\theta_0\|_{L^2(\mathcal{O})}^2-\int_0^T e^{-r(s)}r'(s) \|\theta^\epsilon(s)\|_{L^2(\mathcal{O})}^2\,ds+2\int_0^T e^{-r(s)} \langle F(s),\theta^\epsilon(s)\rangle\,ds\\
&+2\sqrt{\epsilon}\int_0^T\int_{\mathcal{O}}e^{-r(s)}\theta^\epsilon(s)S(s)\,dW(s)dxdydz+\epsilon\int_0^Te^{-r(s)}\|S(s)\|_{\mathcal{L}_2(U_0,H_2)}^2\,ds.
\end{align}
Taking the expectation of both sides of equality \eqref{3.25}, yields
\begin{align}\label{3.26}
\nonumber E\left(\|\theta^\epsilon(T)\|_{L^2(\mathcal{O})}^2e^{-r(T)}\right)-E\|\theta_0\|_{L^2(\mathcal{O})}^2=&-E\int_0^T e^{-r(s)}r'(s) \|\theta^\epsilon(s)\|_{L^2(\mathcal{O})}^2\,ds+2E\int_0^T e^{-r(s)} \langle F(s),\theta^\epsilon(s)\rangle\,ds\\
&+\epsilon E\int_0^Te^{-r(s)}\|S(s)\|_{\mathcal{L}_2(U_0,H_2)}^2\,ds.
\end{align}
We also apply the finite dimensional It\^{o}'s formula to the $H_n$-valued process $\|\theta^\epsilon_n(s)\|_{L^2(\mathcal{O})}^2e^{-r(s)},$ yield
\begin{align}\label{3.27}
\nonumber&\|\theta^\epsilon_n(T)\|_{L^2(\mathcal{O})}^2e^{-r(T)}=\|P_n\theta_0\|_{L^2(\mathcal{O})}^2-\int_0^T e^{-r(s)}r'(s) \|\theta^\epsilon_n(s)\|_{L^2(\mathcal{O})}^2\,ds+2\int_0^T e^{-r(s)} \langle P_nF(\theta^\epsilon_n(s)),\theta^\epsilon_n(s)\rangle\,ds\\
&+2\sqrt{\epsilon}\int_0^T\int_{\mathcal{O}}e^{-r(s)}\theta^\epsilon_n(s)\sigma_n(s,\theta^\epsilon_n(s))\,dW(s)dxdydz+\epsilon\int_0^Te^{-r(s)}\|\sigma_n(s,\theta^\epsilon_n(s))\|_{\mathcal{L}_2(U_0,H_2)}^2\,ds.
\end{align}
Taking the expectation of both hand sides of equality \eqref{3.27} and using twice the identity $|x|^2=|x-y|^2-|y|^2+2(x,y)$ as well as inequality \eqref{3.23}, we obtain
\begin{align}\label{3.28}
\nonumber &E\left(\|\theta^\epsilon_n(T)\|_{L^2(\mathcal{O})}^2e^{-r(T)}\right)-E\|P_n\theta_0\|_{L^2(\mathcal{O})}^2=-E\int_0^T e^{-r(s)}r'(s) \|\theta^\epsilon_n(s)\|_{L^2(\mathcal{O})}^2\,ds\\
\nonumber&+2E\int_0^T e^{-r(s)} \langle P_nF(\theta^\epsilon_n(s)),\theta^\epsilon_n(s)\rangle\,ds+\epsilon E\int_0^Te^{-r(s)}\|\sigma_n(s,\theta^\epsilon_n(s))\|_{\mathcal{L}_2(U_0,H_2)}^2\,ds\\
\nonumber=&-E\int_0^T e^{-r(s)}r'(s) \|\theta^\epsilon_n(s)-\eta(s)\|_{L^2(\mathcal{O})}^2\,ds+2E\int_0^T e^{-r(s)} \langle P_nF(\theta^\epsilon_n(s))-P_nF(\eta(s)),\theta^\epsilon_n(s)-\eta(s)\rangle\,ds\\
\nonumber&+E\int_0^T e^{-r(s)}r'(s) \|\eta(s)\|_{L^2(\mathcal{O})}^2\,ds-2E\int_0^S e^{-r(s)}r'(s) \int_{\mathcal{O}}\theta^\epsilon_n(x,y,z,s)\eta(x,y,z,s)\,dxdydzds\\
\nonumber&+2E\int_0^T e^{-r(s)} \langle P_nF(\eta(s)),\theta^\epsilon_n(s)-\eta(s)\rangle\,ds+2E\int_0^T e^{-r(s)} \langle P_nF(\theta^\epsilon_n(s)),\eta(s)\rangle\,ds\\
\nonumber&+\epsilon E\int_0^Te^{-r(s)}\|\sigma_n(s,\theta^\epsilon_n(s))-\sigma_n(s,\eta(s))\|_{\mathcal{L}_2(U_0,H_2)}^2\,ds+2\epsilon E\int_0^Te^{-r(s)}(\sigma_n(s,\theta^\epsilon_n(s)),\sigma_n(s,\eta(s)))_{\mathcal{L}_2(U_0,H_2)}\,ds\\
\nonumber&-\epsilon E\int_0^Te^{-r(s)}\|\sigma_n(s,\eta(s))\|_{\mathcal{L}_2(U_0,H_2)}^2\,ds\\
\nonumber\leq&E\int_0^T e^{-r(s)}r'(s) \|\eta(s)\|_{L^2(\mathcal{O})}^2\,ds-2E\int_0^T e^{-r(s)}r'(s) \int_{\mathcal{O}}\theta^\epsilon_n(x,y,z,s)\eta(x,y,z,s)\,dxdydzds\\
\nonumber&+2E\int_0^T e^{-r(s)} \langle P_nF(\eta(s)),\theta^\epsilon_n(s)-\eta(s)\rangle\,ds+2E\int_0^T e^{-r(s)} \langle P_nF(\theta^\epsilon_n(s)),\eta(s)\rangle\,ds\\
&+2\epsilon E\int_0^Te^{-r(s)}(\sigma_n(s,\theta^\epsilon_n(s)),\sigma_n(s,\eta(s)))_{\mathcal{L}_2(U_0,H_2)}\,ds-\epsilon E\int_0^Te^{-r(s)}\|\sigma_n(s,\eta(s))\|_{\mathcal{L}_2(U_0,H_2)}^2\,ds.
\end{align}
By lower semi-continuity property of weak convergence and inequality \eqref{3.24} as well as the Lebesgue dominated convergence theorem, we obtain
\begin{align}\label{3.29}
\nonumber&E\left(\|\theta^\epsilon(T)\|_{L^2(\mathcal{O})}^2e^{-r(T)}\right)-E\|\theta_0\|_{L^2(\mathcal{O})}^2\\
\nonumber\leq&\liminf_{n\rightarrow+\infty}\left[E\left(\|\theta^\epsilon_n(T)\|_{L^2(\mathcal{O})}^2e^{-r(T)}\right)-E\|P_n\theta_0\|_{L^2(\mathcal{O})}^2\right]\\
\nonumber\leq&\liminf_{n\rightarrow+\infty}\left[E\int_0^T e^{-r(s)}r'(s) \|\eta(s)\|_{L^2(\mathcal{O})}^2\,ds-2E\int_0^T e^{-r(s)}r'(s) \int_{\mathcal{O}}\theta^\epsilon_n(x,y,z,s)\eta(x,y,z,s)\,dxdydzds\right.
\end{align}
\begin{align}
\nonumber&\left.+2E\int_0^T e^{-r(s)} \langle P_nF(\eta(s)),\theta^\epsilon_n(s)-\eta(s)\rangle\,ds+2E\int_0^T e^{-r(s)} \langle P_nF(\theta^\epsilon_n(s)),\eta(s)\rangle\,ds\right.\\
\nonumber&\left.+2\epsilon E\int_0^Te^{-r(s)}(\sigma_n(s,\theta^\epsilon_n(s)),\sigma_n(s,\eta(s)))_{\mathcal{L}_2(U_0,H_2)}\,ds-\epsilon E\int_0^Te^{-r(s)}\|\sigma_n(s,\eta(s))\|_{\mathcal{L}_2(U_0,H_2)}^2\,ds\right]\\
\nonumber\leq&E\int_0^T e^{-r(s)}r'(s) \|\eta(s)\|_{L^2(\mathcal{O})}^2\,ds-2E\int_0^T e^{-r(s)}r'(s) \int_{\mathcal{O}}\bar{\theta}^\epsilon(x,y,z,s)\eta(x,y,z,s)\,dxdydzds\\
\nonumber&+2E\int_0^T e^{-r(s)} \langle F(\eta(s)),\bar{\theta}^\epsilon(s)-\eta(s)\rangle\,ds+2E\int_0^T e^{-r(s)} \langle F(s),\eta(s)\rangle\,ds\\
&+2\epsilon E\int_0^Te^{-r(s)}(S(s),\sigma(s,\eta(s)))_{\mathcal{L}_2(U_0,H_2)}\,ds-\epsilon E\int_0^Te^{-r(s)}\|\sigma(s,\eta(s))\|_{\mathcal{L}_2(U_0,H_2)}^2\,ds.
\end{align}
From inequality \eqref{3.26}  and inequality \eqref{3.29}, we conclude
\begin{align}\label{3.30}
\nonumber&-E\int_0^T e^{-r(s)}r'(s) \|\bar{\theta}^\epsilon(s)-\eta(s)\|_{L^2(\mathcal{O})}^2\,ds+2E\int_0^T e^{-r(s)} \langle F(s)-F(\eta(s)),\bar{\theta}^\epsilon(s)-\eta(s)\rangle\,ds\\
&+\epsilon E\int_0^Te^{-r(s)}\|S(s)-\sigma(s,\eta(s))\|_{\mathcal{L}_2(U_0,H_2)}^2\,ds\leq0
\end{align}
for every $\mathcal{F}$-progressively measurable process $\eta\in L^2(\Omega\times(0,T);V_2)\cap L^4(\Omega\times (0,T);L^3(\mathcal{O}))\cap L^4(\Omega;L^{\infty}(0,T;H_2)).$ In particular, taking $\eta=\bar{\theta}^\epsilon,$ we obtain $S(s,\omega)=\sigma(s,\bar{\theta}^\epsilon(s,\omega))$ for $dt\otimes d\mathcal{P}$ -a.e. $(s,\omega)\in (0,T)\times\Omega.$

Let $\eta=\bar{\theta}^\epsilon-\lambda \psi$ for any $\psi\in L^{\infty}(\Omega\times (0,T); V_2)$ and any $\lambda>0,$ we obtain
\begin{align}\label{3.31}
-\lambda E\int_0^T e^{-r(s)}r'(s) \|\psi(s)\|_{L^2(\mathcal{O})}^2\,ds+2E\int_0^T e^{-r(s)} \langle F(s)-F(\theta^\epsilon(s)-\lambda\psi(s)),\psi(s)\rangle\,ds\leq0.
\end{align}
Since for any $\lambda>0,$ we have
\begin{align}\label{3.32}
\left|\langle F(\theta^\epsilon(s))-F(\theta^\epsilon(s)-\lambda\psi(s)),\psi(s)\rangle\right|\leq C\lambda\left(\|\psi(s)\|^2+\|\theta^\epsilon(s)\|^2\|\psi(s)\|_{L^2(\mathcal{O})}^2\right),\,\,\,\forall\,\,\,s\in[0,T].
\end{align}
It follows from the Lebesgue Dominated convergence theorem, we have as $\lambda\downarrow 0,$
\begin{align}\label{3.33}
E\int_0^T e^{-r(s)} \langle F(s)-F(\theta^\epsilon(s)-\lambda\psi(s)),\psi(s)\rangle\,ds\rightarrow E\int_0^T e^{-r(s)} \langle F(s)-F(\theta^\epsilon(s)),\psi(s)\rangle\,ds
\end{align}
Let $\lambda\downarrow 0$ on both hand sides of inequality \eqref{3.31}, yields
\begin{align}\label{3.34}
E\int_0^T e^{-r(s)} \langle F(s)-F(\theta^\epsilon(s)),\psi(s)\rangle\,ds\leq0.
\end{align}
Since $\psi$ is arbitrary, this implies that the process $F(s)=F(\theta^\epsilon(s))\in L^2(\Omega\times(0,T);V_2').$ Therefore, problem \eqref{2.1.2} possesses a solution in the sense of Definition \ref{2.3.3}. Moreover, the estimates stated in Theorem \ref{3.1} can be concluded from the lower semi-continuity property of weak convergence.

{\bf Step 2.} General case: $\theta_0\in L^2(\Omega, H_2).$
Taking any sequence $\theta_n(0)\in L^4(\Omega, H_2)$ such that $E\left(\|\theta_n(0)-\theta_0\|_{L^2(\mathcal{O})}^2\right)\rightarrow 0.$ Let $(v^\epsilon_n(t),p^\epsilon_{ns}(t),\theta^\epsilon_n(t)),$ $t\geq 0$ be the solution of the following equation:
\begin{equation}\label{3.35}
\begin{cases}
&\nabla p^\epsilon_{ns}+f{v^\epsilon_n}^\bot+L_1v^\epsilon_n=\int_{-h}^z\nabla \theta^\epsilon_n(x,y,\zeta,t)\,d\zeta,\\
&\int_{-h}^{0}\nabla\cdot v^\epsilon_n(x,y,\zeta,t)\,d\zeta=0,\\
&d\theta^\epsilon_n+B(v^\epsilon_n,\theta^\epsilon_n)\,dt+A_2(\theta^\epsilon_n-\theta^*)\,dt-K_h\Delta \theta^*\,dt=g(x,y,z,t)\,dt+\sqrt{\epsilon}\sigma(t,\theta^\epsilon_n(t))dW(t),\\
&\left.A_{\nu}\frac{\partial v^\epsilon_n}{\partial z}\right|_{\Gamma_u}=\mu, \left.\frac{\partial v^\epsilon_n}{\partial z}\right|_{\Gamma_b}=0, \left.v^\epsilon_n\cdot\vec{n}\right|_{\Gamma_l}=0,\left.\frac{\partial v^\epsilon_n}{\partial\vec{n}}\times\vec{n}\right|_{\Gamma_l}=0,\\
&\theta^\epsilon_n(x,y,z,0)=\theta_n(0).
\end{cases}
\end{equation}
The existence of $(v^\epsilon_n(t),p^\epsilon_{ns}(t),\theta^\epsilon_n(t))$ of problem \eqref{3.35} can be established by step 1. Moreover, from the proof of inequality \eqref{3.8} and Lemma \ref{2.2.2}, we can obtain
\begin{align}\label{3.36}
\nonumber &\sup_n\left(E\left(\sup_{0\leq r\leq T}\|\theta^\epsilon_n(r)\|_{L^2(\mathcal{O})}^2\right)+E\left(\int_0^T\|\theta^\epsilon_n(s)\|^2\,ds\right)+E\left(\int_0^T\|p^\epsilon_{ns}(s)\|_{H^1(\mathcal{M})}^2\,ds\right)+E\left(\int_0^T\|v^\epsilon_n(s)\|_{H^1(\mathcal{O})}^2\,ds\right)\right)\\
\nonumber\leq&C\left(\sup_n (E\|\theta_n(0)\|_{L^2(\mathcal{O})}^2)+\int_0^T\|g(s)\|_{L^2(\mathcal{O})}^2\,ds+T(\|\theta^*\|_{L^2(\mathcal{M})}^2+\|\mu\|_{H^1(\mathcal{M})}^2+1)\right)e^{CKT}\\
=&D_2(T)<+\infty,
\end{align}
which implies that there exist a subsequence (still use the same notation) of $\{(v^\epsilon_n(t),p^\epsilon_{ns}(t),\theta^\epsilon_n(t)):n\geq 1\}$ and a process
$v^\epsilon\in L^2(\Omega\times (0,T);V_1),$ $p^\epsilon_s\in L^2(\Omega\times (0,T);H^1(\mathcal{M})),$ $\theta^\epsilon\in L^2(\Omega\times(0,T);V_2)\cap L^2(\Omega;L^{\infty}(0,T;H_2))$ such that the following hold:
\begin{align*}
&\theta^\epsilon_n\rightharpoonup \theta^\epsilon\;\textit{weakly in}\; L^2(\Omega\times(0,T);V_2),\\
&\theta^\epsilon_n\rightharpoonup \theta^\epsilon\;\textit{weakly star in}\; L^2(\Omega;L^{\infty}(0,T;H_2)),\\
&v^\epsilon_n\rightharpoonup v^\epsilon\;\textit{weakly in}\; L^2(\Omega\times (0,T);V_1),\\
&p^\epsilon_{ns}\rightharpoonup p^\epsilon_s\;\textit{weakly in}\; L^2(\Omega\times (0,T);H^1(\mathcal{M})).
\end{align*}
Based on the above weak convergence property, we can easily prove that
\begin{align*}
&\int_\mathcal{O} \nabla p^\epsilon_s(x,y,t)\cdot\phi\,dxdydz-\int_{\mathcal{O}}\left(\int_{-h}^z
\nabla \theta^\epsilon(x,y,\zeta,t)d\zeta\right)\cdot\phi\,dxdydz+\int_{\mathcal{O}}f{v^\epsilon}^\bot\cdot\phi\,dxdydz+\langle A_1v^\epsilon, \phi\rangle\\
&=\int_{\Gamma_u} \kappa\mu\cdot\phi\, dxdy,\,\,\,\,d\mathcal{P}\otimes dt-a.e.\,\,\,\textit{on}\,\,\,\Omega\times(0,T),
\end{align*}
for any $\phi\in H^1(\mathcal{O}).$

Next, we want to prove that $\theta^\epsilon_n$ also converges to $\theta^\epsilon$ in probability in $L^{\infty}(0,T; H_2)\cap L^2(0,T;V_2).$ For any fixed $R>0,$ define 
\begin{align*}
\tau_R^n=\inf\{t>0:\|\theta^\epsilon_n(t)\|_{L^2(\mathcal{O})}^2+\int_0^t\|\theta^\epsilon_n(s)\|^2\,ds>R\}.
\end{align*}
then $\tau_R^n$ is a stopping time and the following result hold:
\begin{align}\label{3.37}
\mathcal{P}(\tau_R^n\leq T)\leq \mathcal{P}(\sup_{0\leq t\leq T}\|\theta^\epsilon_n(t)\|_{L^2(\mathcal{O})}^2+\int_0^T\|\theta^\epsilon_n(s)\|^2\,ds>R)\leq \frac{D_2(T)}{R}
\end{align}
for all $n\geq 1.$

Put $(v, p_{s},\theta)=(v^\epsilon_n, p^\epsilon_{ns},\theta^\epsilon_n)-(v^\epsilon_m, p^\epsilon_{ms},\theta^\epsilon_m)$ with $\theta(0)=\theta_n(0)-\theta_m(0)$ for any $m,$ $n\geq 1,$ then the following conclusion holds:
\begin{align}\label{3.38}
\|v(t)\|_{H^{\gamma+1}(\mathcal{O})}^2+\|p_{s}(t)\|_{H^{\gamma}(\mathcal{M})}^2\leq \mathcal{K}_1\|\theta(t)\|_{H^\gamma(\mathcal{O})}^2
\end{align}
for $\gamma=0,1.$

We apply the It\^{o}'s formula to the process $\|\theta(s)\|_{L^2(\mathcal{O})}^2$ and use inequality \eqref{3.38} as well as Lemma \ref{2.2.3}, yield
\begin{align*}
&\|\theta(t\land\tau_R^n\land\tau_R^m)\|_{L^2(\mathcal{O})}^2+2\int_0^{t\land\tau_R^n\land\tau_R^m}\|\theta(s)\|^2\,ds\\
=&\|\theta(0)\|_{L^2(\mathcal{O})}^2-2\int_0^{t\land\tau_R^n\land\tau_R^m}b(v(s),\theta^\epsilon_n(s),\theta(s))\,ds+\epsilon\int_0^{t\land\tau_R^n\land\tau_R^m}\|\sigma(s,\theta^\epsilon_n(s))-\sigma(s,\theta^\epsilon_m(s))\|_{\mathcal{L}_2(U_0,H_2)}^2\,ds\\
&+2\sqrt{\epsilon}\int_0^{t\land\tau_R^n\land\tau_R^m}\int_{\mathcal{O}}\theta(s)(\sigma(s,\theta^\epsilon_n(s))-\sigma(s,\theta^\epsilon_m(s)))\,dW(s)dxdydz\\
\leq&\|\theta(0)\|_{L^2(\mathcal{O})}^2+2\mathcal{K}\int_0^{t\land\tau_R^n\land\tau_R^m}\|v(s)\|_{H^1(\mathcal{O})}^{\frac{1}{2}}\|v(s)\|_{H^2(\mathcal{O})}^{\frac{1}{2}}\|\theta(s)\|_{L^2(\mathcal{O})}^{\frac{1}{2}}\|\theta(s)\|^{\frac{1}{2}}\|\theta_n(s)\|\,ds+\epsilon L\int_0^{t\land\tau_R^n\land\tau_R^m}\|\theta(s))\|_{L^2(\mathcal{O})}^2\,ds\\
&+2\sqrt{\epsilon}\int_0^{t\land\tau_R^n\land\tau_R^m}\int_{\mathcal{O}}\theta(s)(\sigma(s,\theta^\epsilon_n(s))-\sigma(s,\theta^\epsilon_m(s)))\,dW(s)dxdydz\\
\leq&\|\theta(0)\|_{L^2(\mathcal{O})}^2+2\mathcal{K}\sqrt{\mathcal{K}_1}\int_0^{t\land\tau_R^n\land\tau_R^m}\|\theta(s)\|_{L^2(\mathcal{O})}\|\theta_n(s)\|\|\theta(s)\|\,ds+\epsilon L\int_0^{t\land\tau_R^n\land\tau_R^m}\|\theta(s))\|_{L^2(\mathcal{O})}^2\,ds\\
&+2\sqrt{\epsilon}\int_0^{t\land\tau_R^n\land\tau_R^m}\int_{\mathcal{O}}\theta(s)(\sigma(s,\theta^\epsilon_n(s))-\sigma(s,\theta^\epsilon_m(s)))\,dW(s)dxdydz\\
\leq&\|\theta(0)\|_{L^2(\mathcal{O})}^2+\int_0^{t\land\tau_R^n\land\tau_R^m}\|\theta(s)\|_{L^2(\mathcal{O})}^2(\mathcal{K}^2\mathcal{K}_1\|\theta_n(s)\|^2+\epsilon L)\,ds+\int_0^{t\land\tau_R^n\land\tau_R^m}\|\theta(s)\|^2\,ds\\
&+2\sqrt{\epsilon}\int_0^{t\land\tau_R^n\land\tau_R^m}\int_{\mathcal{O}}\theta(s)(\sigma(s,\theta^\epsilon_n(s))-\sigma(s,\theta^\epsilon_m(s)))\,dW(s)dxdydz
\end{align*}
for any $m,$ $n\geq 1,$ which implies that
\begin{align*}
&\|\theta(t\land\tau_R^n\land\tau_R^m)\|_{L^2(\mathcal{O})}^2+\int_0^{t\land\tau_R^n\land\tau_R^m}\|\theta(s)\|^2\,ds\\
\leq&\|\theta(0)\|_{L^2(\mathcal{O})}^2+\int_0^{t\land\tau_R^n\land\tau_R^m}\|\theta(s)\|_{L^2(\mathcal{O})}^2(\mathcal{K}^2\mathcal{K}_1\|\theta_n(s)\|^2+\epsilon L)\,ds\\
&+2\sqrt{\epsilon}\int_0^{t\land\tau_R^n\land\tau_R^m}\int_{\mathcal{O}}\theta(s)(\sigma(s,\theta^\epsilon_n(s))-\sigma(s,\theta^\epsilon_m(s)))\,dW(s)dxdydz
\end{align*}
for any $m,$ $n\geq 1,$ it follows from Lemma \ref{2.2.4}  that
\begin{align*}
&\|\theta(t\land\tau_R^n\land\tau_R^m)\|_{L^2(\mathcal{O})}^2+\int_0^{t\land\tau_R^n\land\tau_R^m}\|\theta(s)\|^2\,ds\\
\leq&e^{\mathcal{K}^2\mathcal{K}_1\int_0^{t\land\tau_R^n\land\tau_R^m}\|\theta_n(s)\|^2\,ds+\epsilon LT}\left(\|\theta_n(0)-\theta_m(0)\|_{L^2(\mathcal{O})}^2\right.\\
&\left.+2\sqrt{\epsilon}\int_0^{t\land\tau_R^n\land\tau_R^m}\int_{\mathcal{O}}\theta(s)(\sigma(s,\theta^\epsilon_n(s))-\sigma(s,\theta^\epsilon_m(s)))\,dW(s)dxdydz\right)\\
\leq&e^{\mathcal{K}^2\mathcal{K}_1R+\epsilon LT}\left(\|\theta_n(0)-\theta_m(0)\|_{L^2(\mathcal{O})}^2+2\sqrt{\epsilon}\int_0^{t\land\tau_R^n\land\tau_R^m}\int_{\mathcal{O}}\theta(s)(\sigma(s,\theta^\epsilon_n(s))-\sigma(s,\theta^\epsilon_m(s)))\,dW(s)dxdydz\right).
\end{align*}
By the Burkholder-Davis-Gundy inequality, we obtain 
\begin{align*}
\nonumber&2\sqrt{\epsilon}E\left(\sup_{0\leq r\leq T\land\tau^n_R\land\tau_R^m}\left|\int_0^r\int_{\mathcal{O}}\theta(s)(\sigma(s,\theta^\epsilon_n(s))-\sigma(s,\theta^\epsilon_m(s)))\,dW(s)dxdydz\right|\right)\\
\nonumber\leq &C_1\sqrt{\epsilon}E\left(\left(\int_0^{T\land\tau^n_R\land\tau_R^m}\|\theta(s)\|_{L^2(\mathcal{O})}^2\|\sigma(s,\theta^\epsilon_n(s))-\sigma(s,\theta^\epsilon_m(s))\|_{\mathcal{L}_2(U_0,H_2)}^2\,ds\right)^{\frac{1}{2}}\right)\\
\nonumber\leq &C_1\sqrt{\epsilon}E\left(\sup_{0\leq r\leq T\land\tau^n_R\land\tau_R^m}\|\theta(r)\|_{L^2(\mathcal{O})}\left(\int_0^{T\land\tau^n_R\land\tau_R^m}\|\sigma(s,\theta^\epsilon_n(s))-\sigma(s,\theta^\epsilon_m(s))\|_{\mathcal{L}_2(U_0,H_2)}^2\,ds\right)^{\frac{1}{2}}\right)\\
\leq&C_1\sqrt{\epsilon L}\left(E\sup_{0\leq r\leq T\land\tau^n_R\land\tau_R^m}\|\theta(r)\|_{L^2(\mathcal{O})}^2\right)^{\frac{1}{2}}\left(E\int_0^{T\land\tau^n_R\land\tau_R^m}\|\theta(s))\|_{L^2(\mathcal{O})}^2\,ds\right)^{\frac{1}{2}}.
\end{align*}
Therefore, we conclude the following result:
\begin{align*}
\nonumber&E\left(\sup_{0\leq t\leq T}\|\theta_n(t\land\tau_R^n\land\tau_R^m)-\theta_m(t\land\tau_R^n\land\tau_R^m)\|_{L^2(\mathcal{O})}^2\right)+E\int_0^{T\land\tau_R^n\land\tau_R^m}\|\theta_n(s)-\theta_m(s)\|^2\,ds\\
\leq&e^{\mathcal{K}^2\mathcal{K}_1R+\epsilon LT}\left(\|\theta_n(0)-\theta_m(0)\|_{L^2(\mathcal{O})}^2+C_1\sqrt{\epsilon L}\left(E\sup_{0\leq r\leq T\land\tau^n_R\land\tau_R^m}\|\theta(r)\|_{L^2(\mathcal{O})}^2\right)^{\frac{1}{2}}\left(E\int_0^{T\land\tau^n_R\land\tau_R^m}\|\theta(s))\|_{L^2(\mathcal{O})}^2\,ds\right)^{\frac{1}{2}}\right)\\
\leq&e^{\mathcal{K}^2\mathcal{K}_1R+\epsilon LT}E\left(\|\theta_n(0)-\theta_m(0)\|_{L^2(\mathcal{O})}^2\right)+\frac{1}{2}E\left(\sup_{0\leq r\leq T\land\tau^n_R\land\tau_R^m}\|\theta(r)\|_{L^2(\mathcal{O})}^2\right)\\
&+\frac{1}{2}e^{2\mathcal{K}^2\mathcal{K}_1R+2\epsilon LT}C_1^2\epsilon L\left(E\int_0^{T\land\tau^n_R\land\tau_R^m}\|\theta(s))\|_{L^2(\mathcal{O})}^2\,ds\right).
\end{align*}
It follows from Lemma \ref{2.2.4} that
\begin{align}\label{3.39}
\nonumber&E\left(\sup_{0\leq t\leq T}\|\theta_n(t\land\tau_R^n\land\tau_R^m)-\theta_m(t\land\tau_R^n\land\tau_R^m)\|_{L^2(\mathcal{O})}^2\right)+E\int_0^{T\land\tau_R^n\land\tau_R^m}\|\theta_n(s)-\theta_m(s)\|^2\,ds\\
\leq&D_3(R,T)E\left(\|\theta_n(0)-\theta_m(0)\|_{L^2(\mathcal{O})}^2\right),
\end{align}
where
\begin{align*}
D_3(R,T)=2e^{\mathcal{K}^2\mathcal{K}_1R+\epsilon LT+\frac{1}{4}C_1^2e^{2\mathcal{K}^2\mathcal{K}_1R+2\epsilon LT}}.
\end{align*}

For $\delta>0$ and any $R>0,$ we obtain
\begin{align}\label{3.40}
\nonumber&\mathcal{P}\left(\sup_{0\leq t\leq T}\|\theta_n(t)-\theta_m(t)\|_{L^2(\mathcal{O})}^2>\delta\right)\\
\leq&\mathcal{P}(\tau_R^m\leq T)+\mathcal{P}(\tau_R^n\leq T)+\mathcal{P}\left(\sup_{0\leq t\leq T}\|\theta_n(t\land\tau_R^n\land\tau_R^m)-\theta_m(t\land\tau_R^n\land\tau_R^m)\|_{L^2(\mathcal{O})}^2>\delta\right)
\end{align}
Given an arbitrary small constant $\epsilon>0,$ we deduce from inequality \eqref{3.37} that there exists a positive constant $R_0$ such that for any $R\geq R_0$ and any $m,$ $n\geq 1,$
\begin{align}\label{3.41}
\mathcal{P}(\tau_R^n\leq T)\leq \frac{\epsilon}{3}
\end{align}
and
\begin{align}\label{3.42}
\mathcal{P}(\tau_R^m\leq T)\leq \frac{\epsilon}{3}.
\end{align}
Moreover, it follows from inequality \eqref{3.39} that there exists a natural number $N_0$ such that for all $m,$ $n\geq N_0,$
\begin{align}\label{3.43}
\mathcal{P}\left(\sup_{0\leq t\leq T}\|\theta_n(t\land\tau_R^n\land\tau_R^m)-\theta_m(t\land\tau_R^n\land\tau_R^m)\|_{L^2(\mathcal{O})}^2>\delta\right)\leq \frac{\epsilon}{3}.
\end{align}
Therefore, we conclude that for all $m,$ $n\geq N_0,$
\begin{align}\label{3.44}
\mathcal{P}\left(\sup_{0\leq t\leq T}\|\theta_n(t)-\theta_m(t)\|_{L^2(\mathcal{O})}^2>\delta\right)\leq\epsilon.
\end{align}
Similarly, we can obtain 
\begin{align}\label{3.45}
\mathcal{P}\left(\int_0^T\|\theta_n(s)-\theta_m(s)\|^2\,ds>\delta\right)\leq\epsilon
\end{align}
for all $m,$ $n\geq N_1\geq N_0.$ These inequalities \eqref{3.44}-\eqref{3.45} imply that $(v_n,\theta_n)$ converges to $(v,\theta)$ in probability in $L^{\infty}(0,T; H^1(\mathcal{O})\times H_2)\cap L^2(0,T; H^2(\mathcal{O})\times V_2).$ Finally, we would like to prove that $(v,\theta)$ satisfies equality \eqref{2.3.4}. To this end, it suffices to prove that for any $\psi\in V_2,$
\begin{align}\label{3.46}
\nonumber &\int_{\mathcal{O}}\theta(x,y,z,t)\psi(x,y,z)\,dxdydz=-\int_0^t \langle A_2(\theta(s)-\theta^*),\psi\rangle\,ds-\int_0^t b(v(s),\theta(s),\psi)\,ds
+\langle T_0+K_ht\Delta \theta^*,\psi\rangle\\
&+\int_0^t \langle g(s),\psi\rangle
\,ds+\int_0^t\langle\sigma(s,\theta(s))\,dW(s),\psi\rangle.
\end{align}
Since for each $n\geq 1,$ we have
\begin{align}\label{3.47}
\nonumber &\int_{\mathcal{O}}\theta_n(x,y,z,t)\psi(x,y,z)\,dxdydz=-\int_0^t b(v_n(s),\theta_n(s),\psi)\,ds
+\langle \theta_n(0)+K_ht\Delta \theta^*,\psi\rangle\\
&+\int_0^t \langle g(s),\psi\rangle
\,ds+\int_0^t\langle\sigma(s,\theta_n(s))\,dW(s),\psi\rangle-\int_0^t \langle A_2(\theta_n(s)-\theta^*),\psi\rangle\,ds.
\end{align}
Let $n\rightarrow+\infty,$ thanks to the convergence in probability and also the weak convergence, from the dominated convergence theorem, we conclude that each term in \eqref{3.47} tends to the corresponding term in \eqref{3.46}. Hence, the existence proof of weak solutions for problem \eqref{2.1.2} is completed.

In what follows, we will prove that the solution $(v^\epsilon, p^\epsilon_s,\theta^\epsilon)$ of problem \eqref{2.1.2} is unique. For this purpose, we will use the ‘‘Schmalfuss trick’’ as in \cite{sb1}.

Assume that $(v^\epsilon_1, p^\epsilon_{1s},\theta^\epsilon_1)$ and $(v^\epsilon_1, p^\epsilon_{1s},\theta^\epsilon_1)$ are two solutions to problem \eqref{2.2.1}, put $(v, p_{s},\theta)=(v^\epsilon_1, p^\epsilon_{1s},\theta^\epsilon_1)-(v^\epsilon_2, p^\epsilon_{2s},\theta^\epsilon_2),$ then $\theta(0)=0$ and the following conclusion holds:
\begin{align}\label{3.48}
\|v(t)\|_{H^{\gamma+1}(\mathcal{O})}^2+\|p_{s}(t)\|_{H^{\gamma}(\mathcal{M})}^2\leq \mathcal{K}_1\|\theta(t)\|_{H^\gamma(\mathcal{O})}^2
\end{align}
for $\gamma=0,1.$

Define an auxiliary process $\phi$ by
\begin{align*}
\phi(t):=\exp\{-L_2\int_0^t\|\theta^\epsilon_2(s)\|^2\,ds\},\,\,\,\,\,t\geq 0,
\end{align*}
where $L_2\geq \mathcal{K}^2\mathcal{K}_1.$

Let us introduce the stopping time
\begin{align*}
\tau_R=\inf\{t>0:\|\theta^\epsilon_1(t)\|_{L^2(\mathcal{O})}^2\geq R\}\land \inf\{t>0:\|\theta^\epsilon_2(t)\|_{L^2(\mathcal{O})}^2\geq R\}\land T,
\end{align*}
then we apply the It\^{o}'s formula to the process $\|\theta(s)\|_{L^2(\mathcal{O})}^2\phi(s)$ and use inequality \eqref{3.48} as well as assumption ($A_3$), Lemma \ref{2.2.3}, yields
\begin{align*}
&\phi(t\land\tau_R)\|\theta(t\land\tau_R)\|_{L^2(\mathcal{O})}^2+L_2\int_0^{t\land\tau_R}\phi(s)\|\theta^\epsilon_2(s)\|^2\|\theta(s)\|_{L^2(\mathcal{O})}^2\,ds+2\int_0^{t\land\tau_R}\phi(s)\|\theta(s)\|^2\,ds\\
=&-2\int_0^{t\land\tau_R}b(v(s),\theta^\epsilon_2(s),\theta(s))\phi(s)\,ds+\epsilon\int_0^{t\land\tau_R}\phi(s)\|\sigma(s,\theta^\epsilon_1(s))-\sigma(s,\theta^\epsilon_2(s))\|_{\mathcal{L}_2(U_0,H_2)}^2\,ds\\
&+2\sqrt{\epsilon}\int_0^{t\land\tau_R}\int_{\mathcal{O}}\phi(s)\theta(s)(\sigma(s,\theta^\epsilon_1(s))-\sigma(s,\theta^\epsilon_2(s)))\,dW(s)dxdydz\\
\leq&2\mathcal{K}\int_0^{t\land\tau_R}\|v(s)\|_{H^1(\mathcal{O})}^{\frac{1}{2}}\|v(s)\|_{H^2(\mathcal{O})}^{\frac{1}{2}}\|\theta(s)\|_{L^2(\mathcal{O})}^{\frac{1}{2}}\|\theta(s)\|^{\frac{1}{2}}\|\theta^\epsilon_2(s)\|\phi(s)\,ds+\epsilon L\int_0^{t\land\tau_R}\phi(s)\|\theta(s))\|_{L^2(\mathcal{O})}^2\,ds\\
&+2\sqrt{\epsilon}\int_0^{t\land\tau_R}\int_{\mathcal{O}}\phi(s)\theta(s)(\sigma(s,\theta^\epsilon_1(s))-\sigma(s,\theta^\epsilon_2(s)))\,dW(s)dxdydz\\
\leq&2\mathcal{K}\sqrt{\mathcal{K}_1}\int_0^{t\land\tau_R}\phi(s)\|\theta(s)\|_{L^2(\mathcal{O})}\|\theta^\epsilon_2(s)\|\|\theta(s)\|\,ds+\epsilon L\int_0^{t\land\tau_R}\phi(s)\|\theta(s))\|_{L^2(\mathcal{O})}^2\,ds\\
&+2\sqrt{\epsilon}\int_0^{t\land\tau_R}\int_{\mathcal{O}}\phi(s)\theta(s)(\sigma(s,\theta^\epsilon_1(s))-\sigma(s,\theta^\epsilon_2(s)))\,dW(s)dxdydz\\
\leq&\mathcal{K}^2\mathcal{K}_1\int_0^{t\land\tau_R}\phi(s)\|\theta(s)\|_{L^2(\mathcal{O})}^2\|\theta^\epsilon_2(s)\|^2\,ds+\int_0^{t\land\tau_R}\phi(s)\|\theta(s)\|^2\,ds+\epsilon L\int_0^{t\land\tau_R}\phi(s)\|\theta(s))\|_{L^2(\mathcal{O})}^2\,ds\\
&+2\sqrt{\epsilon}\int_0^{t\land\tau_R}\int_{\mathcal{O}}\phi(s)\theta(s)(\sigma(s,\theta^\epsilon_1(s))-\sigma(s,\theta^\epsilon_2(s)))\,dW(s)dxdydz.
\end{align*}
Taking the expectation of both hand sides of the above inequality, we obtain
\begin{align*}
&E\left(\phi(t\land\tau_R)\|\theta(t\land\tau_R)\|_{L^2(\mathcal{O})}^2\right)+\int_0^{t\land\tau_R}E\left(\phi(s)\|\theta(s)\|^2\right)\,ds\\
\leq&\epsilon L\int_0^{t\land\tau_R}E\left(\phi(s)\|\theta(s))\|_{L^2(\mathcal{O})}^2\right)\,ds.
\end{align*}
It follows from Lemma \ref{2.2.4} that
\begin{align*}
&E\left(\phi(t\land\tau_R)\|\theta(t\land\tau_R)\|_{L^2(\mathcal{O})}^2\right)\leq 0.
\end{align*}
Therefore, the solution $(v^\epsilon, p^\epsilon_s,\theta^\epsilon)$ of problem \eqref{2.1.2} is unique.\\
\qed\hfill

\section{Large deviation principle}
\def\theequation{4.\arabic{equation}}\makeatother
\setcounter{equation}{0}
In this section, we will establish the large deviation principle for problem \eqref{2.1.2}.
\subsection{The large deviations result}
Consider the following three dimensional stochastic planetary geostrophic equations of large-scale ocean circulation:
\begin{equation}\label{4.2.1}
\begin{cases}
&\nabla p^\epsilon_s+f{v^\epsilon}^\bot+L_1v^\epsilon=\int_{-h}^z\nabla \theta^\epsilon(x,y,\zeta,t)\,d\zeta,\\
&\int_{-h}^{0}\nabla\cdot v^\epsilon(x,y,\zeta,t)\,d\zeta=0,\\
&d\theta^\epsilon+B(v^\epsilon, \theta^\epsilon)\,dt+A_2(\theta^\epsilon-\theta^*)\,dt-K_h\Delta \theta^*\,dt=g(x,y,z,t)\,dt+\sqrt{\epsilon}\sigma(t,\theta^\epsilon(t))dW(t),\\
&\left.A_{\nu}\frac{\partial v^\epsilon}{\partial z}\right|_{\Gamma_u}=\mu, \left.\frac{\partial v^\epsilon}{\partial z}\right|_{\Gamma_b}=0, \left.v^\epsilon\cdot\vec{n}\right|_{\Gamma_l}=0,\left.\frac{\partial v^\epsilon}{\partial\vec{n}}\times\vec{n}\right|_{\Gamma_l}=0,\\
&\theta^\epsilon(x,y,z,0)=\theta_0(x,y,z)
\end{cases}
\end{equation}
with $\int_0^T\|g(s)\|_{L^2(\mathcal{O})}^2\,ds<+\infty$ and $\epsilon\in (0,1).$ From Theorem \ref{3.1}, we conclude that there exists a weak solution $(v,p_s,\theta)$ of problem \eqref{2.2.1} with values in $\mathcal{C}([0,T];V_1\times H^1(M)\times H_2)\cap L^2(0,T;H^2(\mathcal{O})\times H^2(M)\times V_2)$ and it is pathwise unique. It follows that (see \cite{ba}) there exists a Borel-measurable function $\Phi^\epsilon:\mathcal{C}([0,T];U)\rightarrow \mathcal{C}([0,T]; H_2)\cap L^2(0,T; V_2)$ such that $\theta^\epsilon(\cdot)=\Phi^\epsilon(W(\cdot))$ a.s.  The aim of this section is to prove the large deviation principle for $\theta^\epsilon.$ The following Lemmas show that the family $\{\Phi^\epsilon\}$ satisfies assumption (A) so that Lemma \ref{2.3.1} can be invoked to prove our main result.
\begin{lemma}\label{4.2.3}
Let the family $\{\Phi^\epsilon\}$ be defined as above. For any $\chi\in \mathcal{A}_R$ with $0<R<+\infty,$ let $\theta_\chi^\epsilon(\cdot)=\Phi^\epsilon(W(\cdot)+\frac{1}{\sqrt{\epsilon}}\int_0^{\cdot}\chi(s)\,ds).$ Then $(v_\chi^\epsilon(\cdot),p_{\chi s}^\epsilon(\cdot),\theta_\chi^\epsilon(\cdot))$ is the unique weak solution of problem
\begin{equation}\label{4.2.2}
\begin{cases}
&\nabla p^\epsilon_{\chi s}+f{v_\chi^\epsilon}^\bot+L_1v_\chi^\epsilon=\int_{-h}^z\nabla \theta_\chi^\epsilon(x,y,\zeta,t)\,d\zeta,\\
&\int_{-h}^{0}\nabla\cdot v_\chi^\epsilon(x,y,\zeta,t)\,d\zeta=0,\\
&d\theta_\chi^\epsilon+B(v_\chi^\epsilon,\theta_\chi^\epsilon)\,dt+A_2(\theta_\chi^\epsilon-\theta^*)\,dt-K_h\Delta \theta^*\,dt\\
=&g(x,y,z,t)\,dt+\sigma(t,\theta_\chi^\epsilon(t))\chi(t)\,dt+\sqrt{\epsilon}\sigma(t,\theta_\chi^\epsilon(t))dW(t),\\
&\left.A_{\nu}\frac{\partial v_\chi^\epsilon}{\partial z}\right|_{\Gamma_u}=\mu, \left.\frac{\partial v_\chi^\epsilon}{\partial z}\right|_{\Gamma_b}=0, \left.v_\chi^\epsilon\cdot\vec{n}\right|_{\Gamma_l}=0,\left.\frac{\partial v_\chi^\epsilon}{\partial\vec{n}}\times\vec{n}\right|_{\Gamma_l}=0,\\
&\theta_\chi^\epsilon(x,y,z,0)=\theta_0(x,y,z)
\end{cases}
\end{equation}
\end{lemma}
{\bf Proof.} Since $\chi\in \mathcal{A}_R,$ $\int_0^T\|\chi(s)\|_{U_0}^2\,ds<+\infty$ a.s., $\tilde{W}(\cdot)=W(\cdot)+\frac{1}{\sqrt{\epsilon}}\int_0^{\cdot} \chi(s)\,ds$ is a Wiener process with covariance form $Q$ under the probability measure
\begin{align*}
d\tilde{\mathcal{P}}_\chi^\epsilon=e^{-\frac{1}{\sqrt{\epsilon}}\int_0^T \chi(s)\,dW(s)-\frac{1}{2\epsilon}\int_0^T\|\chi(s)\|_{U_0}^2\,ds}d\mathcal{P}.
\end{align*}
A Girsanov argument can be used to complete the proof as follows: Let $(v^\epsilon,p_s^\epsilon,\theta^\epsilon)$ be the unique solution of problem \eqref{4.2.1} on $(\Omega,\mathcal{F},\tilde{\mathcal{P}}_\chi^\epsilon)$ with $\tilde{W}$ in place of $W.$ Then $(v^\epsilon,p_s^\epsilon,\theta^\epsilon)$ solves problem \eqref{4.2.2} P-a.s., and $\theta^\epsilon=\Phi^\epsilon(\tilde{W}(\cdot)).$

 If $(v_\chi^\epsilon,p_{\chi s}^\epsilon,\theta_\chi^\epsilon)$ and $(v^\epsilon,p_s^\epsilon,\theta^\epsilon)$ are two solutions of problem \eqref{4.2.2} on $(\Omega,\mathcal{F},\mathcal{P}),$ then $(v_\chi^\epsilon,p_{\chi s}^\epsilon,\theta_\chi^\epsilon)$ and $(v^\epsilon,p_s^\epsilon,\theta^\epsilon)$ will satisfy problem \eqref{4.2.1} on $(\Omega,\mathcal{F},\tilde{\mathcal{P}}_\chi^\epsilon)$ with $\tilde{W}$ in place of $W.$ Thus $(v_\chi^\epsilon,p_{\chi s}^\epsilon,\theta_\chi^\epsilon)=(v^\epsilon,p_s^\epsilon,\theta^\epsilon)$ $\tilde{\mathcal{P}}_\chi^\epsilon$-a.s. so that $(v_\chi^\epsilon,p_{\chi s}^\epsilon,\theta_\chi^\epsilon)=(v^\epsilon,p_s^\epsilon,\theta^\epsilon)$  $\mathcal{P}$-a.s. Thus uniqueness
of solutions to problem \eqref{4.2.2} is obtained.\\
\qed\hfill
\begin{lemma}(see \cite{ba})\label{4.2.4}
Let $\{\chi_n\}$ be a sequence of elements from $\mathcal{A}_R$ for some finite $R > 0.$ Let $\chi_n\rightarrow \chi$ in distribution as $S_R$-valued random elements. Then $\int_0^{\cdot} \chi_n(s)\,ds$ converges in distribution as $\mathcal{C}([0,T];U)$-valued processes to $\int_0^{\cdot} \chi(s)\,ds$ as $n\rightarrow+\infty.$
\end{lemma}

\begin{theorem}\label{4.2.5}
Assume that $\chi\in L^2(0,T; U_0),$ $g\in L^2(0,T;L^2(\mathcal{O}))$ and $\sigma$ satisfies assumptions $(A_1)$-$(A_3).$ Then for any $\theta_0\in H_2,$ there exists a unique weak solution $(v,p_s,\theta)\in \mathcal{C}([0,T];V_1\times H^1(M)\times H_2)\times L^2(0,T;H^2(\mathcal{O})\times H^2(M)\times V_2)$ of problem
\begin{equation}\label{4.2.6}
\begin{cases}
&\nabla p_s+fv^\bot+L_1v=\int_{-h}^z\nabla \theta(x,y,\zeta,t)\,d\zeta,\\
&\int_{-h}^{0}\nabla\cdot v(x,y,\zeta,t)\,d\zeta=0,\\
&\frac{\partial\theta}{\partial t}+B(v,\theta)+A_2(\theta-\theta^*)-K_h\Delta \theta^*=g(x,y,z,t)+\sigma(t,\theta(t))\chi(t),\\
&\left.A_{\nu}\frac{\partial v}{\partial z}\right|_{\Gamma_u}=\mu, \left.\frac{\partial v}{\partial z}\right|_{\Gamma_b}=0, \left.v\cdot\vec{n}\right|_{\Gamma_l}=0,\left.\frac{\partial v}{\partial\vec{n}}\times\vec{n}\right|_{\Gamma_l}=0,\\
&\theta(0)=\theta_0(x,y,z).
\end{cases}
\end{equation}
\end{theorem}
{\bf Proof.} For any $\eta\in\mathcal{C}([0,T];H_2)$ and any fixed $\chi\in L^2(0,T;U_0),$ we conclude from assumption ($A_2$) that $\sigma(t,\eta(t))\chi(t)\in L^2(0,T;H_2)$ and the following result hold:
\begin{align}\label{4.2.7}
\|\sigma(\cdot,\eta)\chi\|_{L^2(0,T;H_2)}^2\leq K(1+\|\eta\|_{\mathcal{C}([0,T]; H_2)}^2)\|\chi\|_{L^2(0,T;U_0)}^2.
\end{align}
Therefore, it follows from the proof of Theorem 3.1 in \cite{ccs} that for any $\theta_0\in H_2$ and any $T>0,$ there exists a unique weak solution $(v^\eta,p_s^\eta,\theta^\eta)\in \mathcal{C}([0,T];V_1\times H^1(\mathcal{M})\times H_2)\times L^2(0,T;H^2(\mathcal{O})\times H^2(\mathcal{M})\times V_2)$ with $\theta_t^\eta\in L^2(0,T;V_2')$ of the following problem
\begin{equation}\label{4.2.8}
\begin{cases}
&\nabla p_s+fv^\bot+L_1v=\int_{-h}^z\nabla \theta(x,y,\zeta,t)\,d\zeta,\\
&\int_{-h}^{0}\nabla\cdot v(x,y,\zeta,t)\,d\zeta=0,\\
&\frac{\partial\theta}{\partial t}+B(v,\theta)+A_2(\theta-\theta^*)-K_h\Delta \theta^*=g(x,y,z,t)+\sigma(t,\eta(t))\chi(t),\\
&\left.A_{\nu}\frac{\partial v}{\partial z}\right|_{\Gamma_u}=\mu, \left.\frac{\partial v}{\partial z}\right|_{\Gamma_b}=0, \left.v\cdot\vec{n}\right|_{\Gamma_l}=0,\left.\frac{\partial v}{\partial\vec{n}}\times\vec{n}\right|_{\Gamma_l}=0,\\
&\theta(0)=\theta_0(x,y,z).
\end{cases}
\end{equation}
Define the operator $S: \mathcal{C}([0,T];H_2)\rightarrow Y=\{\phi\in L^2(0,T;V_2):\phi_t\in L^2(0,T;V_2')\}\subset \mathcal{C}([0,T];H_2)$ by
\begin{align}
\theta^\eta=S(\eta).
\end{align}
In what follows, we will prove that the well-posedness of weak solutions for problem \eqref{4.2.6} by using Banach's fixed point Theorem. To do this, we will show that the operator $S: \mathcal{C}([0,T];H_2)\rightarrow  \mathcal{C}([0,T];H_2)$ on some subset of $\mathcal{C}([0,T];H_2)$ is contractive.

Multiplying the third equation of problem \eqref{4.2.8} by $\theta^\eta$ and integrating  the resulting equality over $\mathcal{O},$ we obtain
\begin{align*}
\nonumber&\frac{1}{2}\frac{d}{dt}\|\theta^\eta(t)\|_{L^2(\mathcal{O})}^2+\|\theta^\eta(t)\|^2\\
\nonumber=&\int_{\mathcal{O}}g(t)\theta^\eta(t)\,dxdydz+\int_{\mathcal{O}}\sigma(t,\eta(t))\chi(t)\theta^\eta(t)\,dxdydz+\beta\int_{\Gamma_u}\theta^*\theta^\eta\,dxdy\\
\nonumber\leq&\|g(t)\|_{L^2(\mathcal{O})}\|\theta^\eta(t)\|_{L^2(\mathcal{O})}+\beta\|\theta^*\|_{L^2(\mathcal{M})}\|\theta^\eta(t)\|_{L^2(\Gamma_u)}+\|\sigma(t,\eta(t))\chi(t)\|_{L^2(\mathcal{O})}\|\theta^\eta(t)\|_{L^2(\mathcal{O})}\\
\leq&C\|g(t)\|_{L^2(\mathcal{O})}^2+C\|\theta^*\|_{L^2(\mathcal{M})}^2+C\|\sigma(t,\eta(t))\chi(t)\|_{L^2(\mathcal{O})}^2+\frac{1}{2}\|\theta^\eta(t)\|^2.
\end{align*}
It follows from Lemma \ref{2.2.4} and inequality \eqref{4.2.7} that
\begin{align}\label{4.2.9}
\nonumber&\sup_{s\in[0,T]}\|\theta^\eta(s)\|_{L^2(\mathcal{O})}^2+\int_0^T\|\theta^\eta(s)\|^2\,ds\\
\nonumber\leq&C\int_0^T\|g(s)\|_{L^2(\mathcal{O})}^2\,ds+CT\|\theta^*\|_{L^2(\mathcal{M})}^2+C\|\sigma(\cdot,\eta)\chi\|_{L^2(0,T;L^2(\mathcal{O}))}^2+\|\theta_0\|_{L^2(\mathcal{O})}^2\\
\leq&C_{13}(\int_0^T\|g(s)\|_{L^2(\mathcal{O})}^2\,ds+T\|\theta^*\|_{L^2(\mathcal{M})}^2)+C_{12}K(1+\|\eta\|_{\mathcal{C}([0,T]; H_2)}^2)\|\chi\|_{L^2(0,T;U_0)}^2+\|\theta_0\|_{L^2(\mathcal{O})}^2
\end{align}
for some fixed positive constants $C_{12}$ and $C_{13}.$

Choosing $T_1>0$ sufficiently small such that 
\begin{align*}
C_{12}K\|\chi\|_{L^2(0,T_1;U_0)}^2\leq\frac{1}{2}.
\end{align*}
Let $R_0(T)=\left(\frac{1}{1-C_{12}K\|\chi\|_{L^2(0,T;U_0)}^2}\left(C_{13}\int_0^T\|g(s)\|_{L^2(\mathcal{O})}^2\,ds+C_{13}T\|\theta^*\|_{L^2(\mathcal{M})}^2+C_{12}K\|\chi\|_{L^2(0,T;U_0)}^2+\|\theta_0\|_{L^2(\mathcal{O})}^2\right)\right)^{\frac{1}{2}}$ for any $T\in(0,T_1]$ and define 
\begin{align*}
B_0(T)=\{\eta\in\mathcal{C}([0,T]; H_2) :\|\eta\|_{\mathcal{C}([0,T]; H_2)}\leq R_0\},
\end{align*}
then we infer from inequality \eqref{4.2.9} that
\begin{align*}
SB_0(T)\subset B_0(T)
\end{align*}
for any $T\in(0,T_1].$

Next, we will prove that the operator $S:B_0(T)\rightarrow B_0(T)$ is contractive. Assume that $\eta_1,$ $\eta_2\in B_0(T)$ and $(v_1,p_{1s},\theta_1),$ $(v_2,p_{2s},\theta_2),$ are the weak solutions for problem \eqref{4.2.8} corresponding to $\eta_1,$ $\eta_2,$ respectively. Put $(v,p_s,\theta)=(v_1,p_{1s},\theta_1)-(v_2,p_{2s},\theta_2)$ and $\eta=\eta_1-\eta_2,$ then $(v,p_s,\theta,\eta)$ satisfies the following problem
\begin{equation}\label{4.2.10}
\begin{cases}
&\nabla p_s+fv^\bot+L_1v=\int_{-h}^z\nabla \theta(x,y,\zeta,t)\,d\zeta,\\
&\int_{-h}^{0}\nabla\cdot v(x,y,\zeta,t)\,d\zeta=0,\\
&\frac{\partial\theta}{\partial t}+B(v_1,\theta)+B(v,\theta_2)+A_2\theta=\sigma(t,\eta_1(t))\chi(t)-\sigma(t,\eta_2(t))\chi(t),\\
&\left.A_{\nu}\frac{\partial v}{\partial z}\right|_{\Gamma_u}=0, \left.\frac{\partial v}{\partial z}\right|_{\Gamma_b}=0, \left.v\cdot\vec{n}\right|_{\Gamma_l}=0,\left.\frac{\partial v}{\partial\vec{n}}\times\vec{n}\right|_{\Gamma_l}=0,\\
&\theta(0)=0
\end{cases}
\end{equation}
and the following estimate hold:
\begin{align}\label{4.2.11}
\|v(t)\|_{H^{\gamma+1}(\mathcal{O})}^2+\|p_{s}(t)\|_{H^{\gamma}(\mathcal{M})}^2\leq \mathcal{K}_1\|\theta(t)\|_{H^\gamma(\mathcal{O})}^2
\end{align}
for $\gamma=0,1.$

Taking the inner product of the third equation of problem \eqref{4.2.10} with $\theta$ in $H_2$ and combining H\"{o}lder inequality with inequality \eqref{4.2.11}, Lemma \ref{2.2.3}, we obtain
\begin{align*}
&\frac{1}{2}\frac{d}{dt}\|\theta(t)\|_{L^2(\mathcal{O})}^2+\|\theta(t)\|^2\\
=&-b(v(t)\,\theta_2(t),\theta(t))+\int_{\mathcal{O}}(\sigma(t,\eta_1(t))\chi(t)-\sigma(t,\eta_2(t))\chi(t))\theta(t)\,dxdydz\\
\leq&\mathcal{K}\|v(t)\|_{H^1(\mathcal{O})}^{\frac{1}{2}}\|v(t)\|_{H^2(\mathcal{O})}^{\frac{1}{2}}\|\theta_2(t)\|\|\theta(t)\|_{L^2(\mathcal{O})}^{\frac{1}{2}}\|\theta(t)\|^{\frac{1}{2}}+\|\sigma(t,\eta_1(t))\chi(t)-\sigma(t,\eta_2(t))\chi(t)\|_{L^2(\mathcal{O})}\|\theta(t)\|_{L^2(\mathcal{O})}\\
\leq&\mathcal{K}\sqrt{\mathcal{K}_1}\|\theta_2(t)\|\|\theta(t)\|_{L^2(\mathcal{O})}\|\theta(t)\|+\|\sigma(t,\eta_1(t))-\sigma(t,\eta_2(t))\|_{L^2(U_0,H_2)}\|\chi(t)\|_{U_0}\|\theta(t)\|_{L^2(\mathcal{O})}\\
\leq&\frac{1}{2}\|\theta(t)\|^2+\mathcal{K}^2\mathcal{K}_1\|\theta_2(t)\|^2\|\theta(t)\|_{L^2(\mathcal{O})}^2+\frac{L}{K_2}\|\eta(t)\|_{L^2(\mathcal{O})}^2\|\chi(t)\|_{U_0}^2.
\end{align*}

We infer from Lemma \ref{2.2.4} that for any $T\in (0,T_1],$
\begin{align}\label{4.2.12}
\nonumber\sup_{s\in[0,T]}\|\theta(s)\|_{L^2(\mathcal{O})}^2+\int_0^T\|\theta(s)\|^2\,ds
\leq &\frac{2L}{K_2}\|\eta\|_{\mathcal{C}([0,T];L^2(\mathcal{O}))}^2\|\chi\|_{L^2(0,T;U_0)}^2e^{2\mathcal{K}^2\mathcal{K}_1\int_0^T\|\theta_2(s)\|^2\,ds}\\
\leq &\frac{2L}{K_2}\|\eta\|_{\mathcal{C}([0,T];L^2(\mathcal{O}))}^2\|\chi\|_{L^2(0,T;U_0)}^2e^{2\mathcal{K}^2\mathcal{K}_1R_0(T)^2}.
\end{align}
Obviously, we can choose $T_0\in (0,T_1]$ such that
\begin{align*}
\max\left\{C_{12}K\|\chi\|_{L^2(0,T_0;U_0)}^2,\frac{2L}{K_2}\|\chi\|_{L^2(0,T_0;U_0)}^2e^{2\mathcal{K}^2\mathcal{K}_1R_0(T_1)^2}\right\}\leq\frac{1}{2},
\end{align*}
then the operator $S:B_0(T_0)\rightarrow B_0(T_0)$ is a contractive mapping. We infer from Banach fixed point Theorem that for any $\theta_0\in H_2$ and $\chi\in L^2(0,T;U_0),$ the operator $S$ has a unique fixed point $\theta=S(\theta)\in B_0(T_0),$ which implies that problem \eqref{4.2.6} possesses a unique weak solution $(v, p_s,\theta)\in \mathcal{C}([0,T_0];V_1\times H^1(M)\times H_2)\times L^2(0,T_0;H^2(\mathcal{O})\times H^2(M)\times V_2)$ on $[0,T_0].$

Multiplying the third equation of problem \eqref{4.2.6} by $\theta$ and integrating  the resulting equality over $\mathcal{O},$ we obtain
\begin{align*}
\nonumber&\frac{1}{2}\frac{d}{dt}\|\theta(t)\|_{L^2(\mathcal{O})}^2+\|\theta(t)\|^2\\
\nonumber=&\int_{\mathcal{O}}g(t)\theta(t)\,dxdydz+\int_{\mathcal{O}}\sigma(t,\theta(t))h(t)\theta^\eta(t)\,dxdydz+\beta\int_{\Gamma_u}\theta^*\theta\,dxdy\\
\nonumber\leq&\|g(t)\|_{L^2(\mathcal{O})}\|\theta(t)\|_{L^2(\mathcal{O})}+\beta\|\theta^*\|_{L^2(\mathcal{M})}\|\theta(t)\|_{L^2(\Gamma_u)}+\|\sigma(t,\theta(t))\chi(t)\|_{L^2(\mathcal{O})}\|\theta(t)\|_{L^2(\mathcal{O})}\\
\nonumber\leq&\|g(t)\|_{L^2(\mathcal{O})}\|\theta(t)\|_{L^2(\mathcal{O})}+\beta\|\theta^*\|_{L^2(\mathcal{M})}\|\theta(t)\|_{L^2(\Gamma_u)}+\|\sigma(t,\theta(t))\|_{L^2(U_0,H_2)}\|\chi(t)\|_{U_0}\|\theta(t)\|_{L^2(\mathcal{O})}\\
\leq&C\|g(t)\|_{L^2(\mathcal{O})}^2+C\|\theta^*\|_{L^2(\mathcal{M})}^2+CK(1+\|\theta(t)\|_{L^2(\mathcal{O})}^2)\|\chi(t)\|_{U_0}^2+\frac{1}{2}\|\theta(t)\|^2,
\end{align*}
which implies that
\begin{align*}
\nonumber&\frac{d}{dt}\|\theta(t)\|_{L^2(\mathcal{O})}^2+\|\theta(t)\|^2\\
\leq&C\|g(t)\|_{L^2(\mathcal{O})}^2+C\|\theta^*\|_{L^2(\mathcal{M})}^2+CK(1+\|\theta(t)\|_{L^2(\mathcal{O})}^2)\|\chi(t)\|_{U_0}^2.
\end{align*}
We conclude from Lemma \ref{2.2.4} that
\begin{align}\label{4.2.13}
\nonumber&\|\theta(t)\|_{L^2(\mathcal{O})}^2+\int_0^t\|\theta(s)\|^2\,ds\\
\leq&C\left(\int_0^t\|g(s)\|_{L^2(\mathcal{O})}^2\,ds+t\|\theta^*\|_{L^2(\mathcal{M})}^2+K\int_0^t\|\chi(s)\|_{U_0}^2\,ds+\|\theta_0\|_{L^2(\mathcal{O})}^2\right)e^{CK\int_0^t\|\chi(s)\|_{U_0}^2\,ds},
\end{align}
which implies that the weak solution of problem \eqref{4.2.6} exists globally.\\
\qed\hfill

For any $\chi\in L^2(0,T; U_0),$ define
\begin{align*}
\theta_\chi=\Phi^0\left(\int_0^{\cdot}\chi(s)\,ds\right),
\end{align*}
where $(v_\chi,p_{\chi s},\theta_\chi)$ is the unique weak solution of problem \eqref{4.2.6}.

\begin{theorem}\label{a}
Assume that $\theta_0\in H_2$ and assumptions ($A_1$)-$(A_4)$ hold. Let $M$ be any fixed finite positive constant and define
\begin{align*}
K_M=\left\{\Phi^0\left(\int_0^{\cdot}\chi(s)\,ds\right): \chi\in S_M\right\}.
\end{align*}
Then the set $K_M$ is compact in $\mathcal{C}([0,T];H_2)\cap L^2(0,T;V_2).$
\end{theorem}
{\bf Proof.} Let $\{\theta_n\}$ be a sequence in $K_M,$ where $(v_n, p_{ns},\theta_n)$ is the weak solution of problem \eqref{4.2.6} with $\chi=\chi_n\in S_M.$ By weak compactness of $S_M,$ there exists a subsequence (still use the same notation) of $\{\chi_n\},$ which weakly converges to a limit $\chi$ in $L^2(0,T;U_0).$ In fact $\chi\in S_M$ as $S_M$ is closed. We now would like to show that the corresponding subsequence (still use the same notation) of $\{\theta_n\}$ converges in $\mathcal{C}([0,T];H_2)\cap L^2(0,T;V_2)$ to $\theta,$ where $(v, p_s,\theta)$ is the weak solution of the following "limit" problem
\begin{equation}\label{4.2.15}
\begin{cases}
&\nabla p_s+fv^\bot+L_1v=\int_{-h}^z\nabla \theta(x,y,\zeta,t)\,d\zeta,\\
&\int_{-h}^{0}\nabla\cdot v(x,y,\zeta,t)\,d\zeta=0,\\
&\frac{\partial\theta}{\partial t}+B(v,\theta)+A_2(\theta-\theta^*)-K_h\Delta \theta^*=g(x,y,z,t)+\sigma(t,\theta(t))\chi(t),\\
&\left.A_{\nu}\frac{\partial v}{\partial z}\right|_{\Gamma_u}=\mu, \left.\frac{\partial v}{\partial z}\right|_{\Gamma_b}=0, \left.v\cdot\vec{n}\right|_{\Gamma_l}=0,\left.\frac{\partial v}{\partial\vec{n}}\times\vec{n}\right|_{\Gamma_l}=0,\\
&\theta(0)=\theta_0(x,y,z).
\end{cases}
\end{equation}
From inequality \eqref{4.2.13}, we know that
\begin{align}\label{4.2.16}
\nonumber&\sup_{s\in[0,T]}\|\theta(s)\|_{L^2(\mathcal{O})}^2+\int_0^T\|\theta(s)\|^2\,ds\\
\nonumber\leq&C\left(\int_0^T\|g(s)\|_{L^2(\mathcal{O})}^2\,ds+T\|\theta^*\|_{L^2(\mathcal{M})}^2+\|\theta_0\|_{L^2(\mathcal{O})}^2+K\int_0^T\|\chi(s)\|_{U_0}^2\,ds\right)e^{CK\int_0^T\|\chi(s)\|_{U_0}^2\,ds}\\
\leq&C\left(\int_0^T\|g(s)\|_{L^2(\mathcal{O})}^2\,ds+T\|\theta^*\|_{L^2(\mathcal{M})}^2+\|\theta_0\|_{L^2(\mathcal{O})}^2+KM\right)e^{CKM}.
\end{align}
Put $(u_n,q_{ns},\eta_n)=(v_n-v,p_{ns}-p_s,\theta_n-\theta),$ then it satisfies the following problem
\begin{equation}\label{4.2.17}
\begin{cases}
&\nabla q_{ns}+f\vec{k}\times u_n+L_1u_n=\int_{-h}^z\nabla \eta_n(x,y,\zeta,t)\,d\zeta,\\
&\int_{-h}^{0}\nabla\cdot u_n(x,y,\zeta,t)\,d\zeta=0,\\
&\frac{\partial\eta_n}{\partial t}+B(v_n,\eta_n)+B(u_n,\theta)+A_2\eta_n=\sigma(t,\theta_n(t))\chi_n(t)-\sigma(t,\theta(t))\chi(t),\\
&\left.A_{\nu}\frac{\partial u_n}{\partial z}\right|_{\Gamma_u}=0, \left.\frac{\partial u_n}{\partial z}\right|_{\Gamma_b}=0, \left.u_n\cdot\vec{n}\right|_{\Gamma_l}=0,\left.\frac{\partial u_n}{\partial\vec{n}}\times\vec{n}\right|_{\Gamma_l}=0,\\
&\theta(0)=0,
\end{cases}
\end{equation}
we need to prove that $\eta_n\rightarrow 0$ in $\mathcal{C}([0,T];H_2)\cap L^2(0,T;V_2)$ as $n\rightarrow+\infty.$

Similarly, we have the following estimate:
\begin{align}\label{4.2.18}
\|u_n(t)\|_{H^{\gamma+1}(\mathcal{O})}^2+\|q_{ns}(t)\|_{H^{\gamma}(\mathcal{M})}^2\leq \mathcal{K}_1\|\eta_n(t)\|_{H^\gamma(\mathcal{O})}^2
\end{align}
for $\gamma=0,1.$

Taking the inner product of the third equation of problem \eqref{4.2.17} with $\eta_n$ in $H_2$ and combining H\"{o}lder inequality with inequality \eqref{4.2.18}, we obtain
\begin{align*}
&\frac{1}{2}\frac{d}{dt}\|\eta_n(t)\|_{L^2(\mathcal{O})}^2+\|\eta_n(t)\|^2\\
=&-b(u_n(t),\theta(t),\eta_n(t))+\int_{\mathcal{O}}(\sigma(t,\theta_n(t))\chi_n(t)-\sigma(t,\theta(t))\chi(t))\eta_n(t)\,dxdydz\\
\leq&\mathcal{K}\|u_n(t)\|_{H^1(\mathcal{O})}^{\frac{1}{2}}\|u_n(t)\|_{H^2(\mathcal{O})}^{\frac{1}{2}}\|\theta(t)\|\|\eta_n(t)\|_{L^2(\mathcal{O})}^{\frac{1}{2}}\|\eta_n(t)\|^{\frac{1}{2}}+\int_{\mathcal{O}}\sigma(t,\theta(t))(\chi_n(t)-\chi(t))\eta_n(t)\,dxdydz\\
&+\|\sigma(t,\theta_n(t))\chi_n(t)-\sigma(t,\theta(t))\chi_n(t)\|_{L^2(\mathcal{O})}\|\eta_n(t)\|_{L^2(\mathcal{O})}\\
\leq&\mathcal{K}\sqrt{\mathcal{K}_1}\|\theta(t)\|\|\eta_n(t)\|_{L^2(\mathcal{O})}\|\eta_n(t)\|+\|\sigma(t,\theta_n(t))-\sigma(t,\theta(t))\|_{L^2(U_0,H_2)}\|\chi_n(t)\|_{U_0}\|\eta_n(t)\|_{L^2(\mathcal{O})}\\
&+\int_{\mathcal{O}}\sigma(t,\theta(t))(\chi_n(t)-\chi(t))\eta_n(t)\,dxdydz\\
\leq&\frac{1}{2}\|\eta_n(t)\|^2+(\mathcal{K}^2\mathcal{K}_1\|\theta(t)\|^2+\frac{L}{K_2}\|\chi_n(t)\|_{U_0}^2)\|\eta_n(t)\|_{L^2(\mathcal{O})}^2+\int_{\mathcal{O}}\sigma(t,\theta(t))(\chi_n(t)-\chi(t))\eta_n(t)\,dxdydz.
\end{align*}

We infer from Lemma \ref{2.2.4} and inequality \eqref{4.2.16} that 
\begin{align}\label{4.2.19}
\nonumber&\sup_{s\in[0,T]}\|\eta_n(s)\|_{L^2(\mathcal{O})}^2+\int_0^T\|\eta_n(s)\|^2\,ds\\
\nonumber\leq &C\left(\sup_{t\in[0,T]}\left|\int_0^t\int_{\mathcal{O}}\sigma(s,\theta(s))(\chi_n(s)-\chi(s))\eta_n(s)\,dxdydzds\right|\right)e^{\int_0^T(\mathcal{K}^2\mathcal{K}_1\|\theta(t)\|^2+\frac{L}{K_2}\|\chi_n(t)\|_{U_0}^2)\,ds}\\
\leq &C\left(\sup_{t\in[0,T]}\left|\int_0^t\int_{\mathcal{O}}\sigma(s,\theta(s))(\chi_n(s)-\chi(s))\eta_n(s)\,dxdydzds\right|\right)e^{\mathcal{K}^2\mathcal{K}_1\left(\int_0^T\|g(s)\|_{L^2(\mathcal{O})}^2\,ds+T\|\theta^*\|_{L^2(\mathcal{M})}^2+KM\right)e^{CKM}+\frac{LM}{K_2}}.
\end{align}
We conclude from inequality \eqref{4.2.13} that there exists a positive constant $\bar{C}$ such that 
\begin{align*}
\sup_{n\geq 1}\left(\sup_{t\in[0,T]}(\|\theta(t)\|_{L^2(\mathcal{O})}^2+\|\theta_n(t)\|_{L^2(\mathcal{O})}^2)+\int_0^T(\|\theta(s)\|^2+\|\theta_n(s)\|^2)\,ds\right)\leq \bar{C}
\end{align*}
In what follows, we will give an estimate on the term $\sup_{t\in[0,T]}\left|\int_0^t\int_{\mathcal{O}}\sigma(s,\theta(s))(\chi_n(s)-\chi(s))\eta_n(s)\,dxdydzds\right|.$ For $N\geq 1$ and $k=0,1,\cdots, 2^N,$ let $t_k=kT2^{-N}$ and let
\begin{align}\label{4.2.20}
\sup_{t\in[0,T]}\left|\int_0^t\int_{\mathcal{O}}\sigma(s,\theta(s))(\chi_n(s)-\chi(s))\eta_n(s)\,dxdydzds\right|\leq\sum_{i=1}^5 I_{n,N}^i,
\end{align}
where
\begin{align*}
I_{n,N}^1=&\sum_{k=1}^{2^N}\int_{t_{k-1}}^{t_k}\left|\int_{\mathcal{O}}\sigma(s,\theta(s))(\chi_n(s)-\chi(s))(\eta_n(s)-\eta_n(t_k))\,dxdydz\right|\,ds,\\
I_{n,N}^2=&\sum_{k=1}^{2^N}\int_{t_{k-1}}^{t_k}\left|\int_{\mathcal{O}}(\sigma(s,\theta(s)-\sigma(t_k,\theta(s)))(\chi_n(s)-\chi(s))\eta_n(t_k)\,dxdydz\right|\,ds,\\
I_{n,N}^3=&\sum_{k=1}^{2^N}\int_{t_{k-1}}^{t_k}\left|\int_{\mathcal{O}}(\sigma(t_k,\theta(s))-\sigma(t_k,\theta(t_k)))(\chi_n(s)-\chi(s))\eta_n(t_k)\,dxdydz\right|\,ds,\\
I_{n,N}^4=&\sup_{1\leq k\leq 2^N}\sup_{t_{k-1}\leq t\leq t_k}\left|\int_{\mathcal{O}}\sigma(t_k,\theta(t_k))\int_{t_{k-1}}^t(\chi_n(s)-\chi(s))\,ds\,\eta_n(t_k)\,dxdydz\right|,\\
I_{n,N}^5=&\sum_{k=1}^{2^N}\left|\int_{t_{k-1}}^{t_k}\int_{\mathcal{O}}\sigma(t_k,\theta(t_k))(\chi_n(s)-\chi(s))\eta_n(t_k)\,dxdydzds\right|.
\end{align*}
It follows from Cauchy-Schwarz's inequality, assumptions ($A_2$), ($A_4$), H\"{o}lder's inequality and the proof of Theorem \ref{5.1} that 
\begin{align}\label{4.2.21}
\nonumber I_{n,N}^1=&\sum_{k=1}^{2^N}\int_{t_{k-1}}^{t_k}\left|\int_{\mathcal{O}}\sigma(s,\theta(s))(\chi_n(s)-\chi(s))(\eta_n(s)-\eta_n(t_k))\,dxdydz\right|\,ds\\
\nonumber\leq&\sum_{k=1}^{2^N}\int_{t_{k-1}}^{t_k}\|\sigma(s,\theta(s))\|_{\mathcal{L}_2(U_0,H_2)}\|\chi_n(s)-\chi(s)\|_{U_0}\|\eta_n(s)-\eta_n(t_k))\|_{L^2(\mathcal{O})}\,ds\\
\nonumber\leq&\sqrt{K(1+\bar{C})}\left(\sum_{k=1}^{2^N}\int_{t_{k-1}}^{t_k}\|\chi_n(s)-\chi(s)\|_{U_0}\|^2\,ds\right)^{\frac{1}{2}}\left(\sum_{k=1}^{2^N}\int_{t_{k-1}}^{t_k}\|\eta_n(s)-\eta_n(t_k))\|_{L^2(\mathcal{O})}^2\,ds\right)^{\frac{1}{2}}\\
\nonumber\leq&2\sqrt{K(1+\bar{C})M}\left(\sum_{k=1}^{2^N}\int_{t_{k-1}}^{t_k}\|\theta_n(s)-\theta_n(t_k))\|_{L^2(\mathcal{O})}^2\,ds\right)^{\frac{1}{2}}+2\sqrt{K(1+\bar{C})M}\left(\sum_{k=1}^{2^N}\int_{t_{k-1}}^{t_k}\|\theta(s)-\theta(t_k))\|_{L^2(\mathcal{O})}^2\,ds\right)^{\frac{1}{2}}\\
\leq&\frac{4\sqrt{\mathcal{D}K(1+\bar{C})M}}{2^{\frac{N}{4}}},
\end{align}

\begin{align}\label{4.2.22}
\nonumber I_{n,N}^2=&\sum_{k=1}^{2^N}\int_{t_{k-1}}^{t_k}\left|\int_{\mathcal{O}}(\sigma(s,\theta(s)-\sigma(t_k,\theta(s)))(\chi_n(s)-\chi(s))\eta_n(t_k)\,dxdydz\right|ds\\
\nonumber\leq&\sum_{k=1}^{2^N}\int_{t_{k-1}}^{t_k}\|\sigma(s,\theta(s)-\sigma(t_k,\theta(s))\|_{\mathcal{L}_2(U_0,H_2)}\|\chi_n(s)-\chi(s)\|_{U_0}\|\eta_n(t_k)\|_{L^2(\mathcal{O})}\,ds\\
\nonumber\leq&L_1\sum_{k=1}^{2^N}\int_{t_{k-1}}^{t_k}(1+\|\theta(s)\|_{L^2(\mathcal{O})})(s-t_k)^\gamma\|\chi_n(s)-\chi(s)\|_{U_0}\|\eta_n(t_k)\|_{L^2(\mathcal{O})}\,ds\\
\nonumber\leq&\frac{2L_1\sqrt{\bar{C}}}{2^{\gamma N}}\left(\sum_{k=1}^{2^N}\int_{t_{k-1}}^{t_k}(1+\|\theta(s)\|_{L^2(\mathcal{O})})^2\,ds\right)^{\frac{1}{2}}\left(\sum_{k=1}^{2^N}\int_{t_{k-1}}^{t_k}\|\chi_n(s)-\chi(s)\|_{U_0}^2\,ds\right)^{\frac{1}{2}}\\
\leq&\frac{4L_1\sqrt{2TM(1+\bar{C})\bar{C}}}{2^{\gamma N}},
\end{align}

\begin{align}\label{4.2.23}
\nonumber I_{n,N}^3=&\sum_{k=1}^{2^N}\int_{t_{k-1}}^{t_k}\left|\int_{\mathcal{O}}(\sigma(t_k,\theta(s))-\sigma(t_k.\theta(t_k))(\chi_n(s)-\chi(s))\eta_n(t_k)\,dxdydz\right|ds\\
\nonumber\leq&\sum_{k=1}^{2^N}\int_{t_{k-1}}^{t_k}\|\sigma(t_k,\theta(s))-\sigma(t_k,\theta(t_k))\|_{\mathcal{L}_2(U_0,H_2)}\|\chi_n(s)-\chi(s)\|_{U_0}\|\eta_n(t_k)\|_{L^2(\mathcal{O})}\,ds\\
\nonumber\leq&2\sqrt{L\bar{C}}\sum_{k=1}^{2^N}\int_{t_{k-1}}^{t_k}\|\theta(s)-\theta(t_k)\|\|_{L^2(\mathcal{O})}\|\chi_n(s)-\chi(s)\|_{U_0}\,ds
\end{align}
\begin{align}
\nonumber\leq&2\sqrt{L\bar{C}}\left(\sum_{k=1}^{2^N}\int_{t_{k-1}}^{t_k}\|\theta(s)-\theta(t_k)\|_{L^2(\mathcal{O})}^2\,ds\right)^{\frac{1}{2}}\left(\sum_{k=1}^{2^N}\int_{t_{k-1}}^{t_k}\|\chi_n(s)-\chi(s)\|_{U_0}^2\,ds\right)^{\frac{1}{2}}\\
\leq&\frac{4\sqrt{\mathcal{D}LM\bar{C}}}{2^{\frac{N}{4}}},
\end{align}

\begin{align}\label{4.2.24}
\nonumber I_{n,N}^4=&\sup_{1\leq k\leq 2^N}\sup_{t_{k-1}\leq t\leq t_k}\left|\int_{\mathcal{O}}\sigma(t_k,\theta(t_k))\int_{t_{k-1}}^t(\chi_n(s)-\chi(s))\,ds\,\eta_n(t_k)\,dxdydz\right|\\
\nonumber \leq&\sup_{1\leq k\leq 2^N}\sup_{t_{k-1}\leq t\leq t_k}\left(\|\sigma(t_k,\theta(t_k))\|_{\mathcal{L}_2(U_0,H_2)}\left\|\int_{t_{k-1}}^t(\chi_n(s)-\chi(s))\,ds\right\|_{U_0}\|\eta_n(t_k)\|_{L^2(\mathcal{O})}\right)\\
\nonumber\leq&2\sqrt{K\bar{C}(1+\bar{C})}\sup_{1\leq k\leq 2^N}\left(\int_{t_{k-1}}^{t_k}\|\chi_n(s)-\chi(s)\|_{U_0}\,ds\right)\\
\nonumber\leq&\frac{2\sqrt{K\bar{C}(1+\bar{C})}}{2^{\frac{N}{2}}}\left(\int_0^T\|\chi_n(s)-\chi(s)\|_{U_0}^2\,ds\right)^\frac{1}{2}\\
\leq&\frac{4\sqrt{KM\bar{C}(1+\bar{C})}}{2^{\frac{N}{2}}}.
\end{align}
For any fixed $N$ and any $k=1,2,\cdots,2^N,$ from the weak convergence of $\chi_n$ to $\chi$ in $L^2(0,T;U_0)$ as $n\rightarrow+\infty,$ we deduce that the term $\int_{t_{k-1}}^{t_k}(\chi_n(s)-\chi(s))\,ds$ weakly converge to $0$ in $U_0.$ Since $\sigma(t_k,\theta(t_k))$ is a compact operator, we infer that for fixed $k,$ the sequence $\sigma(t_k,\theta(t_k))\int_{t_{k-1}}^{t_k}(\chi_n(s)-\chi(s))\,ds$ strongly converge to $0$ in $H_2$ as $n\rightarrow+\infty.$ But $\sup_{n\geq 1,1\leq k\leq 2^N}\|\eta_n(t_k)\|_{L^2(\mathcal{O})}\leq 2\sqrt{\bar{C},}$ which implies that 
\begin{align}\label{4.2.25}
\lim_{n\rightarrow+\infty}I_{n,N}^5=0.
\end{align}
Thus, we deduce from inequalities \eqref{4.2.19}-\eqref{4.2.25} that there exists a positive constant $\mathcal{K}_2$ independent of $n,$ $N$ such that for any integer $N\geq 1,$
\begin{align*}
\limsup_{n\rightarrow+\infty}\left(\sup_{s\in[0,T]}\|\eta_n(s)\|_{L^2(\mathcal{O})}^2+\int_0^T\|\eta_n(s)\|^2\,ds\right)\leq \mathcal{K}_22^{-(\gamma\land\frac{1}{4})N},
\end{align*}
which implies that
\begin{align*}
\lim_{n\rightarrow+\infty}\left(\sup_{s\in[0,T]}\|\eta_n(s)\|_{L^2(\mathcal{O})}^2+\int_0^T\|\eta_n(s)\|^2\,ds\right)=0.
\end{align*}
Therefore, for every sequence $\{\theta_n\}$ in $K_M,$ there exists a subsequence $\{\theta_{n_k}\}$ which converges to some element $\theta_\chi\in K_M$ in $\mathcal{C}([0,T];H_2)\cap L^2(0,T;V_2),$ i.e., the set $K_M$ is compact in $\mathcal{C}([0,T];H_2)\cap L^2(0,T;V_2).$\\
\qed\hfill


\begin{theorem}\label{4.2.26}
Under the assumptions ($A_1$)-$(A_4).$ Let $\{\chi^\epsilon\}_{\epsilon>0}\subset\mathcal{A}_M$ for some fixed finite positive constant $M<+\infty.$ Assume that $\chi^\epsilon$ converges to $\chi$ in distribution as $S_M$-valued random elements, then
\begin{align*}
\Phi^\epsilon\left(W_{\cdot}+\frac{1}{\sqrt{\epsilon}}\int_0^{\cdot} \chi^\epsilon(s)\,ds\right)\rightarrow \Phi^0\left(\int_0^{\cdot}\chi(s)\,ds\right)
\end{align*}
in distribution as $\epsilon\rightarrow 0.$
\end{theorem}
{\bf Proof.}  Since $\mathcal{A}_M$ is a Polish space, by the Skorokhod representation theorem, we can construct processes $(\tilde{\chi}_\epsilon,\tilde{\chi},\tilde{W})$ such that the joint distribution of $(\tilde{\chi}_\epsilon,\tilde{W})$ is the same as that of $(\chi_\epsilon,W)$, the distribution of $\tilde{\chi}$ coincides with that of $\chi,$ and $\tilde{\chi}_\epsilon\rightarrow\tilde{\chi}$ a.s., in the (weak) topology of $S_M.$ Hence a.s. for every $t\in[0,T],$ $\int_0^t\tilde{\chi}_\epsilon(s)\,ds-\int_0^t\tilde{\chi}(s)\,ds\rightarrow 0$ weakly in $U_0.$ Let $\theta^\epsilon(\cdot)=\Phi^\epsilon(W(\cdot)+\frac{1}{\sqrt{\epsilon}}\int_0^{\cdot}\chi^\epsilon(s)\,ds).$ By the Girsanov theorem and the uniqueness of solution for problem \eqref{2.2.1}, we know that $(v^\epsilon(\cdot),p_{s}^\epsilon(\cdot),\theta^\epsilon(\cdot))$ is the unique weak solution of problem
\begin{equation}\label{4.2.27}
\begin{cases}
&\nabla p^\epsilon_{s}+f{v^\epsilon}^\bot+L_1v^\epsilon=\int_{-h}^z\nabla \theta^\epsilon(x,y,\zeta,t)\,d\zeta,\\
&\int_{-h}^{0}\nabla\cdot v^\epsilon(x,y,\zeta,t)\,d\zeta=0,\\
&d\theta^\epsilon+B(v^\epsilon,\theta^\epsilon)\,dt+A_2(\theta^\epsilon-\theta^*)\,dt-K_h\Delta \theta^*\,dt\\
&=g(x,y,z,t)\,dt+\sigma(t,\theta^\epsilon(t))\chi^\epsilon(t)\,dt+\sqrt{\epsilon}\sigma(t,\theta^\epsilon(t))dW(t),\\
&\left.A_{\nu}\frac{\partial v^\epsilon}{\partial z}\right|_{\Gamma_u}=\mu, \left.\frac{\partial v^\epsilon}{\partial z}\right|_{\Gamma_b}=0, \left.v^\epsilon\cdot\vec{n}\right|_{\Gamma_l}=0,\left.\frac{\partial v^\epsilon}{\partial\vec{n}}\times\vec{n}\right|_{\Gamma_l}=0,\\
&\theta^\epsilon(x,y,z,0)=\theta_0(x,y,z).
\end{cases}
\end{equation}
Now, we need to prove $\theta^\epsilon\rightarrow \theta=\Phi^0(\int_0^{\cdot}\chi(s)\,ds)$ in distribution as $\epsilon\rightarrow 0.$

Let $(u, q_s,\eta)=(v^\epsilon-v, p_s^\epsilon-p_s,\theta^\epsilon-\theta),$ then the following estimate holds:
\begin{align}\label{4.2.28}
\|u(t)\|_{H^{\gamma+1}(\mathcal{O})}^2+\|q_{s}(t)\|_{H^{\gamma}(\mathcal{M})}^2\leq \mathcal{K}_1\|\eta(t)\|_{H^\gamma(\mathcal{O})}^2
\end{align}
for $\gamma=0,1.$

Applying It\^{o} formula to $\|\eta(t)\|_{L^2(\mathcal{O})}^2$  and combining assumption ($A_3$) with inequality \eqref{4.2.28}, Young's inequality, H\"{o}lder's inequality, Lemma \ref{2.2.3}, yield
\begin{align*}
&d\|\eta(t)\|_{L^2(\mathcal{O})}^2+2\|\eta(t)\|^2\,dt\\
=&-2b(u(t),\theta(t),\eta(t))\,dt+\epsilon\|\sigma(t,\theta^\epsilon(t))\|_{\mathcal{L}_2(U_0,H_2)}^2\,dt+2\sqrt{\epsilon}\int_{\mathcal{O}}\eta(t)\sigma(t,\theta^\epsilon(t))\,dW(t)dxdydz\\
&+2\int_{\mathcal{O}}(\sigma(t,\theta^\epsilon(t))\chi^\epsilon(t)-\sigma(t,\theta(t))\chi(t))\eta(t)\,dxdydzdt\\
\leq&\mathcal{K}\|\theta(t)\|\|\eta(t)\|_{L^2(\mathcal{O})}\|\eta(t)\|+2\|\sigma(t,\theta^\epsilon(t))-\sigma(t,\theta(t))\|_{L^2(U_0,H_2)}\|\chi^\epsilon(t)\|_{U_0}\|\eta(t)\|_{L^2(\mathcal{O})}\\
&+2\sqrt{\epsilon}\int_{\mathcal{O}}\eta(t)\sigma(t,\theta^\epsilon(t))\,dW(t)dxdydz+2\int_{\mathcal{O}}\sigma(t,\theta(t))(\chi^\epsilon(t)-\chi(t))\eta(t)\,dxdydzdt\\
&+2\epsilon\|\sigma(t,\theta^\epsilon(t))-\sigma(t,\theta(t))\|_{\mathcal{L}_2(U_0,H_2)}^2\,dt+2\epsilon\|\sigma(t,\theta(t))\|_{\mathcal{L}_2(U_0,H_2)}^2\,dt\\
\leq&\|\eta(t)\|^2+C(1+\|\theta(t)\|^2+L\|\chi^\epsilon(t)\|_{U_0}^2)\|\eta(t)\|_{L^2(\mathcal{O})}^2+2K\epsilon(1+\|\theta(t))\|_{L^2(\mathcal{O})}^2)\,dt\\
&+2\sqrt{\epsilon}\int_{\mathcal{O}}\eta(t)\sigma(t,\theta^\epsilon(t))\,dW(t)dxdydz+2\int_{\mathcal{O}}\sigma(t,\theta(t))(\chi^\epsilon(t)-\chi(t))\eta(t)\,dxdydzdt,
\end{align*}
which implies that
\begin{align}\label{4.2.29}
\nonumber&\|\eta(t)\|_{L^2(\mathcal{O})}^2+\int_0^t\|\eta(s)\|^2\,ds\\
\nonumber\leq&C\int_0^t(1+\|\theta(s)\|^2+L\|\chi^\epsilon(s)\|_{U_0}^2)\|\eta(s)\|_{L^2(\mathcal{O})}^2\,ds+2K\epsilon\int_0^t(1+\|\theta(s))\|_{L^2(\mathcal{O})}^2)\,ds\\
&+2\sqrt{\epsilon}\int_0^t\int_{\mathcal{O}}\eta(s)\sigma(s,\theta^\epsilon(s))\,dW(s)dxdydz+2\int_0^t\int_{\mathcal{O}}\sigma(s,\theta(s))(\chi^\epsilon(s)-\chi(s))\eta(t)\,dxdydzds.
\end{align}
For any $R>0$ and any $t\in[0,T],$ define
\begin{align*}
G_{R,\epsilon}(t)=\left\{\omega:\left(\sup_{s\in[0,t]}\|\theta(s,\omega)\|_{L^2(\mathcal{O})^2}+\int_0^t\|\theta(s,\omega)\|^2\,ds\right)\lor \left(\sup_{s\in[0,t]}\|\theta^\epsilon(s,\omega)\|_{L^2(\mathcal{O})^2}+\int_0^t\|\theta^\epsilon(s,\omega)\|^2\,ds\right)\leq R\right\}.
\end{align*}
It follows from Lemma \ref{2.2.4} that on $G_{R,\epsilon}(T),$
\begin{align*}
\nonumber&\sup_{s\in [0,T]}\|\eta(t)\|_{L^2(\mathcal{O})}^2+\int_0^T\|\eta(s)\|^2\,ds\\
\nonumber\leq&\left(2K\epsilon T(1+R)+2\sqrt{\epsilon}\sup_{r\in [0,T]}\left|\int_0^r\int_{\mathcal{O}}\eta(s)\sigma(s,\theta^\epsilon(s))\,dW(s)dxdydz\right|\right.\\
&\left.+2\sup_{r\in[0,T]}\left|\int_0^r\int_{\mathcal{O}}\sigma(s,\theta(s))(\chi^\epsilon(s)-\chi(s))\eta(t)\,dxdydzds\right|\right)e^{C(T+R+LM)},
\end{align*}
which entails that
\begin{align}\label{4.2.30}
\nonumber&E\left(I_{G_{R,\epsilon}(T)}\sup_{s\in [0,T]}\|\eta(t)\|_{L^2(\mathcal{O})}^2\right)+E\left(I_{G_{R,\epsilon}(T)}\int_0^T\|\eta(s)\|^2\,ds\right)\\
\nonumber\leq&\left(2K\epsilon T(1+R)+2\sqrt{\epsilon}E\left(I_{G_{R,\epsilon}(T)}\sup_{r\in [0,T]}\left|\int_0^r\int_{\mathcal{O}}\eta(s)\sigma(s,\theta^\epsilon(s))\,dW(s)dxdydz\right|\right)\right.\\
&\left.+2E\left(I_{G_{R,\epsilon}(T)}\sup_{r\in[0,T]}\left|\int_0^r\int_{\mathcal{O}}\sigma(s,\theta(s))(\chi^\epsilon(s)-\chi(s))\eta(t)\,dxdydzds\right|\right)\right)e^{C(T+R+LM)}.
\end{align}
From the Burkholder-Davis-Gundy inequality and $G_{R,\epsilon}(T)\subset G_{R,\epsilon}(t)$ for any $t\in[0,T],$ we conclude
\begin{align}\label{4.2.31}
\nonumber&2\sqrt{\epsilon}E\left(I_{G_{R,\epsilon}(T)}\sup_{0\leq r\leq T}\left|\int_0^r\int_{\mathcal{O}}\eta(s)\sigma(s,\theta^\epsilon(s))\,dW(s)dxdydz\right|\right)\\
\nonumber\leq&2\sqrt{\epsilon}E\left(\sup_{0\leq r\leq T}\left|\int_0^rI_{G_{R,\epsilon}(s)}\int_{\mathcal{O}}\eta(s)\sigma(s,\theta^\epsilon(s))\,dW(s)dxdydz\right|\right)\\
\nonumber\leq &C\sqrt{\epsilon}E\left(\left(\int_0^TI_{G_{R,\epsilon}(s)}\|\eta(s)\|_{L^2(\mathcal{O})}^2\|\sigma(s,\theta^\epsilon(s))\|_{\mathcal{L}_2(U_0,H_2)}^2\,ds\right)^{\frac{1}{2}}\right)
\end{align}
\begin{align}
\nonumber\leq &C\sqrt{\epsilon}E\left(4R\int_0^TI_{G_{R,\epsilon}(s)}\|\sigma(s,\theta^\epsilon(s))\|_{\mathcal{L}_2(U_0,H_2)}^2\,ds\right)^{\frac{1}{2}}\\
\leq&C\sqrt{\epsilon KRT(1+R)}.
\end{align}
It follows from inequalities \eqref{4.2.30}-\eqref{4.2.31} that
\begin{align}\label{4.2.32}
\nonumber&E\left(I_{G_{R,\epsilon}(T)}\sup_{s\in [0,T]}\|\eta(t)\|_{L^2(\mathcal{O})}^2\right)+E\left(I_{G_{R,\epsilon}(T)}\int_0^T\|\eta(s)\|^2\,ds\right)\\
\nonumber\leq&2\left(E\left(I_{G_{R,\epsilon}(T)}\sup_{r\in[0,T]}\left|\int_0^r\int_{\mathcal{O}}\sigma(s,\theta(s))(\chi^\epsilon(s)-\chi(s))\eta(t)\,dxdydzds\right|\right)\right)e^{C(T+R+LM)}\\
&+\left(2K\epsilon T(1+R)+\sqrt{\epsilon KRT(1+R)}\right)e^{C(T+R+LM)}.
\end{align}
In what follows, we will give an estimate on the term $E\left(I_{G_{R,\epsilon}(T)}\sup_{t\in[0,T]}\left|\int_0^t\int_{\mathcal{O}}\sigma(s,\theta(s))(\chi^\epsilon(s)-\chi(s))\eta(s)\,dxdydzds\right|\right)$ as in the proof of Theorem \ref{a}. For $N\geq 1$ and $k=0,1,\cdots, 2^N,$ let $t_k=kT2^{-N}$ and let
\begin{align}\label{4.2.33}
E\left(I_{G_{R,\epsilon}(T)}\sup_{t\in[0,T]}\left|\int_0^t\int_{\mathcal{O}}\sigma(s,\theta(s))(\chi^\epsilon(s)-\chi(s))\eta(s)\,dxdydzds\right|\right)\leq\sum_{i=1}^5 I_{\epsilon,N}^i,
\end{align}
where
\begin{align*}
I_{\epsilon,N}^1=&E\left(I_{G_{R,\epsilon}(T)}\sum_{k=1}^{2^N}\int_{t_{k-1}}^{t_k}\left|\int_{\mathcal{O}}\sigma(s,\theta(s))(\chi^\epsilon(s)-\chi(s))(\eta(s)-\eta(t_k))\,dxdydz\right|\,ds\right),\\
I_{\epsilon,N}^2=&E\left(I_{G_{R,\epsilon}(T)}\sum_{k=1}^{2^N}\int_{t_{k-1}}^{t_k}\left|\int_{\mathcal{O}}(\sigma(s,\theta(s)-\sigma(t_k,\theta(s)))(\chi^\epsilon(s)-\chi(s))\eta(t_k)\,dxdydz\right|\,ds\right),\\
I_{\epsilon,N}^3=&E\left(I_{G_{R,\epsilon}(T)}\sum_{k=1}^{2^N}\int_{t_{k-1}}^{t_k}\left|\int_{\mathcal{O}}(\sigma(t_k,\theta(s))-\sigma(t_k,\theta(t_k)))(\chi^\epsilon(s)-\chi(s))\eta(t_k)\,dxdydz\right|\,ds\right),\\
I_{\epsilon,N}^4=&E\left(I_{G_{R,\epsilon}(T)}\sup_{1\leq k\leq 2^N}\sup_{t_{k-1}\leq t\leq t_k}\left|\int_{\mathcal{O}}\sigma(t_k,\theta(t_k))\int_{t_{k-1}}^t(\chi^\epsilon(s)-\chi(s))\,ds\,\eta(t_k)\,dxdydz\right|\right),\\
I_{\epsilon,N}^5=&E\left(I_{G_{R,\epsilon}(T)}\sum_{k=1}^{2^N}\left|\int_{\mathcal{O}}\sigma(t_k,\theta(t_k))\int_{t_{k-1}}^{t_k}(\chi^\epsilon(s)-\chi(s))\,ds\,\eta(t_k)\,dxdydz\right|\right).
\end{align*}
It follows from Cauchy-Schwarz's inequality, assumptions ($A_2$), ($A_4$), H\"{o}lder's inequality and Theorem \ref{5.1} that 
\begin{align}\label{4.2.34}
\nonumber I_{\epsilon,N}^1=&E\left(I_{G_{R,\epsilon}(T)}\sum_{k=1}^{2^N}\int_{t_{k-1}}^{t_k}\left|\int_{\mathcal{O}}\sigma(s,\theta(s))(\chi^\epsilon(s)-\chi(s))(\eta(s)-\eta(t_k))\,dxdydz\right|\,ds\right)\\
\nonumber\leq&E\left(I_{G_{R,\epsilon}(T)}\sum_{k=1}^{2^N}\int_{t_{k-1}}^{t_k}\|\sigma(s,\theta(s))\|_{\mathcal{L}_2(U_0,H_2)}\|\chi^\epsilon(s)-\chi(s)\|_{U_0}\|\eta(s)-\eta(t_k))\|_{L^2(\mathcal{O})}\,ds\right)\\
\nonumber\leq&\sqrt{K(1+R)}E\left(I_{G_{R,\epsilon}(T)}\sum_{k=1}^{2^N}\int_{t_{k-1}}^{t_k}\|\chi^\epsilon(s)-\chi(s)\|_{U_0}\|^2\,ds\right)^{\frac{1}{2}}E\left(I_{G_{R,\epsilon}(T)}\sum_{k=1}^{2^N}\int_{t_{k-1}}^{t_k}\|\eta(s)-\eta(t_k))\|_{L^2(\mathcal{O})}^2\,ds\right)^{\frac{1}{2}}\\
\nonumber\leq&2\sqrt{K(1+R)M}E\left(I_{G_{R,\epsilon}(T)}\sum_{k=1}^{2^N}\int_{t_{k-1}}^{t_k}\|\theta^\epsilon(s)-\theta^\epsilon(t_k))\|_{L^2(\mathcal{O})}^2\,ds\right)^{\frac{1}{2}}\\
\nonumber&+2\sqrt{K(1+R)M}E\left(I_{G_{R,\epsilon}(T)}\sum_{k=1}^{2^N}\int_{t_{k-1}}^{t_k}\|\theta(s)-\theta(t_k))\|_{L^2(\mathcal{O})}^2\,ds\right)^{\frac{1}{2}}\\
\leq&\frac{4\sqrt{\mathcal{D}K(1+R)M}}{2^{\frac{N}{4}}},
\end{align}

\begin{align}\label{4.2.35}
\nonumber I_{\epsilon,N}^2=&E\left(I_{G_{R,\epsilon}(T)}\sum_{k=1}^{2^N}\int_{t_{k-1}}^{t_k}\left|\int_{\mathcal{O}}(\sigma(s,\theta(s)-\sigma(t_k,\theta(s)))(\chi^\epsilon(s)-\chi(s))\eta(t_k)\,dxdydz\right|ds\right)\\
\nonumber\leq&E\left(I_{G_{R,\epsilon}(T)}\sum_{k=1}^{2^N}\int_{t_{k-1}}^{t_k}\|\sigma(s,\theta(s)-\sigma(t_k,\theta(s))\|_{\mathcal{L}_2(U_0,H_2)}\|\chi^\epsilon(s)-\chi(s)\|_{U_0}\|\eta(t_k)\|_{L^2(\mathcal{O})}\,ds\right)\\
\nonumber\leq&L_1E\left(I_{G_{R,\epsilon}(T)}\sum_{k=1}^{2^N}\int_{t_{k-1}}^{t_k}(1+\|\theta(s)\|_{L^2(\mathcal{O})})(s-t_k)^\gamma\|\chi^\epsilon(s)-\chi(s)\|_{U_0}\|\eta(t_k)\|_{L^2(\mathcal{O})}\,ds\right)\\
\nonumber\leq&\frac{2L_1\sqrt{R}}{2^{\gamma N}}E\left(I_{G_{R,\epsilon}(T)}\sum_{k=1}^{2^N}\int_{t_{k-1}}^{t_k}(1+\|\theta(s)\|_{L^2(\mathcal{O})})^2\,ds\right)^{\frac{1}{2}}E\left(I_{G_{R,\epsilon}(T)}\sum_{k=1}^{2^N}\int_{t_{k-1}}^{t_k}\|\chi^\epsilon(s)-\chi(s)\|_{U_0}^2\,ds\right)^{\frac{1}{2}}\\
\leq&\frac{4L_1\sqrt{2TM(1+R)R}}{2^{\gamma N}},
\end{align}

\begin{align}\label{4.2.36}
\nonumber I_{\epsilon,N}^3=&E\left(I_{G_{R,\epsilon}(T)}\sum_{k=1}^{2^N}\int_{t_{k-1}}^{t_k}\left|\int_{\mathcal{O}}(\sigma(t_k,\theta(s))-\sigma(t_k.\theta(t_k))(\chi^\epsilon(s)-\chi(s))\eta(t_k)\,dxdydz\right|ds\right)\\
\nonumber\leq&E\left(I_{G_{R,\epsilon}(T)}\sum_{k=1}^{2^N}\int_{t_{k-1}}^{t_k}\|\sigma(t_k,\theta(s))-\sigma(t_k,\theta(t_k))\|_{\mathcal{L}_2(U_0,H_2)}\|\chi^\epsilon(s)-\chi(s)\|_{U_0}\|\eta(t_k)\|_{L^2(\mathcal{O})}\,ds\right)
\end{align}
\begin{align}
\nonumber\leq&2\sqrt{LR}E\left(I_{G_{R,\epsilon}(T)}\sum_{k=1}^{2^N}\int_{t_{k-1}}^{t_k}\|\theta(s)-\theta(t_k)\|\|_{L^2(\mathcal{O})}\|\chi^\epsilon(s)-\chi(s)\|_{U_0}\,ds\right)\\
\nonumber\leq&2\sqrt{LR}E\left(I_{G_{R,\epsilon}(T)}\sum_{k=1}^{2^N}\int_{t_{k-1}}^{t_k}\|\theta(s)-\theta(t_k)\|_{L^2(\mathcal{O})}^2\,ds\right)^{\frac{1}{2}}E\left(I_{G_{R,\epsilon}(T)}\sum_{k=1}^{2^N}\int_{t_{k-1}}^{t_k}\|\chi^\epsilon(s)-\chi(s)\|_{U_0}^2\,ds\right)^{\frac{1}{2}}\\
\leq&\frac{4\sqrt{\mathcal{D}LMR}}{2^{\frac{N}{4}}}
\end{align}

\begin{align}\label{4.2.37}
\nonumber I_{\epsilon,N}^4=&E\left(I_{G_{R,\epsilon}(T)}\sup_{1\leq k\leq 2^N}\sup_{t_{k-1}\leq t\leq t_k}\left|\int_{\mathcal{O}}\sigma(t_k,\theta(t_k))\int_{t_{k-1}}^t(\chi^\epsilon(s)-\chi(s))\,ds\,\eta(t_k)\,dxdydz\right|\right)\\
\nonumber \leq&\sup_{1\leq k\leq 2^N}\sup_{t_{k-1}\leq t\leq t_k}\left(E\left(I_{G_{R,\epsilon}(T)}\|\sigma(t_k,\theta(t_k))\|_{\mathcal{L}_2(U_0,H_2)}\left\|\int_{t_{k-1}}^t(\chi^\epsilon(s)-\chi(s))\,ds\right\|_{U_0}\|\eta(t_k)\|_{L^2(\mathcal{O})}\right)\right)\\
\nonumber\leq&2\sqrt{KR(1+R)}\sup_{1\leq k\leq 2^N}E\left(I_{G_{R,\epsilon}(T)}\int_{t_{k-1}}^{t_k}\|\chi^\epsilon(s)-\chi(s)\|_{U_0}\,ds\right)\\
\nonumber\leq&\frac{2\sqrt{KR(1+R)}}{2^{\frac{N}{2}}}E\left(I_{G_{R,\epsilon}(T)}\int_0^T\|\chi^\epsilon(s)-\chi(s)\|_{U_0}^2\,ds\right)^\frac{1}{2}\\
\leq&\frac{4\sqrt{KMR(1+R)}}{2^{\frac{N}{2}}}.
\end{align}
From the weak convergence of $\chi^\epsilon$ to $\chi,$ we know that for any $a,$ $b\in[0,T]$ with $a<b,$ the integral $\int_a^b\chi^\epsilon(s)\,ds\rightarrow\int_a^b\chi(s)\,ds$ weakly in $U_0$ as $\epsilon\rightarrow 0.$ Therefore, for any $\theta\in H_2,$ the operator $\sigma(\theta)$ is compact from $U_0$ to $H_2,$ which entails that $\sigma(\theta)\int_a^b\chi^\epsilon(s)\,ds\rightarrow\sigma(\theta)\int_a^b\chi(s)\,ds$ strongly in $H_2$ as $\epsilon\rightarrow0.$ Hence, a.s. for fixed $N,$ 
\begin{align*}
I_{G_{R,\epsilon}(T)}\sum_{k=1}^{2^N}\left|\int_{\mathcal{O}}\sigma(t_k,\theta(t_k))\int_{t_{k-1}}^{t_k}(\chi^\epsilon(s)-\chi(s))\,ds\,\eta(t_k)\,dxdydz\right|\rightarrow 0
\end{align*}
as $\epsilon\rightarrow 0$ and
\begin{align*}
I_{G_{R,\epsilon}(T)}\sum_{k=1}^{2^N}\left|\int_{\mathcal{O}}\sigma(t_k,\theta(t_k))\int_{t_{k-1}}^{t_k}(\chi^\epsilon(s)-\chi(s))\,ds\,\eta(t_k)\,dxdydz\right|\leq4\sqrt{KMRT(1+R)}.
\end{align*}
It follows from the dominated convergence theorem that for any fixed $N,$ 
\begin{align}\label{4.2.38}
\lim_{\epsilon\rightarrow0}I_{\epsilon,N}^5=0.
\end{align}

Define $\tau_{\tilde{R}}^\epsilon=\{t>0:\|\theta^\epsilon(t)\|_{L^2(\mathcal{O})}^2+\int_0^t\|\theta^\epsilon(s)\|^2\,ds>{\tilde{R}}\},$ applying It\^{o} formula to $\|\theta^\epsilon(t)\|_{L^2(\mathcal{O})}^2$  and combining assumption ($A_3$) with Young's inequality, H\"{o}lder's inequality, yield
\begin{align}\label{4.2.39}
\nonumber&\|\theta^\epsilon(t)\|_{L^2(\mathcal{O})}^2+2\int_0^t\|\theta^\epsilon(s)\|^2\,ds\\
\nonumber=&2\beta\int_0^t\int_{\Gamma_u}\theta^*\theta^\epsilon(s)\,dxdyds+2\int_0^t\int_{\mathcal{O}}(g(x,y,z,s)+\sigma(s,\theta^\epsilon(s))\chi^\epsilon(s))\theta^\epsilon(x,y,z,s)\,dxdydzds+\|\theta_0\|_{L^2(\mathcal{O})}^2\\
\nonumber&+2\sqrt{\epsilon}\int_0^t\int_{\mathcal{O}}\theta^\epsilon(x,y,z,s)\sigma(s,\theta^\epsilon(s))dW(s)\,dxdydz+\epsilon\int_0^t\|\sigma(s,\theta^\epsilon(s))\|_{\mathcal{L}_2(U_0,H_2)}^2\,ds\\
\nonumber\leq&2\beta\int_0^t\|\theta^*\|_{L^2(\mathcal{M})}\|\theta^\epsilon(s)\|_{L^2(\Gamma_u)}\,ds+2\int_0^t\|g(s)\|_{L^2(\mathcal{O})}\|\theta^\epsilon(s)\|_{L^2(\mathcal{O})}\,ds\\
\nonumber&+2\int_0^t\|\sigma(s,\theta^\epsilon(s))\|_{\mathcal{L}_2(U_0,H_2)}\|\chi^\epsilon(s))\|_{U_0}\|\theta^\epsilon(s)\|_{L^2(\mathcal{O})}\,ds\\
\nonumber&+2\sqrt{\epsilon}\int_0^t\int_{\mathcal{O}}\theta^\epsilon(x,y,z,s)\sigma(t,\theta^\epsilon(s))dW(s)\,dxdydz+K\int_0^t(1+\|\theta^\epsilon(s)\|_{L^2(\mathcal{O})}^2)\,ds+\|\theta_0\|_{L^2(\mathcal{O})}^2\\
\nonumber\leq&C\int_0^t(\|\theta^*\|_{L^2(\mathcal{M})}^2+\|g(s)\|_{L^2(\mathcal{O})}^2+\|\chi^\epsilon(s))\|_{U_0}^2+1)\,ds+C\int_0^t(1+\|\chi^\epsilon(s))\|_{U_0}^2)\|\theta^\epsilon(s)\|_{L^2(\mathcal{O})}^2\,ds\\
&+\int_0^t\|\theta^\epsilon(s)\|^2\,ds+2\sqrt{\epsilon}\int_0^t\int_{\mathcal{O}}\theta^\epsilon(x,y,z,s)\sigma(s,\theta^\epsilon(s))dW(s)\,dxdydz+\|\theta_0\|_{L^2(\mathcal{O})}^2.
\end{align}
Taking the supremum up to time $T\land\tau^\epsilon_{\tilde{R}}$ in inequality \eqref{4.2.39}, we obtain
\begin{align}\label{4.2.40}
\nonumber &\sup_{0\leq t\leq T\land\tau^\epsilon_{\tilde{R}}}\left(\|\theta^\epsilon(t)\|_{L^2(\mathcal{O})}^2\right)+\int_0^{T\land\tau^\epsilon_{\tilde{R}}}\|\theta^\epsilon(s)\|^2\,ds\\
\nonumber\leq&\|\theta_0\|_{L^2(\mathcal{O})}^2+C\int_0^T\|g(s)\|_{L^2(\mathcal{O})}^2\,ds+2\sqrt{\epsilon}\sup_{0\leq r\leq T\land\tau^\epsilon_{\tilde{R}}}\left|\int_0^r\int_{\mathcal{O}}\theta^\epsilon_n(s)\sigma(s,\theta^\epsilon(s))\,dW(s)dxdydz\right|\\
&+C\int_0^{T\land\tau^\epsilon_{\tilde{R}}}(1+\|\chi^\epsilon(s))\|_{U_0}^2)\sup_{0\leq r\leq s}\|\theta^\epsilon(r)\|_{L^2(\mathcal{O})}^2\,ds+CT(1+\|\theta^*\|_{L^2(\mathcal{M})}^2)+CM.
\end{align}
We deduce from Lemma \ref{2.2.4} that
\begin{align}\label{4.2.41}
\nonumber &\sup_{0\leq t\leq T\land\tau^\epsilon_{\tilde{R}}}\left(\|\theta^\epsilon(t)\|_{L^2(\mathcal{O})}^2\right)+\int_0^{T\land\tau^\epsilon_{\tilde{R}}}\|\theta^\epsilon(s)\|^2\,ds\\
\nonumber\leq&\left(\|\theta_0\|_{L^2(\mathcal{O})}^2+C\int_0^T\|g(s)\|_{L^2(\mathcal{O})}^2\,ds+2\sqrt{\epsilon}\sup_{0\leq r\leq T\land\tau^\epsilon_{\tilde{R}}}\left|\int_0^r\int_{\mathcal{O}}\theta^\epsilon_n(s)\sigma(s,\theta^\epsilon(s))\,dW(s)dxdydz\right|\right.\\
\nonumber&\left.+CT(1+\|\theta^*\|_{L^2(\mathcal{M})}^2)+CM\right)e^{C\int_0^{T\land\tau^\epsilon_{\tilde{R}}}(1+\|\chi^\epsilon(s))\|_{U_0}^2)\,ds}\\
\nonumber\leq&\left(\|\theta_0\|_{L^2(\mathcal{O})}^2+C\int_0^T\|g(s)\|_{L^2(\mathcal{O})}^2\,ds+2\sqrt{\epsilon}\sup_{0\leq r\leq T\land\tau^\epsilon_{\tilde{R}}}\left|\int_0^r\int_{\mathcal{O}}\theta^\epsilon_n(s)\sigma(s,\theta^\epsilon(s))\,dW(s)dxdydz\right|\right.\\
&\left.+CT(1+\|\theta^*\|_{L^2(\mathcal{M})}^2)+CM\right)e^{C(T+M)}.
\end{align}
From the Burkholder-Davis-Gundy inequality, we conclude 
\begin{align}\label{4.2.42}
\nonumber&2\sqrt{\epsilon}E\left(\sup_{0\leq r\leq T\land\tau^\epsilon_{\tilde{R}}}\left|\int_0^r\int_{\mathcal{O}}\theta^\epsilon(s)\sigma(s,\theta^\epsilon(s))\,dW(s)dxdydz\right|\right)\\
\nonumber\leq &C\sqrt{\epsilon}E\left(\left(\int_0^{T\land\tau^\epsilon_{\tilde{R}}}\|\theta^\epsilon(s)\|_{L^2(\mathcal{O})}^2\|\sigma(s,\theta^\epsilon(s))\|_{\mathcal{L}_2(U_0,H_2)}^2\,ds\right)^{\frac{1}{2}}\right)\\
\nonumber\leq &C\sqrt{\epsilon}E\left(\sup_{0\leq r\leq T\land\tau^\epsilon_{\tilde{R}}}\|\theta^\epsilon(r)\|_{L^2(\mathcal{O})}\left(\int_0^{T\land\tau^\epsilon_{\tilde{R}}}\|\sigma(s,\theta^\epsilon(s))\|_{\mathcal{L}_2(U_0,H_2)}^2\,ds\right)^{\frac{1}{2}}\right)\\
\nonumber\leq &\frac{1}{2}E\left(\sup_{0\leq r\leq T\land\tau^\epsilon_{\tilde{R}}}\|\theta^\epsilon(r)\|_{L^2(\mathcal{O})}^2\right)+C\epsilon E\left(\int_0^{T\land\tau^\epsilon_{\tilde{R}}}\|\sigma(s,\theta^\epsilon(s))\|_{\mathcal{L}_2(U_0,H_2)}^2\,ds\right)\\
\leq &\frac{1}{2e^{C(T+M)}}E\left(\sup_{0\leq r\leq T\land\tau^\epsilon_{\tilde{R}}}\|\theta^\epsilon(r)\|_{L^2(\mathcal{O})}^2\right)+CK\epsilon E\left(\int_0^{T\land\tau^\epsilon_{\tilde{R}}}\sup_{0\leq r\leq s}\|\theta^\epsilon(r))\|_{L^2(\mathcal{O})}^2\,ds\right)+CKT.
\end{align}
Taking the expectation on both hand sides of inequality \eqref{4.2.41}, it follows from inequality \eqref{4.2.42} and Lemma \ref{2.2.4} that
\begin{align}\label{4.2.43}
\nonumber &E\left(\sup_{0\leq t\leq T\land\tau^\epsilon_{\tilde{R}}}\left(\|\theta^\epsilon(t)\|_{L^2(\mathcal{O})}^2\right)\right)+E\left(\int_0^{T\land\tau^\epsilon_{\tilde{R}}}\|\theta^\epsilon(s)\|^2\,ds\right)\\
\leq&\left(\|\theta_0\|_{L^2(\mathcal{O})}^2+C\int_0^T\|g(s)\|_{L^2(\mathcal{O})}^2\,ds+CT(1+\|\theta^*\|_{L^2(\mathcal{M})}^2)+CM\right)e^{C(T+M+KTe^{C(T+M)})}.
\end{align}
Let $\tilde{R}\rightarrow+\infty$ in both side of inequality \eqref{4.2.43}, we obtain
\begin{align}\label{4.2.44}
\nonumber &E\left(\sup_{0\leq t\leq T}\left(\|\theta^\epsilon(t)\|_{L^2(\mathcal{O})}^2\right)+\int_0^T\|\theta^\epsilon(s)\|^2\,ds\right)\\
\leq&\left(\|\theta_0\|_{L^2(\mathcal{O})}^2+C\int_0^T\|g(s)\|_{L^2(\mathcal{O})}^2\,ds+CT(1+\|\theta^*\|_{L^2(\mathcal{M})}^2)+CM\right)e^{C(T+M+KTe^{C(T+M)})}.
\end{align}
Inequality \eqref{4.2.16} and inequality \eqref{4.2.41} imply that 
\begin{align*}
\mathcal{P}\left(G_{\epsilon, R}(T)^c\right)\leq &\mathcal{P}\left(\left\{\omega:\left(\sup_{s\in[0,t]}\|\theta(s,\omega)\|_{L^2(\mathcal{O})^2}+\int_0^t\|\theta(s,\omega)\|^2\,ds\right)>R\right\}\right)\\
&+\mathcal{P}\left(\left\{\omega:\left(\sup_{s\in[0,t]}\|\theta^\epsilon(s,\omega)\|_{L^2(\mathcal{O})^2}+\int_0^t\|\theta^\epsilon(s,\omega)\|^2\,ds\right)>R\right\}\right)\\
\leq&\frac{1}{R}\left(\|\theta_0\|_{L^2(\mathcal{O})}^2+C\int_0^T\|g(s)\|_{L^2(\mathcal{O})}^2\,ds+CT(1+\|\theta^*\|_{L^2(\mathcal{M})}^2)+CM\right)e^{C(T+M+KTe^{C(T+M)})}\\
&+\frac{1}{R}C\left(\int_0^T\|g(s)\|_{L^2(\mathcal{O})}^2\,ds+T\|\theta^*\|_{L^2(\mathcal{M})}^2+\|\theta_0\|_{L^2(\mathcal{O})}^2+KM\right)e^{CKM},
\end{align*}
which entails that for any $\delta>0,$ there exists a positive constant $R_0$ such that for any $R\geq R_0$ and any $\epsilon\in(0,1],$
\begin{align}\label{4.2.45}
\mathcal{P}\left(G_{\epsilon, R}(T)^c\right)\leq\frac{\delta}{2}.
\end{align}
For any given $\alpha>0,$ we may choose some integer $N_0\geq 1$ large enough such that for any $R\geq R_0$ and any $N\geq N_0,$
\begin{align}\label{4.2.46}
2e^{C(T+R+LM)}\sum_{j=1}^4 I_{\epsilon,N}^j\leq \frac{\alpha\delta}{6}.
\end{align}
Then for any fixed $R\geq R_0$ and any fixed $N\geq N_0$ there exists $\epsilon_0>0$ such that for any $\epsilon\in(0,\epsilon_0],$
\begin{align}\label{4.2.47}
2e^{C(T+R+LM)}I_{\epsilon,N}^5\leq \frac{\alpha\delta}{6}
\end{align}
and
\begin{align}\label{4.2.48}
\left(2K\epsilon T(1+R)+\sqrt{\epsilon KRT(1+R)}\right)e^{C(T+R+LM)}\leq \frac{\alpha\delta}{6}.
\end{align}
We infer from inequalities \eqref{4.2.46}-\eqref{4.2.48} that for any $\epsilon\in(0,\epsilon_0],$
\begin{align}\label{4.2.49}
E\left(I_{G_{R,\epsilon}(T)}\sup_{s\in [0,T]}\|\eta(t)\|_{L^2(\mathcal{O})}^2\right)+E\left(I_{G_{R,\epsilon}(T)}\int_0^T\|\eta(s)\|^2\,ds\right)\leq\frac{\alpha\delta}{2}.
\end{align}
For any $\alpha>0,$ it follows from Markov inequality and inequalities \eqref{4.2.45}, \eqref{4.2.49} that for any $\epsilon\in(0,\epsilon_0],$
\begin{align*}
&\mathcal{P}\left(\sup_{s\in [0,T]}\|\eta(t)\|_{L^2(\mathcal{O})}^2+\int_0^T\|\eta(s)\|^2\,ds\geq\delta\right)\\
\leq &\frac{1}{\delta}E\left(I_{G_{R,\epsilon}(T)}\sup_{s\in [0,T]}\|\eta(t)\|_{L^2(\mathcal{O})}^2+I_{G_{R,\epsilon}(T)}\int_0^T\|\eta(s)\|^2\,ds\right)+\mathcal{P}\left(G_{\epsilon, R}(T)^c\right)\\
\leq&\alpha.
\end{align*}
\qed\hfill

With the above results in hands, we immediately obtain the following large deviation theorem.
\begin{theorem}
Let $\{(v^\epsilon,p_s^\epsilon,\theta^\epsilon)\}$ be the unique weak solution of problem
\begin{equation}\label{4.2.44}
\begin{cases}
&\nabla p^\epsilon_{s}+f{v^\epsilon}^\bot+L_1v^\epsilon=\int_{-h}^z\nabla \theta^\epsilon(x,y,\zeta,t)\,d\zeta,\\
&\int_{-h}^{0}\nabla\cdot v^\epsilon(x,y,\zeta,t)\,d\zeta=0,\\
&d\theta^\epsilon+B(v^\epsilon,\theta^\epsilon)\,dt+A_2(\theta^\epsilon-\theta^*)\,dt-K_h\Delta \theta^*\,dt=g(x,y,z,t)\,dt+\sqrt{\epsilon}\sigma(t,\theta^\epsilon(t))dW(t),\\
&\left.A_{\nu}\frac{\partial v^\epsilon}{\partial z}\right|_{\Gamma_u}=\mu, \left.\frac{\partial v^\epsilon}{\partial z}\right|_{\Gamma_b}=0, \left.v^\epsilon\cdot\vec{n}\right|_{\Gamma_l}=0,\left.\frac{\partial v^\epsilon}{\partial\vec{n}}\times\vec{n}\right|_{\Gamma_l}=0,\\
&\theta^\epsilon(x,y,z,0)=\theta_0(x,y,z)\in H_2.
\end{cases}
\end{equation}
Then $\{\theta^\epsilon\}$ satisfies the Laplace principle in $\mathcal{C}([0,T];H_2)\cap L^2(0,T;V_2)$ with a good rate function
\begin{align}
I(\theta)=\inf_{\chi\in L^2(0,T;U_0): \theta=\Phi^0(\int_0^\cdot \chi(s)\,ds)}\left\{\frac{1}{2}\int_0^T\|\chi(s)\|_{U_0}^2\,ds\right\}
\end{align}
with the convention that the infimum of an empty set is infinity.
\end{theorem}
{\bf Proof.} From Lemma \ref{2.3.1}, Theorem \ref{a} and Theorem \ref{4.2.20}, we conclude that $\{\theta^\epsilon\}$ satisfies the Laplace principle which is equivalent to the large deviation principle in $\mathcal{C}([0,T];H_2)\cap L^2(0,T;V_2)$ with the above-mentioned rate function.\\
\qed\hfill

\section{Appendix}
\def\theequation{A.\arabic{equation}}\makeatother
\setcounter{equation}{0}
\subsection*{A.1. Proof of Lemma \ref{2.2.3}}
For any $v\in H^2(\mathcal{O})\cap V_1$ and $\theta,$ $\eta\in V_2,$ we infer from H\"{o}lder's inequality and interpolation inequality that
\begin{align*}
|b(v,\theta,\eta)|\leq&\left|\int_{\mathcal{O}}\left(v(x,y,z)\cdot\nabla\theta(x,y,z)-\left(\int_{-h}^z\nabla\cdot v(x,y,\zeta)\,d\zeta\right)\frac{\partial \theta(x,y,z)}{\partial z}\right)\eta(x,y,z)\,dxdydz\right|\\
\leq&\|v\|_{L^6(\mathcal{O})}\|\nabla\theta\|_{L^2(\mathcal{O})}\|\eta\|_{L^3(\mathcal{O})}+\int_M\left(\int_{-h}^0|\nabla v(x,y,\zeta)|\,d\zeta\right)\left(\int_{-h}^0|\partial_z \theta(x,y,z)||\eta(x,y,z)|\,dz\right)dxdy\\
\leq&\|v\|_{L^6(\mathcal{O})}\|\nabla\theta\|_{L^2(\mathcal{O})}\|\eta\|_{L^3(\mathcal{O})}+\sqrt{h}\int_M\|\nabla v(x,y)\|_{L^2(-h,0)}\|\partial_z \theta(x,y)\|_{L^2(-h,0)}\|\eta(x,y)\|_{L^2(-h,0)}dxdy\\
\leq&\|v\|_{L^6(\mathcal{O})}\|\nabla\theta\|_{L^2(\mathcal{O})}\|\eta\|_{L^3(\mathcal{O})}+\sqrt{h}\|\nabla v\|_{L^4(M;L^2(-h,0))}\|\partial_z \theta\|_{L^2(M;L^2(-h,0))}\|\eta\|_{L^4(M;L^2(-h,0))}\\
\leq&\|v\|_{L^6(\mathcal{O})}\|\nabla\theta\|_{L^2(\mathcal{O})}\|\eta\|_{L^3(\mathcal{O})}+\sqrt{h}\|\nabla v\|_{L^2(-h,0;L^4(M))}\|\partial_z \theta\|_{L^2(\mathcal{O})}\|\eta\|_{L^2(-h,0;L^4(M))}\\
\leq&C\sqrt{h}\left\|\|\nabla v\|_{L^2(M)}^{\frac{1}{2}}\|\nabla v\|_{H^1(M)}^{\frac{1}{2}}\right\|_{L^2(-h,0)}\|\partial_z \theta\|_{L^2(\mathcal{O})}\left\|\|\eta\|_{L^2(M)}^{\frac{1}{2}}\|\nabla\eta\|_{L^2(M)}^{\frac{1}{2}}\right\|_{L^2(-h,0)}\\
&+\|v\|_{L^6(\mathcal{O})}\|\nabla\theta\|_{L^2(\mathcal{O})}\|\eta\|_{L^3(\mathcal{O})}\\
\leq&\|v\|_{H^1(\mathcal{O})}\|\theta\|\|\eta\|_{L^2(\mathcal{O})}^{\frac{1}{2}}\|\eta\|^{\frac{1}{2}}+C\sqrt{h}\|\nabla v\|_{L^2(\mathcal{O})}^{\frac{1}{2}}\|\nabla v\|_{H^1(\mathcal{O})}^{\frac{1}{2}}\|\partial_z \theta\|_{L^2(\mathcal{O})}\|\eta\|_{L^2(\mathcal{O})}^{\frac{1}{2}}\|\nabla\eta\|_{L^2(\mathcal{O})}^{\frac{1}{2}}\\
\leq&C\|v\|_{H^1(\mathcal{O})}^{\frac{1}{2}}\|v\|_{H^2(\mathcal{O})}^{\frac{1}{2}}\|\theta\|\|\eta\|_{L^2(\mathcal{O})}^{\frac{1}{2}}\|\eta\|^{\frac{1}{2}}.
\end{align*}
\qed\hfill

\subsection*{A.2. Proof of Lemma \ref{2.2.4}}

Let
\begin{align*}
\tilde{Y}(t)=\int_0^ta(s)Y(s)\,ds,
\end{align*}
then we obtain
\begin{align*}
\frac{d}{dt}\tilde{Y}(t)=&a(t)Y(t)\\
\leq &a(t)\tilde{Y}(t)+a(t)Z(t).
\end{align*}
We deduce from the classical Gronwall inequality and the fact that $Z(t)$ is non-negative, non-decreasing function  that
\begin{align*}
\tilde{Y}(t)\leq &\tilde{Y}(0)e^{\int_0^ta(r)\,dr}+\int_0^ta(s)Z(s)e^{\int_s^ta(r)\,dr}\,ds\\
\leq&Z(t)\int_0^ta(s)e^{\int_s^ta(r)\,dr}\,ds\\
=&Z(t)(e^{\int_0^ta(r)\,dr}-1),
\end{align*}
which entails that
\begin{align*}
Y(t)+\int_0^tX(s)\,ds\leq&\tilde{Y}(t)+Z(t)\\
\leq&Z(t)e^{\int_0^ta(r)\,dr}.
\end{align*}
\qed\hfill

Finally, we will establish a technical lemma used to prove the large deviation principle which studies time increments of the solution to the stochastic control equation. For $N\geq 1$ and $k=0,1,\cdots, 2^N,$ let $t_k=kT2^{-N}.$ Given $M>0,$ $\chi\in\mathcal{A}_M,$ $\epsilon\geq 0$ is small enough, let $(v_\chi^\epsilon,p_{\chi s}^\epsilon,\theta_\chi^\epsilon)$ be the weak solution of the following problem:
\begin{equation}\label{5.1}
\begin{cases}
&\nabla p^\epsilon_{\chi s}+f{v_\chi^\epsilon}^\bot+L_1v_\chi^\epsilon=\int_{-h}^z\nabla \theta_\chi^\epsilon(x,y,\zeta,t)\,d\zeta,\\
&\int_{-h}^{0}\nabla\cdot v_\chi^\epsilon(x,y,\zeta,t)\,d\zeta=0,\\
&d\theta_\chi^\epsilon+B(v_\chi^\epsilon,\theta_\chi^\epsilon)\,dt+A_2(\theta_\chi^\epsilon-\theta^*)\,dt-K_h\Delta \theta^*\,dt=g(x,y,z,t)\,dt+\sigma(t,\theta_\chi^\epsilon(t))\chi(t)\,dt+\sqrt{\epsilon}\sigma(t,\theta_\chi^\epsilon(t))dW(t),\\
&\left.A_{\nu}\frac{\partial v_\chi^\epsilon}{\partial z}\right|_{\Gamma_u}=\mu, \left.\frac{\partial v_\chi^\epsilon}{\partial z}\right|_{\Gamma_b}=0, \left.v_\chi^\epsilon\cdot\vec{n}\right|_{\Gamma_l}=0,\left.\frac{\partial v_\chi^\epsilon}{\partial\vec{n}}\times\vec{n}\right|_{\Gamma_l}=0,\\
&\theta_\chi^\epsilon(x,y,z,0)=\theta_0(x,y,z)
\end{cases}
\end{equation}
and for any $t\in[0,T],$ define
\begin{align*}
G_R(t)=\left\{\omega:\sup_{r\in[0,t]}\|\theta_\chi^\epsilon(s,\omega)\|_{L^2(\mathcal{O})}^2+\int_0^t\|\theta_\chi^\epsilon(s,\omega)\|^2\leq R\right\}.
\end{align*}
Then we have the following conclusion:
\begin{theorem}\label{5.2}
Assume that ($A_1$)-($A_3$) hold and let $M,$ $R>0.$ Let $\theta_0\in H_2$ and $(v_\chi^\epsilon,p_{\chi s}^\epsilon,\theta_\chi^\epsilon)$ is the weak solution of problem \eqref{5.1}. Then there exists a positive constant $\mathcal{D}=\mathcal{D}(\epsilon,K,M,\|\mu\|_{H^1(\mathcal{M})},\|g\|_{L^2(0,T;L^2(\mathcal{O}))},R,T)$ such that for any $\chi\in\mathcal{A}_M,$ $\epsilon\in[0,\epsilon_0],$ 
\begin{align*}
\sum_{k=1}^{2^N}E\left(I_{G_R(T)}\int_{t_{k-1}}^{t_k}\|\theta_\chi^\epsilon(s)-\theta^\epsilon_\chi(t_k)\|_{L^2(\mathcal{O})}^2\,ds\right)\leq \mathcal{D}2^{-\frac{N}{2}}.
\end{align*}
\end{theorem}
{\bf Proof.} Let $\chi\in\mathcal{A}_M,$ we deduce from the It\^{o} formula that for any $k\in[1,2^N]$ and $t_{k-1}\leq s\leq t_k$ for any $k\in[1,2^N],$
\begin{align*}
&\|\theta_\chi^\epsilon(s)-\theta^\epsilon_\chi(t_k)\|_{L^2(\mathcal{O})}^2
=2\int_s^{t_k}\int_{\mathcal{O}}g(x,y,z,r)(\theta_\chi^\epsilon(r)-\theta^\epsilon_\chi(s))\,dxdydzdr\\
&+2\int_s^{t_k}\int_{\mathcal{O}}\sigma(r,\theta_\chi^\epsilon(r))\chi(r)(\theta_\chi^\epsilon(r)-\theta^\epsilon_\chi(s))\,dxdydzdr\\
&+2\sqrt{\epsilon}\int_s^{t_k}\int_{\mathcal{O}}\sigma(r,\theta_\chi^\epsilon(r))dW(r)(\theta_\chi^\epsilon(r)-\theta^\epsilon_\chi(s))\,dxdydz
-2\int_s^{t_k}\int_{\mathcal{O}}A_2(\theta_\chi^\epsilon(r)-\theta^*)(\theta_\chi^\epsilon(r)-\theta^\epsilon_\chi(s))\,dxdydzdr\\
&-2\int_s^{t_k}b(v_\chi^\epsilon(r),\theta_\chi^\epsilon(r),\theta_\chi^\epsilon(r)-\theta^\epsilon_\chi(s))\,dr+2K_h\int_s^{t_k}\int_{\mathcal{O}}\Delta\theta^*(\theta_\chi^\epsilon(r)-\theta^\epsilon_\chi(s))\,dxdydzdr+\epsilon\int_s^{t_k}\|\sigma(r,\theta_\chi^\epsilon(r))\|_{\mathcal{L}_2(U_0,H_2)}^2\,dr.
\end{align*}

In what follows, we will estimate the each term of the right hand side of the above inequality step by step. Clearly, $G_R(T)\subset G_R(r)$ for any $r\in[0,T],$ which implies that $\|\theta_\chi^\epsilon(r)\|_{L^2(\mathcal{O})}+\|\theta_\chi^\epsilon(s)\|_{L^2(\mathcal{O})}\leq 2\sqrt{R} $ for any $0\leq s\leq r\leq T.$ To begin with, it follows from H\"{o}lder's inequality that
\begin{align}\label{5.3}
\nonumber&2\sum_{k=1}^{2^N}E\left(I_{G_R(T)}\int_{t_{k-1}}^{t_k}\left|\int_s^{t_k}\int_{\mathcal{O}}g(x,y,z,r)(\theta_\chi^\epsilon(r)-\theta^\epsilon_\chi(s))\,dxdydzdr\right|\,ds\right)\\
\nonumber\leq&2\sum_{k=1}^{2^N}E\left(\int_{t_{k-1}}^{t_k}\int_s^{t_k}I_{G_R(r)}\|g(r)\|_{L^2(\mathcal{O})}\|\theta_\chi^\epsilon(r)-\theta^\epsilon_\chi(s)\|_{L^2(\mathcal{O})}\,dr\,ds\right)\\
\nonumber\leq&4\sqrt{R}\sum_{k=1}^{2^N}\int_{t_{k-1}}^{t_k}\int_s^{t_k}\|g(r)\|_{L^2(\mathcal{O})}\,dr\,ds\\
\leq&\frac{4T\sqrt{RT}}{2^N}\|g\|_{L^2(0,T;L^2(\mathcal{O}))}.
\end{align}
We infer from  H\"{o}lder's inequality and assumption ($A_2$) that
\begin{align}\label{5.4}
\nonumber&2\sum_{k=1}^{2^N}E\left(I_{G_R(T)}\int_{t_{k-1}}^{t_k}\left|\int_s^{t_k}\int_{\mathcal{O}}\sigma(r,\theta_\chi^\epsilon(r))\chi(r)(\theta_\chi^\epsilon(r)-\theta^\epsilon_\chi(s))\,dxdydzdr\right|\,ds\right)\\
\nonumber\leq&2\sum_{k=1}^{2^N}E\left(\int_{t_{k-1}}^{t_k}\int_s^{t_k}I_{G_R(r)}\|\sigma(r,\theta_\chi^\epsilon(r))\|_{\mathcal{L}_2(U_0,H_2)}\|\chi(r)\|_{U_0}\|\theta_\chi^\epsilon(r)-\theta^\epsilon_\chi(s)\|_{L^2(\mathcal{O})}\,dr\,ds\right)\\
\nonumber\leq&4\sqrt{KR(1+R)}\sum_{k=1}^{2^N}\int_{t_{k-1}}^{t_k}\int_s^{t_k}\|\chi(r)\|_{U_0}\,dr\,ds\\
\leq&\frac{4T\sqrt{KM(1+R)RT}}{2^N}.
\end{align}
From the Burkholder-Davis-Gundy inequality and assumption ($A_2$), we conclude 
\begin{align}\label{5.5}
\nonumber&2\sqrt{\epsilon}\sum_{k=1}^{2^N}E\left(I_{G_R(T)}\int_{t_{k-1}}^{t_k}\left|\int_s^{t_k}\int_{\mathcal{O}}\sigma(r,\theta_\chi^\epsilon(r))dW(r)(\theta_\chi^\epsilon(r)-\theta^\epsilon_\chi(s))\,dxdydz\right|\,ds\right)\\
\nonumber\leq&2\sqrt{\epsilon}\sum_{k=1}^{2^N}\int_{t_{k-1}}^{t_k}E\left(\int_s^{t_k}I_{G_R(r)}\|\sigma(r,\theta_\chi^\epsilon(r))\|_{\mathcal{L}_2(U_0,H_2)}^2\|\theta_\chi^\epsilon(r)-\theta^\epsilon_\chi(s)\|_{L^2(\mathcal{O})}^2\,dr\right)^{\frac{1}{2}}\,ds\\
\leq&\frac{4T\sqrt{K(1+R)RT\epsilon}}{2^\frac{N}{2}}.
\end{align}
Applying H\"{o}lder's inequality and Young's inequality, yields 
\begin{align}\label{5.6}
\nonumber&-2\sum_{k=1}^{2^N}E\left(I_{G_R(T)}\int_{t_{k-1}}^{t_k}\int_s^{t_k}\int_{\mathcal{O}}A_2(\theta_\chi^\epsilon(r)-\theta^*)(\theta_\chi^\epsilon(r)-\theta^\epsilon_\chi(s))\,dxdydzdr\right.\\
\nonumber&\left.-K_h\int_s^{t_k}\int_{\mathcal{O}}\Delta\theta^*(\theta_\chi^\epsilon(r)-\theta^\epsilon_\chi(s))\,dxdydzdr\,ds\right)\\
\nonumber\leq
&2\sum_{k=1}^{2^N}\int_{t_{k-1}}^{t_k}E\left(I_{G_R(T)}\int_s^{t_k}(-\|\theta_\chi^\epsilon(r)\|^2+\|\theta_\chi^\epsilon(r)\|\|\theta^\epsilon_\chi(s)\|+\beta\|\theta^*\|_{L^2(\mathcal{M})}\|\theta^\epsilon_\chi(s)\|_{L^2(\Gamma_u)}\right.\\
\nonumber&\left.+\beta\|\theta^*\|_{L^2(\mathcal{M})}\|\theta^\epsilon_\chi(r)\|_{L^2(\Gamma_u)}\,dr\right)\,ds\\
\leq&\frac{C}{2^N}.
\end{align}
It follows from assumption ($A_2$) that
\begin{align}\label{5.7}
\nonumber&\epsilon\sum_{k=1}^{2^N}E\left(I_{G_R(T)}\int_{t_{k-1}}^{t_k}\int_s^{t_k}\|\sigma(r,\theta_\chi^\epsilon(r))\|_{\mathcal{L}_2(U_0,H_2)}^2\,dr\,ds\right)\\
\nonumber\leq&\epsilon\sum_{k=1}^{2^N}E\left(\int_{t_{k-1}}^{t_k}\int_s^{t_k}I_{G_R(r)}\|\sigma(r,\theta_\chi^\epsilon(r))\|_{\mathcal{L}_2(U_0,H_2)}^2\,dr\,ds\right)\\
\leq&\frac{\epsilon K(1+R)T^2}{2^N}.
\end{align}
Finally, we conclude from Lemma \ref{2.2.3} and H\"{o}lder's inequality that
\begin{align}\label{5.8}
\nonumber&2\sum_{k=1}^{2^N}E\left(I_{G_R(T)}\int_{t_{k-1}}^{t_k}\left|\int_s^{t_k}b(v_\chi^\epsilon,\theta_\chi^\epsilon,\theta_\chi^\epsilon(r)-\theta^\epsilon_\chi(s))\,dr\right|\,ds\right)\\
\nonumber\leq&\mathcal{K}\sum_{k=1}^{2^N}E\left(I_{G_R(T)}\int_{t_{k-1}}^{t_k}\int_s^{t_k}\|v_\chi^\epsilon(r)\|_{H^1(\mathcal{O})}^{\frac{1}{2}}\|v_\chi^\epsilon(r)\|_{H^2(\mathcal{O})}^{\frac{1}{2}}\| \theta_\chi^\epsilon(s)\|\|\theta^\epsilon_\chi(r)\|_{L^2(\mathcal{O})}^{\frac{1}{2}}\|\theta^\epsilon_\chi(r)\|^{\frac{1}{2}}\,dr\right)\,ds\\
\nonumber\leq&\mathcal{K}\sum_{k=1}^{2^N}\left(E\left(I_{G_R(T)}\int_{t_{k-1}}^{t_k}\|v_\chi^\epsilon(r)\|_{H^1(\mathcal{O})}^{\frac{1}{2}}\|v_\chi^\epsilon(r)\|_{H^2(\mathcal{O})}^{\frac{1}{2}}\|\theta^\epsilon_\chi(r)\|_{L^2(\mathcal{O})}^{\frac{1}{2}}\|\theta^\epsilon_\chi(r)\|^{\frac{1}{2}}\,dr\int_{t_{k-1}}^{t_k}\| \theta_\chi^\epsilon(s)\|\,ds\right)\right)\\
\leq&\frac{\mathcal{K}\sqrt{R}(R+\|\mu\|_{H^1(M)}^2)}{2^\frac{N}{2}}.
\end{align}
Therefore, we deduce from inequalities \eqref{5.3}-\eqref{5.8} that the desired inequality holds.\\
\qed\hfill

\section*{Acknowledgement}
This work was supported by the National Science Foundation of China Grant (11401459,11871389), the Nat- ural Science Foundation of Shaanxi Province (2018JM1012) and the Fundamental Research Funds for the Central Universities (xjj2018088).

\bibliographystyle{abbrv}
\bibliography{BIB}

\end{document}